\documentclass[12pt,letterpaper,reqno]{amsart}
\usepackage{amsmath}
\usepackage{amssymb}
\usepackage{amsthm}
\usepackage{indentfirst}
\usepackage{xspace}
\usepackage{multirow}
\usepackage{hyperref}
\usepackage{xcolor}
\usepackage{calc}
\usepackage{float}
\usepackage{makecell}
\definecolor{darkblue}{rgb}{0,0,0.4}
\hypersetup{pdfauthor={Emily Dumas and Andrew
    Neitzke},pdftitle={Opers and nonabelian Hodge: numerical studies},colorlinks=true,urlcolor=darkblue,linkcolor=darkblue,citecolor=darkblue}

\definecolor{mpl1}{HTML}{1F77B4}
\definecolor{mpl2}{HTML}{FF7F0E}

\usepackage{verbatim}
\usepackage[letterpaper,margin=25mm,marginparwidth=23mm]{geometry}
\usepackage{thmtools}
\usepackage[all]{xy}
\usepackage{capt-of}
\usepackage{enumitem}
\usepackage{graphicx}
\usepackage{siunitx}

\usepackage{booktabs}

\graphicspath{{figures/}}

\usepackage{marginnote}

\newlength{\floatpagetop}

\newlength{\bibitemsep}
\newlength{\bibparskip}
\let\oldthebibliography\thebibliography
\renewcommand\thebibliography[1]{%
  \oldthebibliography{#1}%
  \setlength{\parskip}{\bibitemsep}%
  \setlength{\itemsep}{\bibparskip}%
}

\setlist{itemsep = 0.25em, topsep = 0.25em}

\declaretheoremstyle[spaceabove=0.25cm,spacebelow=0.25cm,notefont=\normalfont\bfseries, notebraces={(}{)}]{theorem}
\declaretheoremstyle[spaceabove=0.25cm,spacebelow=0.25cm,bodyfont=\normalfont,notefont=\normalfont\bfseries, notebraces={(}{)}]{noital}
\declaretheoremstyle[spaceabove=0.25cm,spacebelow=0.25cm,bodyfont=\normalfont\color{darkgreen},notefont=\normalfont\bfseries, notebraces={(}{)}]{green}
\declaretheoremstyle[spaceabove=0.25cm,spacebelow=0.25cm,bodyfont=\normalfont,notefont=\normalfont\bfseries,qed=$\qedsymbol$,notebraces={(}{)}]{proofstyle}

\numberwithin{equation}{section}

\makeatletter
\def\paragraph{\smallskip\@startsection{paragraph}{4}%
  \z@\z@{-\fontdimen2\font}%
  {\normalfont\itshape}}
\makeatother

\newcommand{\cB}{\ensuremath{\mathcal B}}

\newcommand{\cI}{\ensuremath{\mathcal I}}

\newcommand{\cF}{\ensuremath{\mathcal F}}

\newcommand{\cM}{\ensuremath{\mathcal M}}

\newcommand{\cO}{\ensuremath{\mathcal O}}

\newcommand{\cX}{\ensuremath{\mathcal X}}

\newcommand{\R}{\ensuremath{\mathbb R}}
\newcommand{\C}{\ensuremath{\mathbb C}}
\newcommand{\PP}{\ensuremath{\mathbb P}}
\newcommand{\Z}{\ensuremath{\mathbb Z}}

\newcommand{\half}{\ensuremath{\frac{1}{2}}}

\newcommand{\st}{\: | \:}

\newcommand{\oper}{\mathrm{op}}
\newcommand{\higgs}{\mathrm{H}}
\newcommand{\kahler}{K\"ahler\xspace}
\newcommand{\hk}{hyperk\"ahler\xspace}

\newcommand{\I}{{\mathrm i}}
\newcommand{\e}{{\mathrm e}}
\newcommand{\de}{\mathrm{d}}

\newcommand{\inst}{{\mathrm {inst}}}

\newcommand{\SL}{{\mathrm {SL}}}
\newcommand{\PSL}{{\mathrm {PSL}}}

\newcommand{\CP}{{\mathbb {CP}}}

\newcommand{\abs}[1]{\lvert#1\rvert}
\newcommand{\norm}[1]{\lVert#1\rVert}
\newcommand{\IP}[1]{\langle#1\rangle}
\newcommand{\eps}{\epsilon}

\newcommand{\rmsf}{\mathrm{sf}}

\newcommand{\ti}[1]{\textit{#1}}

\renewcommand{\leq}{\leqslant}
\renewcommand{\geq}{\geqslant}

\DeclareMathOperator{\im}{Im}
\DeclareMathOperator{\re}{Re}

\DeclareMathOperator{\End}{End}

\usepackage{stackengine}
\usepackage{scalerel}
\newcommand\pluscross{\scalerel*{\stackinset{c}{}{c}{}{$\boldsymbol{+}$}{$\boldsymbol{\times}$}}{+}}

\newcommand{\DE}{\mathrm{DE}}
\newcommand{\IEQ}{\mathrm{IEQ}}

\newcommand{\inswide}[1]{\makebox[\textwidth][c]{\includegraphics[width=1.13\textwidth]{#1}}}

\newcommand{\operfigs}{Figures \ref{fig:A1A2operth0.0}--\ref{fig:A2A2operth0.1}\xspace}
\newcommand{\hitchinfigs}{Figures \ref{fig:A1A2th0.0}--\ref{fig:A2A2th0.1}\xspace}

\begin{document}


\setcounter{page}{1}

\title[Opers and nonabelian Hodge: numerical studies]{Opers and nonabelian Hodge: numerical studies}
\author[E. Dumas]{Emily Dumas}
\author[A. Neitzke]{Andrew Neitzke}
\date{January 24, 2025 (First version: July 1, 2020).}

{\abstract{We present numerical experiments that test the predictions of a conjecture of Gaiotto-Moore-Neitzke and Gaiotto concerning the monodromy map for opers, the non-abelian Hodge correspondence, and the restriction of the \hk $L^2$ metric to the Hitchin section.  These experiments are conducted in the setting of polynomial holomorphic differentials on the complex plane, where the predictions take the form of conjectural formulas for the Stokes data and the metric tensor.  Overall, the results of our experiments support the conjecture.}}

\maketitle

\section{Introduction}

In this paper we present and discuss numerical experiments that compute two maps that arise naturally in Teichm\"uller theory: The non-abelian Hodge correspondence (NAHC) for the Hitchin section, and the monodromy map for opers (which is a particular case of the Riemann-Hilbert correspondence).
Both of these maps associate monodromy data to a tuple of holomorphic differentials on a Riemann surface.
Furthermore, these maps are expected to be asymptotic in a certain sense (e.g.~the conjecture of \cite{Katzarkov2015} and results of \cite{Simpson1991,Wolf1991,Loftin2007,Acosta2016,Dumas2017}).

The motivation for our experiments is a statement we dub the \ti{twistor Riemann-Hilbert conjecture},
which asserts that these monodromy maps can also be computed by solving a system of coupled integral equations.
Particularly in the case of the NAHC, such a description would be remarkable in that the integral equations involve only contour integrals of holomorphic functions, and do not involve the solution of any partial differential equation.
Moreover, the conjectural integral equation suggests a rather simple iterative strategy for computing the map which in the cases we study converges very rapidly and is computationally inexpensive. 

The twistor Riemann-Hilbert conjecture, in the form that we investigate it here,
is formulated in the works \cite{Gaiotto:2008cd,Gaiotto:2009hg,Gaiotto2012,Gaiotto2014}, originally motivated
by considerations of supersymmetric quantum field theory.
Many special cases of the conjecture had been discovered earlier; in particular,
for the monodromy of $\SL_2\C$-opers with potential of the form $z^n + c$, the conjecture appeared 
in the context of the ODE/IM correspondence pioneered in \cite{Dorey:1998pt} 
(see e.g. \cite{Dorey:2000kq,Dorey2007} for reviews of the ODE/IM correspondence).
More recently, statements related to the twistor Riemann-Hilbert conjecture for the NAHC 
have been studied in the form of the massive ODE/IM correspondence, beginning with 
\cite{Lukyanov_2010}.
The twistor Riemann-Hilbert conjecture is also closely related to the exact WKB method (and its extension
to higher rank); this is an extensive literature which we cannot review here, 
but see e.g.~\cite{Hollands2019} for a description
close to our point of view in this paper, and references therein.
Numerical investigations of special cases of 
the twistor Riemann-Hilbert conjecture (or closely related
statements) for $\SL_2\C$-opers have been 
described in e.g. \cite{dorey2004aspects,Dorey2007,Gaiotto2014,Ito2018,grassi2019nonperturbative}, and 
in
\cite{Dorey_2000} for $\SL_3\C$-opers. In particular, \cite{Ito2018} gives results of a numerical test of
a version of the conjecture for $\SL_2\C$-opers in 
the same examples we consider below, involving slightly different
quantities than we compute in this paper.
We are not aware of any numerical investigations of the twistor Riemann-Hilbert conjecture for the NAHC.

To test the twistor Riemann-Hilbert conjecture, we developed software to compute the monodromy maps directly, and using the conjecturally equivalent integral equations, and here we compare the results.
We only consider the case of polynomial Higgs bundles and opers on the complex plane.
This case is more amenable to computation than the case of a compact surface, though since the plane is simply-connected, there is no monodromy in the classical sense.
Instead, the monodromy-type invariants relevant to the correspondences are the Stokes data of the connections, as discussed in \cite{MR1864833}.
The bulk of this paper is thus devoted to discussing two methods for computing Stokes data (one of them conjectural) and comparing the results of numerical experiments with these methods.

In addition to allowing us to study the monodromy maps themselves, our implementation of the direct and integral equation approaches easily extends to compute the \hk metric on the moduli space of Higgs bundles, restricted to a \kahler metric on the Hitchin section.
While the integral equation side of this picture again gives a formula that is only conjectural, it is particularly appealing because it implies specific asymptotics of the metric which are not evident from its original definition.
Indeed, the numerical evidence supporting the conjectural formula for the \hk metric reported in this paper was part of the original inspiration for the authors' work in \cite{Dumas2018}, where an asymptotic formula with exponentially decaying error term was established for the Hitchin section of a compact surface in rank $2$.
(This exponential decay improved on the earlier asymptotic results of Mazzeo, Swoboda, Weiss, and Witt \cite{Mazzeo2019} in rank $2$, which had a polynomially decaying error term; exponential decay was later established in all ranks by Fredrickson \cite{Fredrickson2020}.)

\subsection{Concrete predictions}

Because a full description of the families of connections, monodromy maps, and integral equations is somewhat lengthy, we defer that to the subsequent sections.
But to give an indication of what the conjectural picture looks like, here is a concrete geometric consequence of it that is supported by our numerical experiments.

\smallskip

{\textit {Predictions for harmonic maps.}}
For any polynomial $P(z)$ with complex coefficients there exists a harmonic map $h : \C \to \mathbb{H}^2$, unique up to isometry, with Hopf differential $-4 P(z) \de z^2$ and such that the Riemannian metric $\abs{\partial h}^2 \abs{\de z}^2$ on $\C$ is complete \cite{WA94}.
Furthermore, the image of this map is the interior of an ideal polygon with $d+2$ vertices, where $d = \deg(P)$ \cite{Han95}.
Exactly which ideal polygon appears depends on the coefficients of $P$; for example, $P(z) = z^d$ gives the regular $(d+2)$-gon.
Except for a few symmetric examples like this one, no explicit formula is known which describes the polygon in terms of the polynomial $P$.

To consider a specific example, we might ask which isometry class of ideal pentagons in $\mathbb{H}^2$ corresponds to the cubic polynomial
\begin{equation}
	P(z) = z^3-1.
\end{equation}
An ideal pentagon is determined up to isometry by two real-valued invariants, which can be taken to be cross ratios of any two $4$-tuples of the vertices (considered as elements of $\mathbb{RP}^1$).
In the case of the pentagon associated to $z^3-1$, one can use the fact that the polynomial has real coefficients to show that the corresponding pentagon has a reflection symmetry about a geodesic passing through one of the ideal vertices.
This reduces the problem of characterizing the shape of the pentagon to determining 
a single real number;
for our purposes it will be convenient to take this parameter as $X := 1 + \chi(v_1,v_2,v_3,v_4)$, where $\chi(a,b,c,d) = \frac{(a-b)(c-d)}{(a-d)(c-b)}$ is the cross ratio and $v_1, \ldots, v_5$ are the ideal vertices, with $v_1$ fixed by the reflection symmetry.

Assuming the twistor Riemann-Hilbert conjecture, this invariant $X$ can be computed by solving the following integral equation.
First, define two exponentially decaying kernel functions $K_\pm$ by 
\begin{equation}
K_\pm(t) = A_\pm \cdot \frac{\cosh(t)}{2 \cosh(2t) \pm 1}
\end{equation}
where the real constants $A_\pm$ are chosen so that $\int_{-\infty}^{\infty} K_\pm(t)  \, \de t= 1$.
Also define the constant
\begin{equation}
M = \sqrt{3 \pi} \frac{\Gamma(4/3)}{\Gamma(11/6)}
\end{equation}
and let $x_0$ denote the function on $\R$,
\begin{equation}
x_0(t) = -2 M \cosh(t).
\end{equation}
Conjecturally, there is a unique smooth function $x : \R \to \R$ that grows as $O(x_0)$ and which satisfies the integral equation
\begin{equation}
\label{eqn:ieq-explicit-example}
x = x_0 - \big( K_+ \ast \log(1+\exp(x)) \big) 
\end{equation}
where $f \ast g$ denotes the convolution of $f$ and $g$ on $\R$.
In terms of this solution, the prediction for the cross ratio invariant is that
\begin{equation}
 X = \exp(y(0)) \;\; \text{ where } \;\;
 y = \frac{\sqrt{3}}{2} \, x_0 - \big(K_- \ast \log(1 + \exp(x))\big).
\end{equation}
In our experiments, computing a numerical solution of \eqref{eqn:ieq-explicit-example} leads to a predicted value $X \approx 0.006415703$, while computing a numerical approximation of the harmonic map itself gives $X \approx 0.006415699$ with an estimated error of $1.8 \times 10^{-8}$.
The method for solving the integral equation is based on writing \eqref{eqn:ieq-explicit-example} as a fixed point equation, $x = \cF(x)$, and then locating a fixed point by taking a limit of iterates $\cF^n(x_0)$.
This scheme converges very rapidly because $x$ is ultimately a very small correction to $x_0$, with $\|x-x_0\|_\infty < 0.001$ while $\inf |x_0| = 2M \approx 5.83$.
In comparison, our implementation of the direct approach to computing the harmonic map is rather expensive in computation time and memory, e.g.~requiring minutes to hours of computation time depending on the desired accuracy.

While we have described this example in terms of harmonic maps, it has an equivalent formulation in terms of Higgs bundles in the Hitchin section.
In that interpretation, which is the one used in the body of the paper, the isometry type of the ideal polygon in $\mathbb{H}^2$ corresponds to the Stokes data of the flat connection corresponding to a rank-$2$ bundle with Higgs field of the form $\left(\begin{smallmatrix}0 & P(z)\\1 & 0\end{smallmatrix}\right) \de z$.
This construction is detailed in \autoref{sec:families}, and the specific families considered in our experiments are introduced in \autoref{sec:experiments}.
In the terminology of the latter section, the harmonic maps problem phrased above corresponds to the $(A_1,A_2)$ family with parameters $c=1$, $\Lambda=0$, $\vartheta=0$, and $R=1$.
The values of $X$ given above appear in \autoref{fig:A1A2th0.0} over $R=1$, with markers {\color{mpl1}$\boldsymbol{+}$} and {\color{mpl1}$\boldsymbol{\times}$} for the direct harmonic map computation and integral equation prediction, respectively.
(The close agreement between the two methods creates the appearance of a single marker {\color{mpl1}\pluscross}.)

In addition to these rank-$2$ Higgs bundles, our experiments also consider their natural generalization to rank $3$ and polynomial cubic differentials; here a geometric interpretation can be given in terms of affine spheres, and in this interpretation our experiments involve computing the map $\alpha$ that is the main object of study in \cite{DW15}.

\smallskip

{\textit{Predictions for opers.}}
A minor variation on the rank-$2$ example described above illustrates the experimental study of opers in \autoref{subsec:results-opers}:
Rather than considering the Hopf differential of a harmonic map, we consider the holomorphic immersion $f : \C \to \mathbb{CP}^1$ with Schwarzian derivative $2 P(z) \, \de z^2$, which is unique up to composition with a linear fractional map.
Such a map (with $P(z)$ polynomial of degree $d$) distinguishes a cyclically ordered configuration of $d+2$ points on $\CP^1$ which are its asymptotic values.
The conjecture of \cite{Gaiotto2014} then expresses cross ratios of $4$-tuples of these points in terms of the solution of an integral equation that is a minor modification of \eqref{eqn:ieq-explicit-example}.
(This particular example of the twistor Riemann-Hilbert 
conjecture was first discovered as part of 
the ODE/IM correspondence \cite{Dorey:1998pt}.)
Again, our discussion of this example in the body of the paper uses a bundle description, where the map $f$ is replaced with the equivalent data of a $\SL_2\C$-oper over $\C$.
This is a certain type of flat holomorphic connection, whose connection $1$-form can be taken to have the form $\left(\begin{smallmatrix}0 & P(z)\\1 & 0\end{smallmatrix}\right) \de z$, and computing the Stokes data of this connection is equivalent to computing the asymptotic values of $f$.
In this case, we compare the integral equation predictions to the more direct approach of computing the Stokes data from the parallel transport operator of the flat connection (which is obtained by solving an ordinary differential equation).
Here again, we study the natural generalization of this picture to rank $3$.

\subsection{Summary and interpretation of results}

\begin{figure}[H]
  \inswide{{paperfigA1A2operth0.0small}.pdf}
  \hspace*{-3mm}\inswide{{paperfigA1A2operth0.0large}.pdf}\\
  \vspace*{4mm}\includegraphics[width=0.8\textwidth]{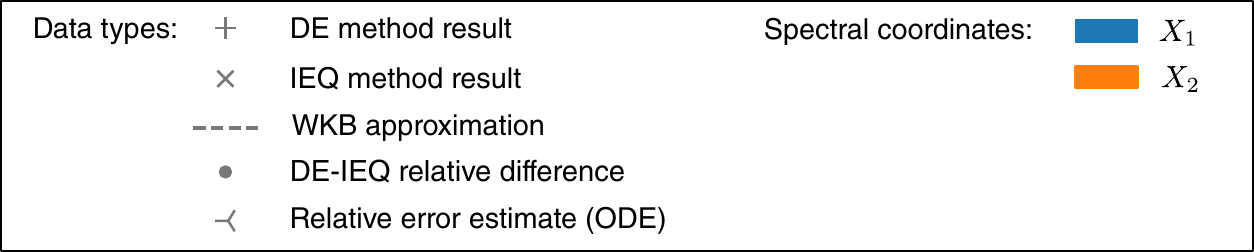}%
  \caption{Oper comparison: $(A_1,A_2)$ family at base point ($\Lambda=0$, $c=1$) with $\vartheta=0$.}
  \label{fig:A1A2operth0.0}
\end{figure}

\begin{figure}[H]
  \inswide{{paperfigA1A2operth0.1small}.pdf}
  \hspace*{-3mm}\inswide{{paperfigA1A2operth0.1large}.pdf}\\
  \vspace*{4mm}\includegraphics[width=0.8\textwidth]{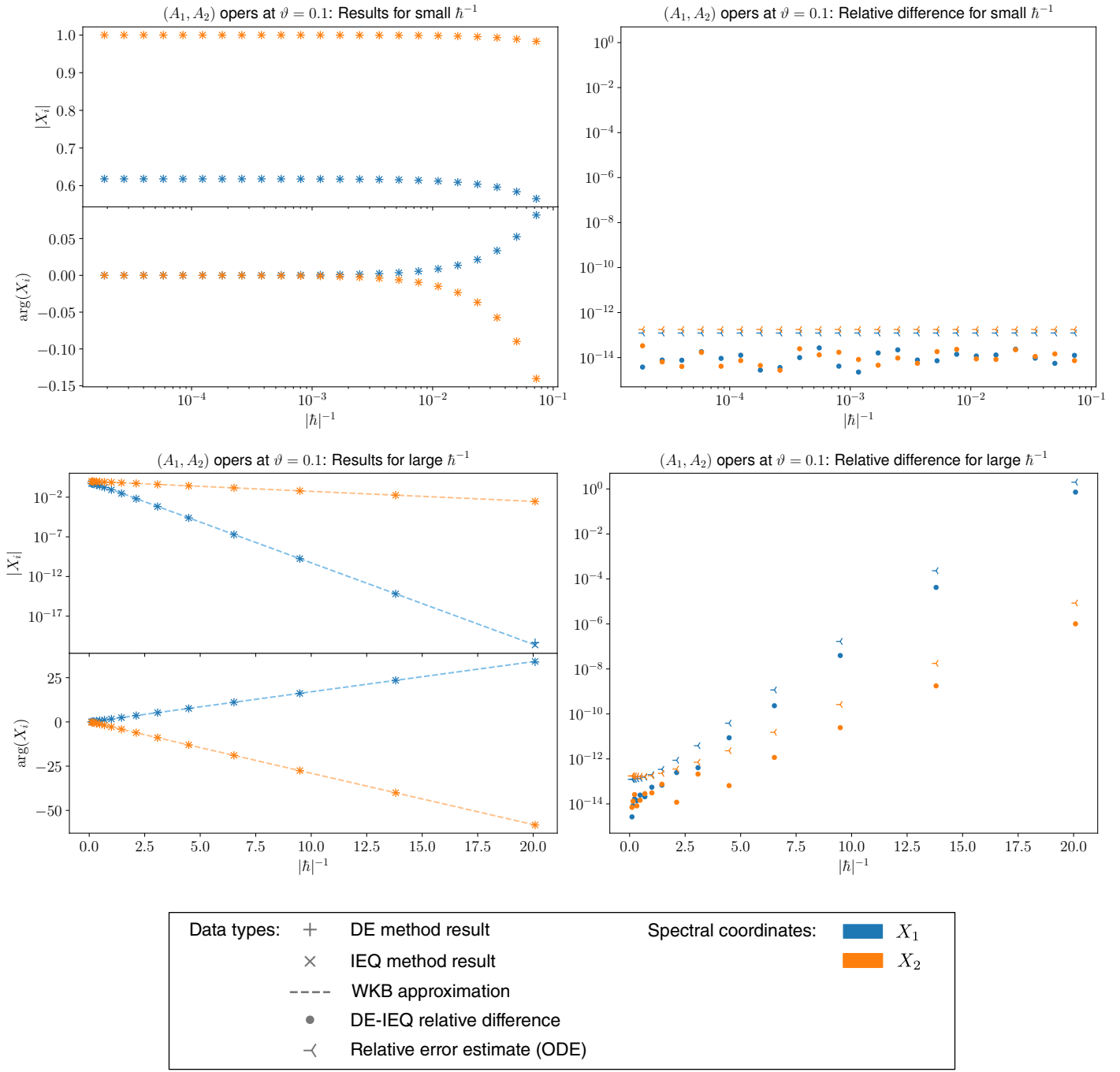}%
  \caption{Oper comparison: $(A_1,A_2)$ family at base point ($\Lambda=0$, $c=1$) with $\vartheta=0.1$.}
  \label{fig:A1A2operth0.1}
\end{figure}

\begin{figure}[H]
  \inswide{{paperfigA1A2operth0.0-deformedsmall}.pdf}
  \hspace*{0mm}\inswide{{paperfigA1A2operth0.0-deformedlarge}.pdf}\\
  \vspace*{4mm}\includegraphics[width=0.8\textwidth]{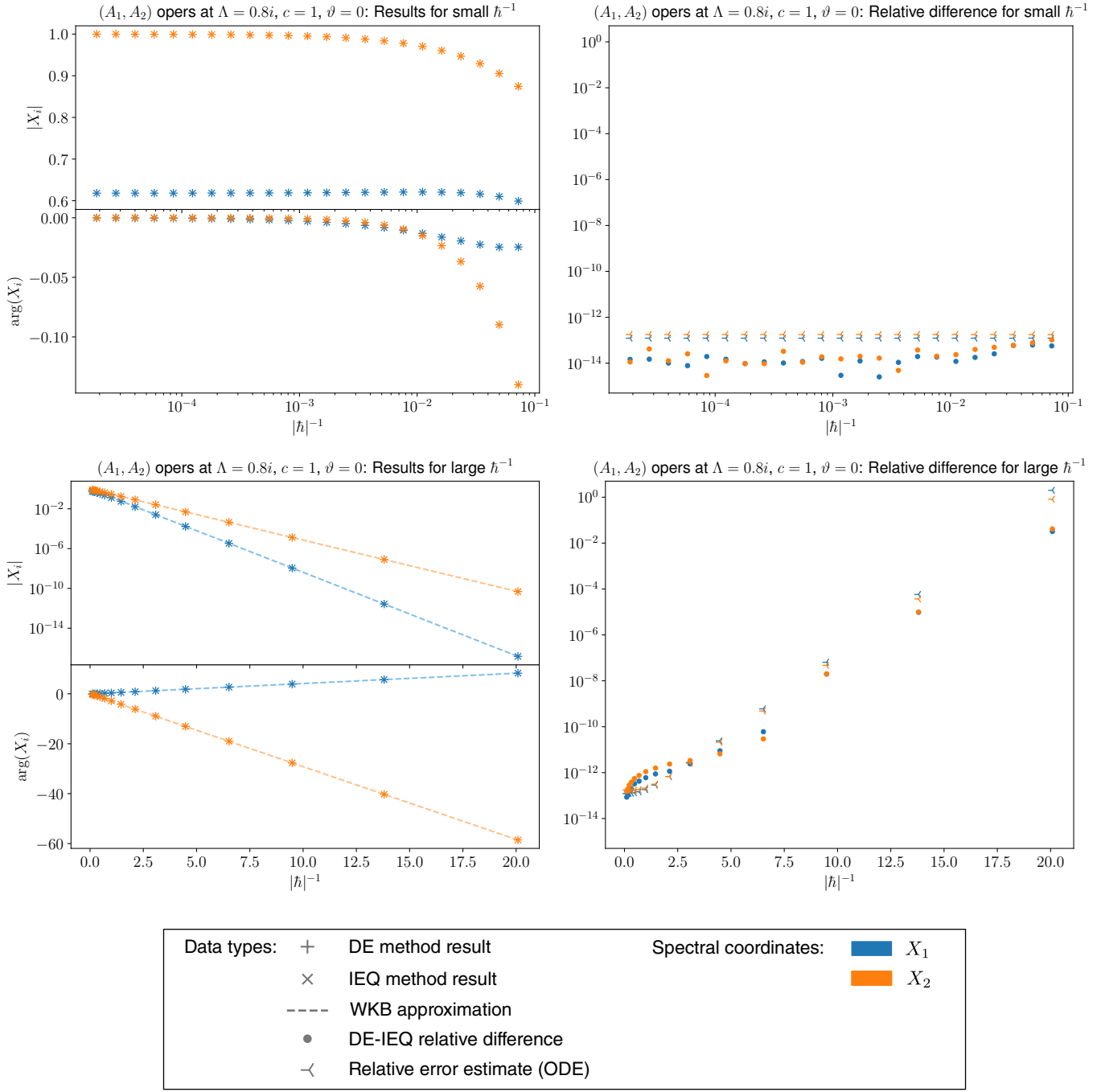}%
  \caption{Oper comparison: $(A_1,A_2)$ family at $\Lambda=0.8\I$, $c=1$ with $\vartheta=0$.}
  \label{fig:A1A2deformedoper}
\end{figure}

\begin{figure}[H]
  \inswide{{paperfigA1A3operth0.1small}.pdf}
  \inswide{{paperfigA1A3operth0.1large}.pdf}\\
  \vspace*{4mm}\includegraphics[width=0.8\textwidth]{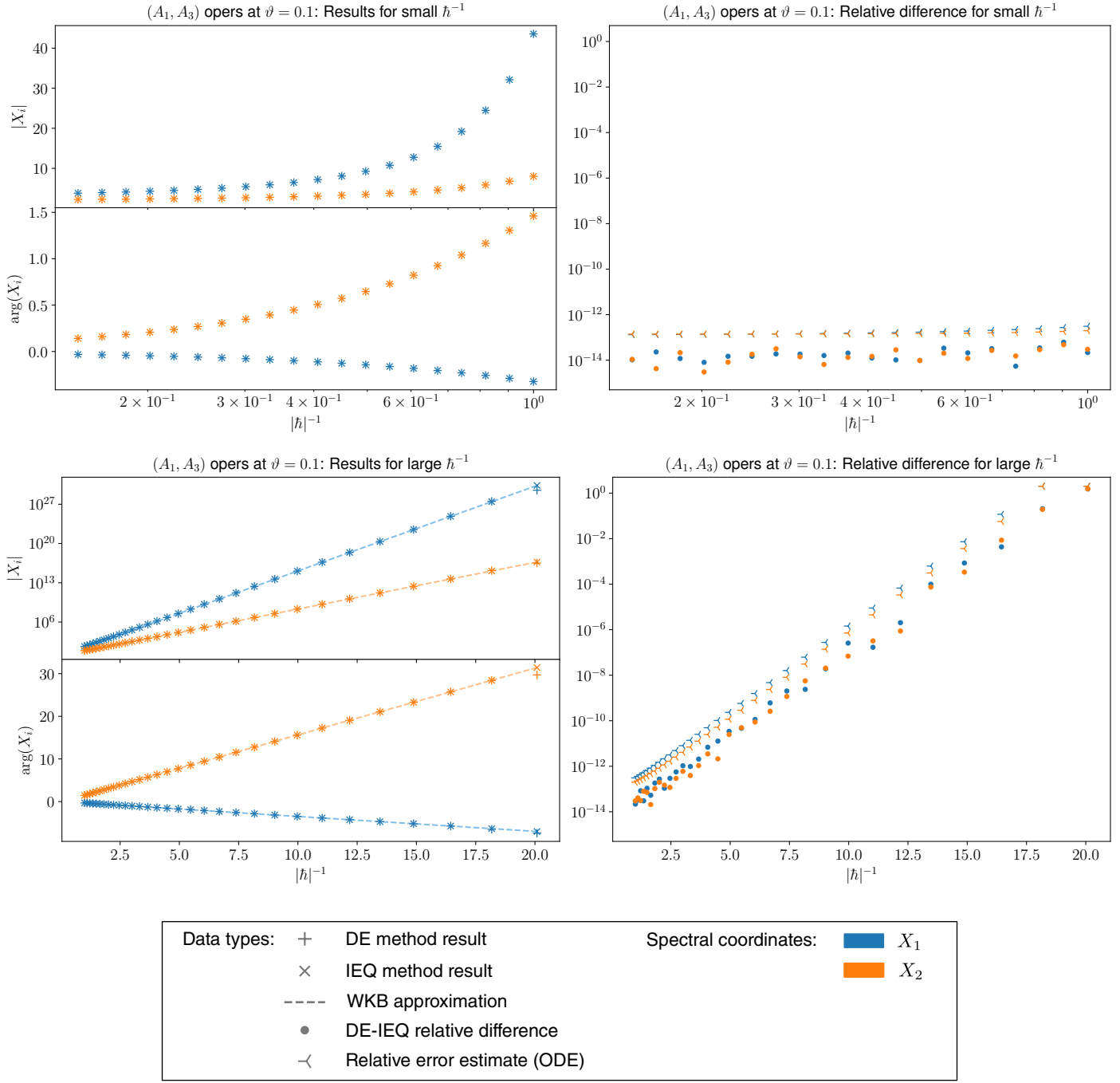}%
  \caption{Oper comparison: $(A_1,A_3)$ example with $\vartheta=0.1$. }
  \label{fig:A1A3operth0.1}
\end{figure}

\begin{figure}[H]
  \inswide{{paperfigA2A1operth0.0small}.pdf}
  \inswide{{paperfigA2A1operth0.0large}.pdf}\\
  \vspace*{4mm}\includegraphics[width=0.8\textwidth]{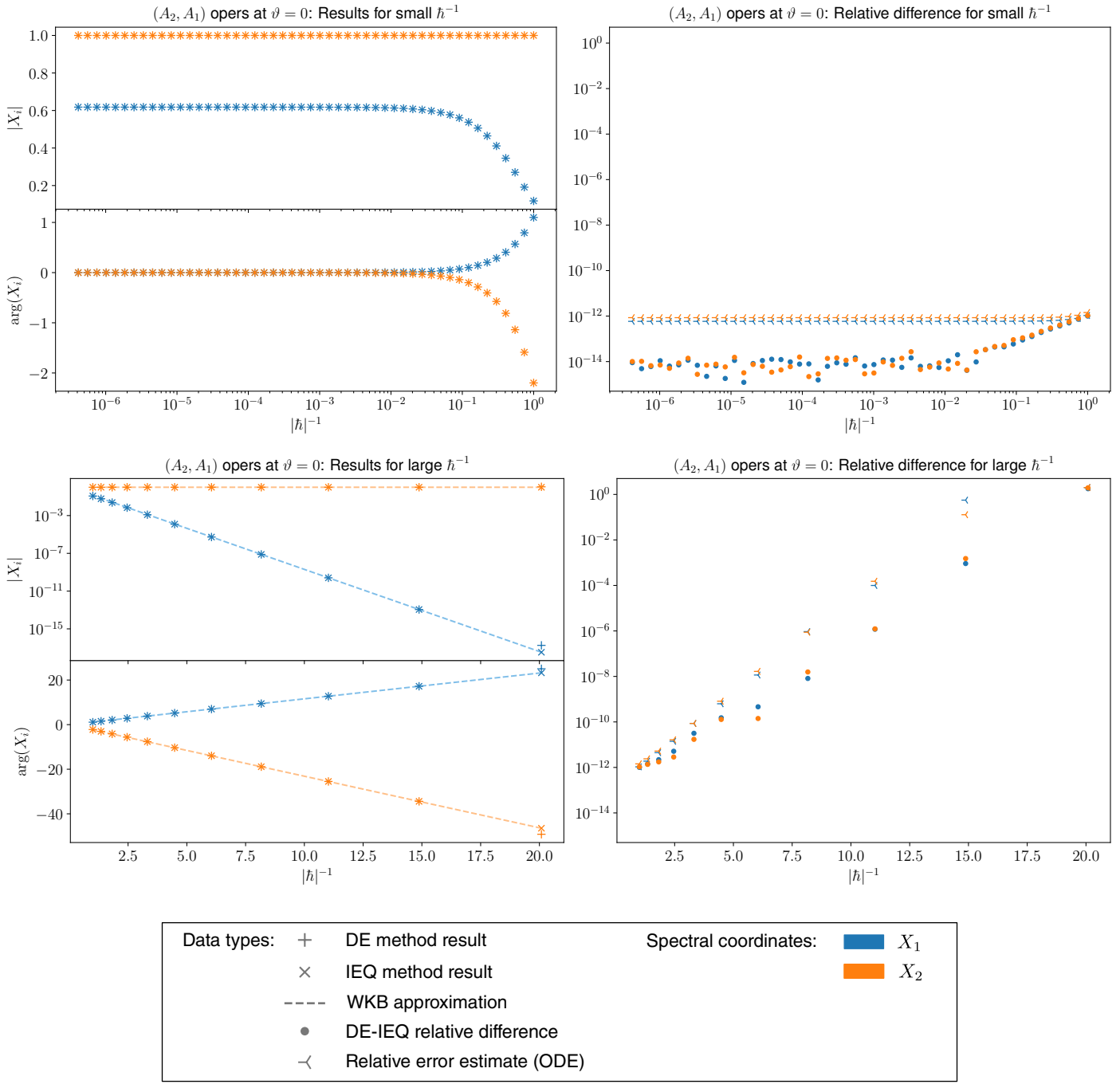}%
  \caption{Oper comparison: $(A_2,A_1)$ family at base point ($\Lambda=0$) with $\vartheta=0$.}
  \label{fig:A2A1operth0.0}
\end{figure}

\begin{figure}[H]
  \inswide{{paperfigA2A1deformedoperth0.0small}.pdf}
  \hspace*{-2mm}\inswide{{paperfigA2A1deformedoperth0.0large}.pdf}\\
  \vspace*{4mm}\includegraphics[width=0.8\textwidth]{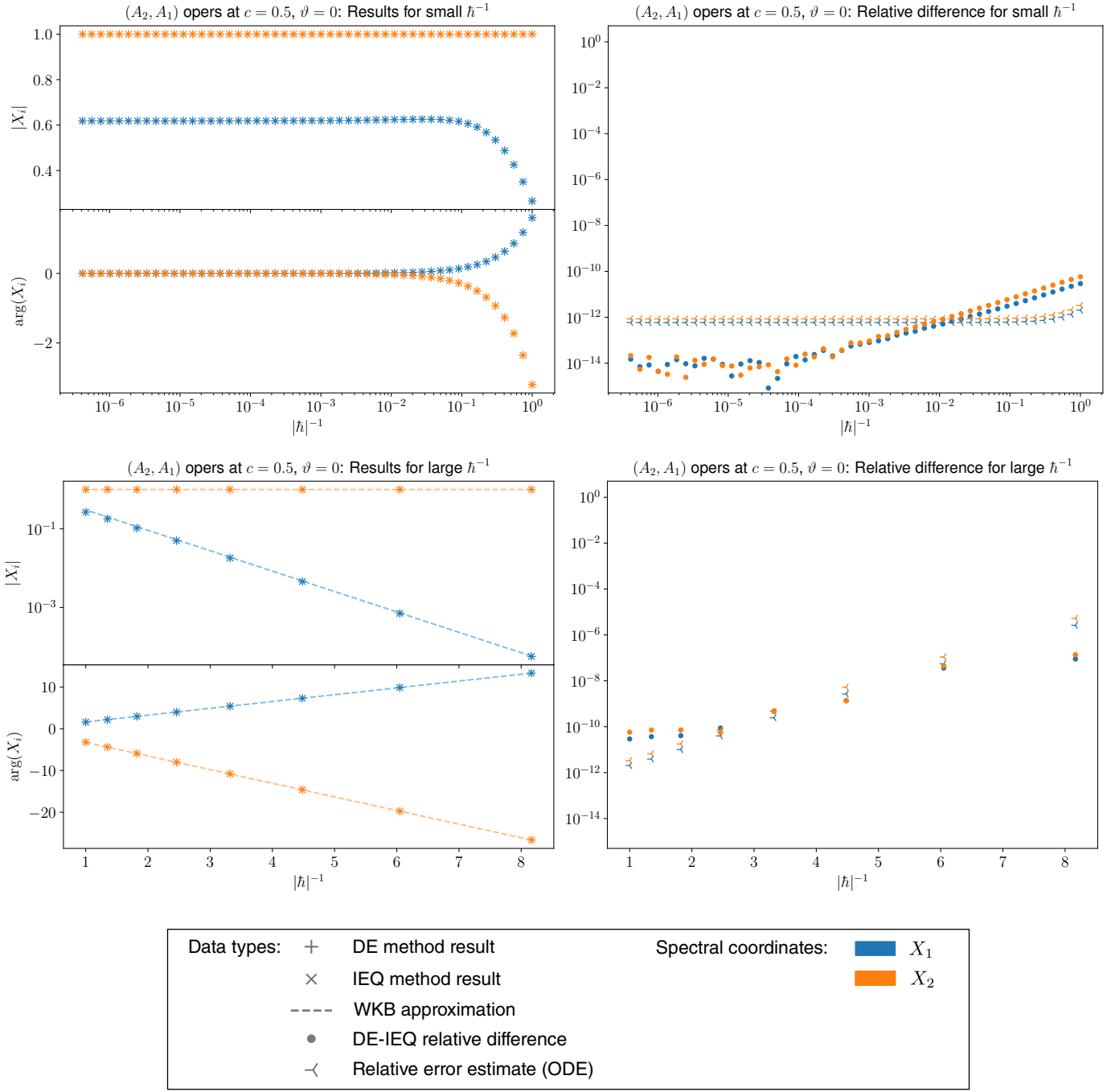}%
  \caption{Oper comparison: $(A_2,A_1)$ family at $c=\frac{1}{2}$ with $\vartheta=0$.}
  \label{fig:A2A1deformedoperth0.0}
\end{figure}

\begin{figure}[H]
  \inswide{{paperfigA2A2operth0.1small}.pdf}
  \hspace*{-4mm}\inswide{{paperfigA2A2operth0.1large}.pdf}\\
  \vspace*{4mm}\includegraphics[width=0.8\textwidth]{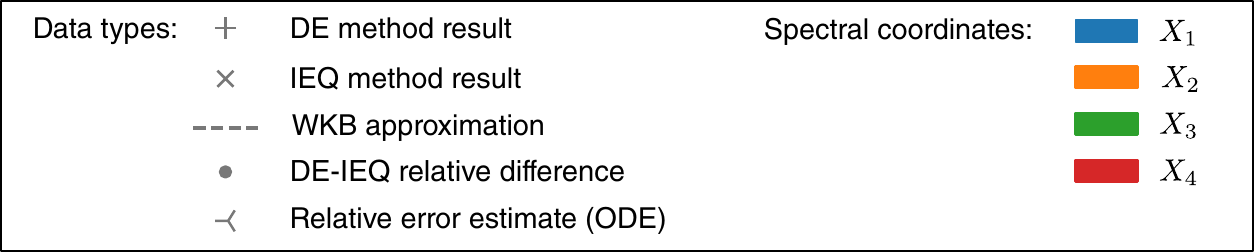}%
  \caption{Oper comparison: $(A_2,A_2)$ example with $\vartheta=0.1$.  The coordinates $X_3$ and $X_4$ are associated to pure flavor charges, so the WKB asymptotic is conjecturally an exact formula in those cases.}
  \label{fig:A2A2operth0.1}
\end{figure}

\begin{figure}[H]
  \inswide{{paperfigA1A2th0.0small}.pdf}
  \hspace*{-4mm}\inswide{{paperfigA1A2th0.0large}.pdf}\\
  \vspace*{4mm}\includegraphics[width=0.8\textwidth]{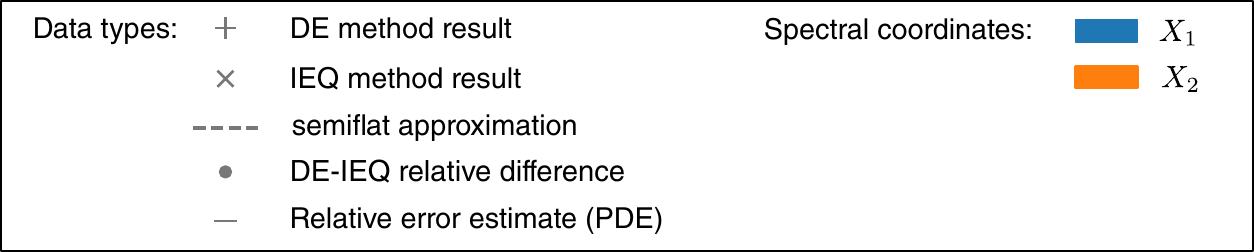}%
  \caption{Hitchin section comparison: $(A_1,A_2)$ family at base point ($\Lambda=0$, $c=1$) with $\vartheta=0$.}
  \label{fig:A1A2th0.0}
\end{figure}

\begin{figure}[H]
  \inswide{{paperfigA1A2th0.1small}.pdf}
  \hspace*{-3mm}\inswide{{paperfigA1A2th0.1large}.pdf}\\
  \vspace*{4mm}\includegraphics[width=0.8\textwidth]{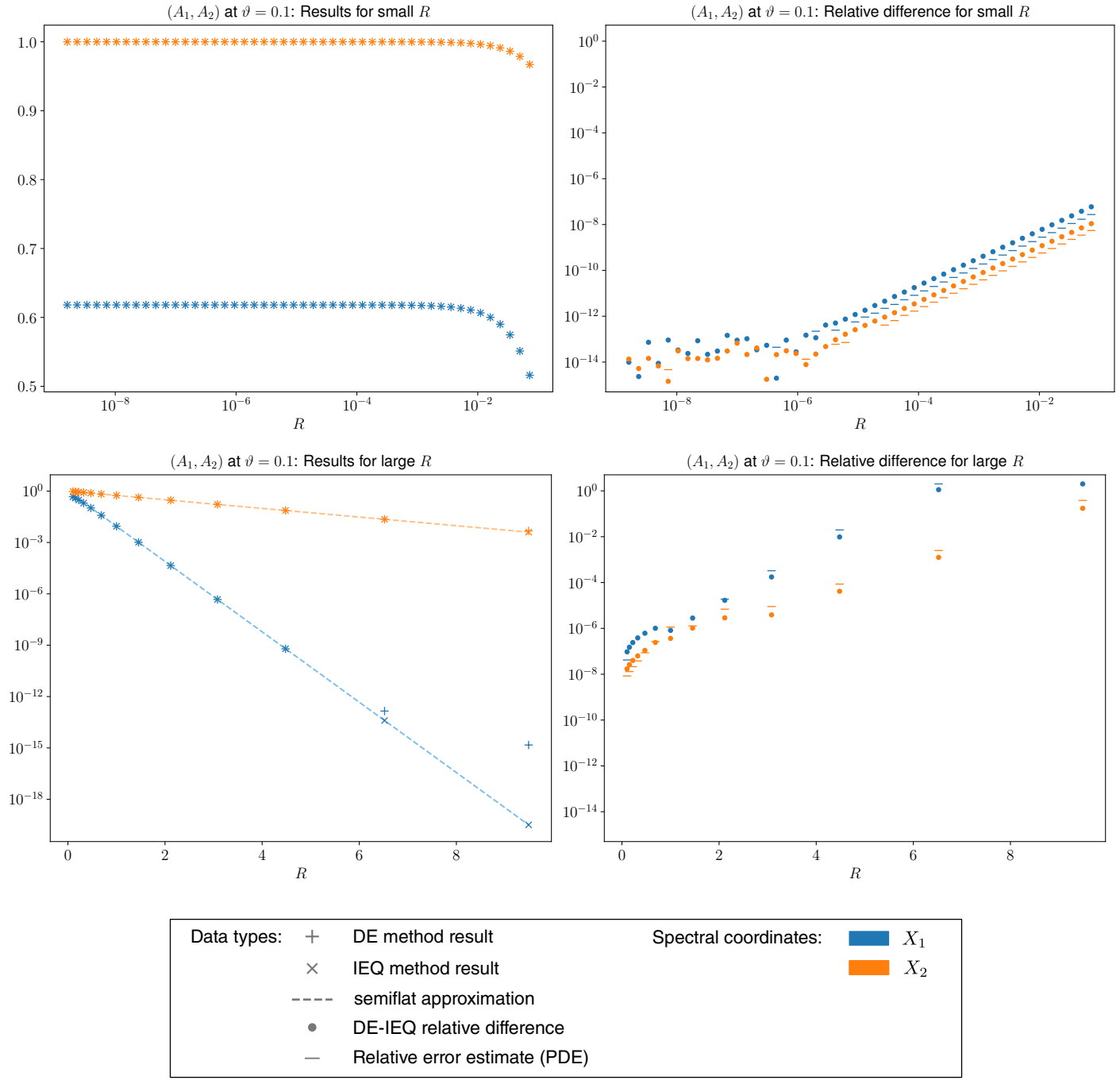}%
  \caption{Hitchin section comparison: $(A_1,A_2)$ family at base point ($\Lambda=0$, $c=1$) with $\vartheta=0.1$.}
  \label{fig:A1A2th0.1}
\end{figure}

\begin{figure}[H]
  \inswide{{paperfigA1A2th0.0-deformedsmall}.pdf}
  \hspace*{-3mm}\inswide{{paperfigA1A2th0.0-deformedlarge}.pdf}\\
  \vspace*{4mm}\includegraphics[width=0.8\textwidth]{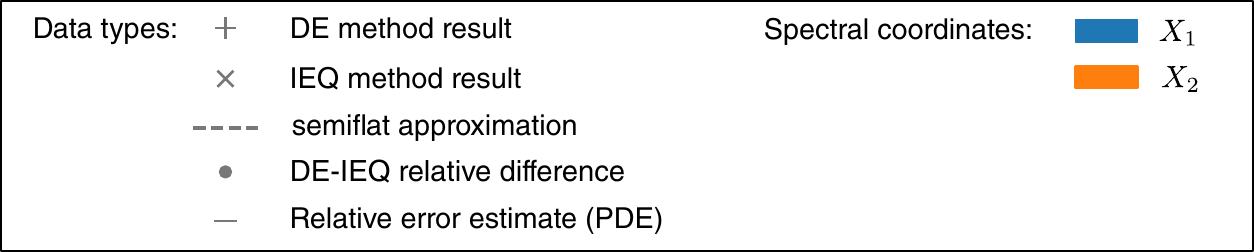}%
  \caption{Hitchin section comparison: $(A_1,A_2)$ family at $\Lambda=0.8\I$, $c=1$ with $\vartheta=0$.}
  \label{fig:A1A2deformed}
\end{figure}

\begin{figure}[H]
  \inswide{{paperfigA1A3th0.1small}.pdf}
  \hspace*{0mm}\inswide{{paperfigA1A3th0.1large}.pdf}\\
  \vspace*{4mm}\includegraphics[width=0.8\textwidth]{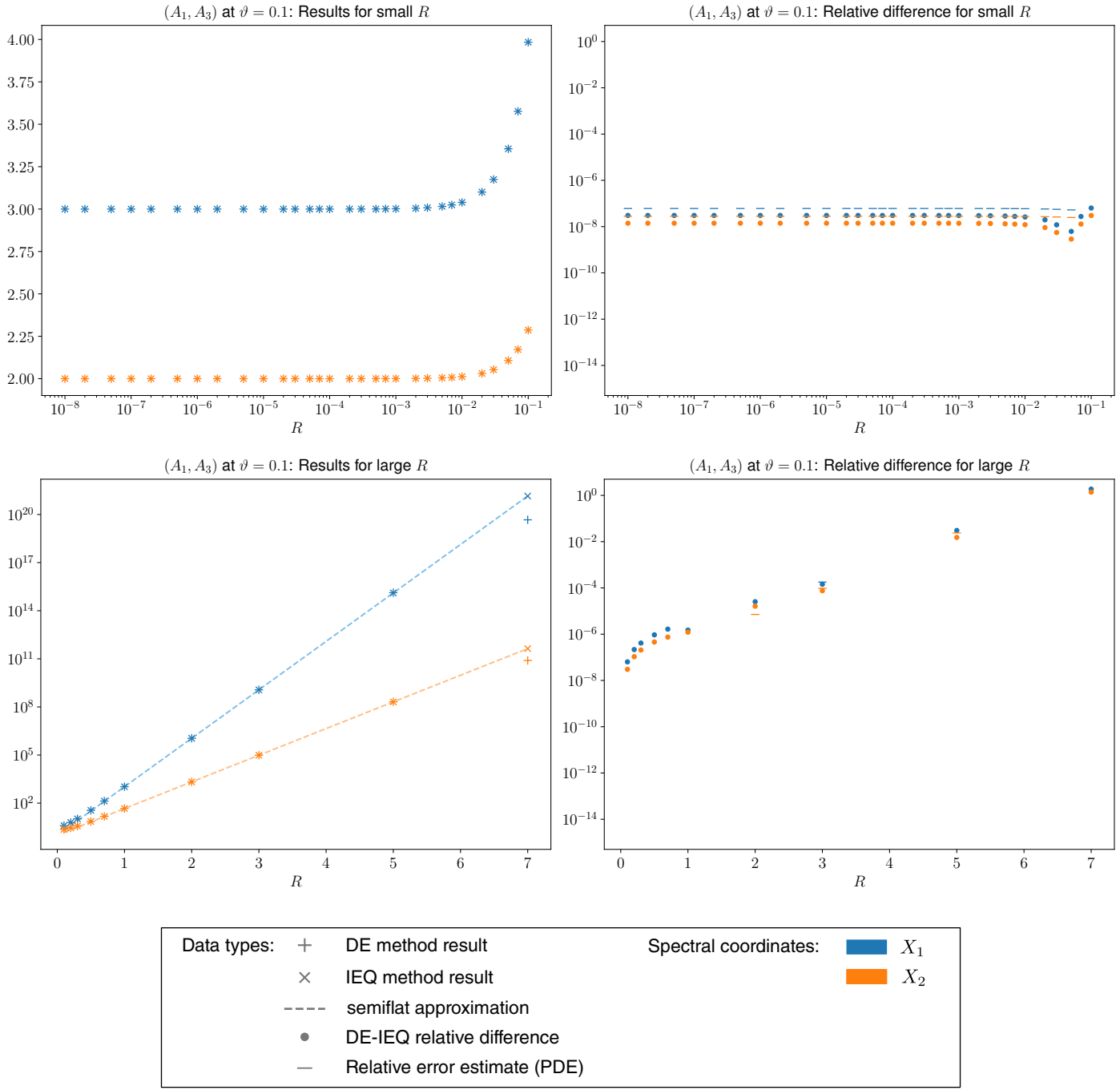}%
  \caption{Hitchin section comparison: $(A_1,A_3)$ example with $\vartheta=0.1$.}
  \label{fig:A1A3th0.1}
\end{figure}

\begin{figure}[H]
  \inswide{{paperfigA2A1th0.0small}.pdf}
  \hspace*{-2mm}\inswide{{paperfigA2A1th0.0large}.pdf}\\
  \vspace*{4mm}\includegraphics[width=0.8\textwidth]{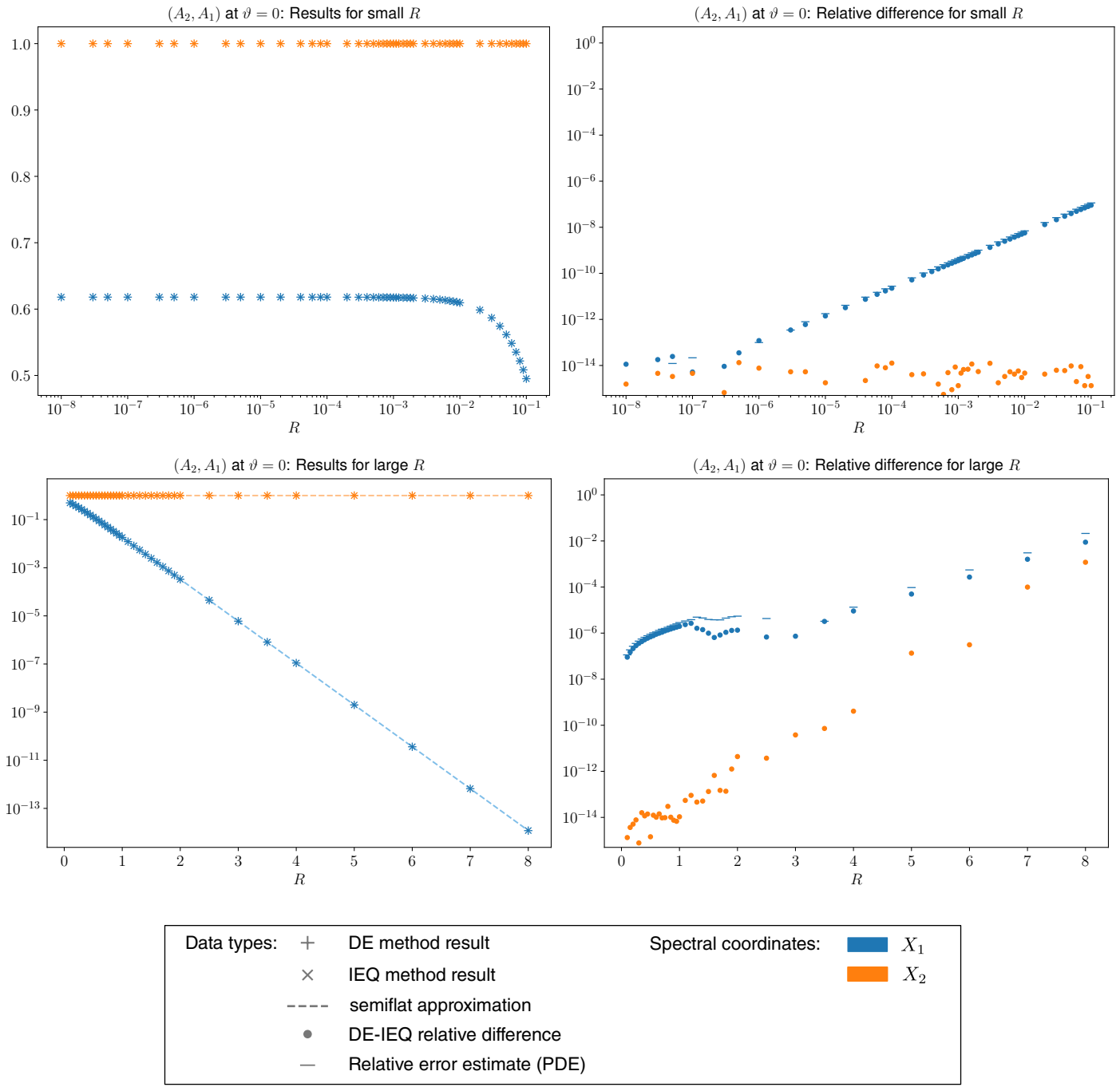}%
  \caption{Hitchin section comparison: $(A_2,A_1)$ family at base point ($c=0$) with $\vartheta=0$.}
  \label{fig:A2A1th0.0}
\end{figure}

\begin{figure}[H]
  \inswide{{paperfigA2A2th0.1small}.pdf}
  \hspace*{-2mm}\inswide{{paperfigA2A2th0.1large}.pdf}\\
  \vspace*{4mm}\includegraphics[width=0.8\textwidth]{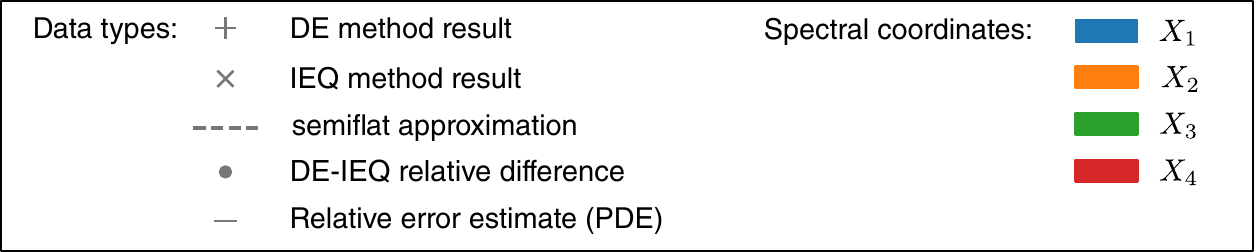}%
  \caption{Hitchin section comparison: $(A_2,A_2)$ example with $\vartheta=0.1$. The coordinates $X_3$ and $X_4$ are associated to pure flavor charges, so the semiflat asymptotic is conjecturally an exact formula in those cases.}
  \label{fig:A2A2th0.1}
\end{figure}

The main experiments we report in this paper involve computing and comparing Stokes data for $13$ one-parameter families of connections. The results are summarized in Figures \ref{fig:A1A2operth0.0}-\ref{fig:A2A2operth0.1}
for opers, and Figures \ref{fig:A1A2th0.0}-\ref{fig:A2A2th0.1} for the Hitchin section.
For one of these families we also compute the restriction of the \hk metric to the Hitchin section, for which
the results are summarized in Figure \ref{fig:paperfighkmetric}.
In general we believe that the results support the twistor Riemann-Hilbert conjecture.
Of course, this does not mean that we find exact agreement between the direct method and the
integral equation method; rather, it means that we believe that the difference we
do find can be accounted for by numerical error.

Each one-parameter family of connections which we study 
involves a positive real-valued parameter ($\abs{\hbar}^{-1}$ for opers, $R$ for the Hitchin section) with the property that the expected numerical error in the direct method 
grows rapidly with $R$ or $\abs{\hbar}^{-1}$.
It must therefore be expected that the difference between the results from the two methods may exhibit the same type of growth, even if the twistor Riemann-Hilbert conjecture holds, and this is indeed
what we find: in general the difference is small for small values of the parameter,
and grows as the parameter is increased.

Though we do not conduct a complete numerical analysis of both methods, we do analyze certain sources of numerical error, and ultimately conclude that our experiments do not provide any strong candidates for counterexamples to the conjecture.
Correspondingly, the breadth and variety of examples we have studied without finding an apparent counterexample may be seen as evidence toward the conjecture.

\subsection{Code and data}
\label{subsec:code}

All of our experiments were performed using implementations of the direct and integral equation methods we developed in Python.
The datasets resulting from our experiments, which were used to generate the plots and figures in this paper, are available at \cite{data}.
The source code for our program, with installation instructions and some documentation of the interfaces, is available at \cite{code}.
The code includes scripts to reproduce our experiments from scratch (taking $\sim\!130$ CPU-days on a fast machine in mid-2020) or to regenerate the tables and plots using the prepared dataset (which is of course much faster).

\subsection{Outline}

\mbox{}

\autoref{sec:families} introduces the Hitchin section and the family of opers (in the meromorphic case we consider), their associated Stokes data, and the \hk metric.

\autoref{sec:ieq} describes the conjectural integral equations for the Stokes data.

\autoref{sec:experiments} lists the specific connection families that we study, and reports the results of our experiments with Stokes data.

\autoref{sec:hk-experiments} reports the results of our experiments with the \hk $L^2$ metric.

\autoref{sec:gallery} is a small gallery of images related to the experiments of the previous sections.

\autoref{sec:implementation} gives a more detailed description of the computational methods used to produce the results reported in Sections \ref{sec:experiments}--\ref{sec:hk-experiments}, including e.g.~the specific parameter values (grid sizes, tolerances, etc.) used in the calculations.

\autoref{sec:sample} shows an example calculation using our code.

\autoref{sec:discussion} discusses the results of our experiments and the outlook for extensions of this work in the future.

\subsection{Acknowledgements}

The authors thank Philip Boalch, Qiongling Li, Marcos Marino, Rafe Mazzeo, and David Nicholls for helpful conversations related to this work, and the anonymous referees for helpful corrections and suggestions.
They also thank the Mathematical Sciences Research Institute for supporting their participation in the Fall 2019 program ``Holomorphic Differentials in Mathematics and Physics'' where some of this work was completed.
Some of the numerical experiments reported here were conducted using a computer cluster managed by the Advanced Cyberinfrastructure for Education and Research (ACER) group at the University of Illinois at Chicago.
We thank the Yale Center for Research Computing for guidance and use of the research computing infrastructure, specifically the Grace cluster.
The authors were supported by the U.S.~National Science Foundation, through individual grants DMS 1709877 (DD) and DMS 1711692 (AN), and their participation in the 2019 MSRI program was supported by DMS 1440140 (AN) and DMS 1107452, 1107263, 1107367, ``RNMS: GEometric structures And Representation varieties'' (the GEAR Network) (DD).

\section{Families of connections and monodromy data}
\label{sec:families}

\subsection{Hitchin base}

Unless otherwise indicated, all of the bundles we consider are vector bundles equipped with a holomorphic structure.
We will study certain families of connections on bundles over $\C$, and denote the rank of the bundle by $N$.
Later we focus on the cases $N=2,3$, though for the moment our discussion is general.

The connections we consider are associated to tuples of holomorphic differentials on $\C$.
We denote the canonical line bundle of $\C$ by $K$, and let $Q_k \subset H^0(\C,K^k)$ denote the space of holomorphic $k$-differentials on $\C$ of the form $\phi = P(z) \,\de z^k$ for $P$ a polynomial; equivalently $Q_k$ is the space of sections of $K^k$ which extend meromorphically to $\hat{\C}$.
Define
\begin{equation}
  \cB = \bigoplus_{k=2}^{N} Q_k.
\end{equation}
We write a typical point $w \in \cB$ as $w = (\phi_2, \ldots, \phi_N)$ and always denote by $P_k$ the polynomial such that $\phi_k = P_k(z) \, \de z^k$.

While $\cB$ is infinite-dimensional, the following finite-dimensional affine subspaces will be our main focus.
Let $\cB_d$ denote the set of $w \in \cB$ such that
\begin{equation}
  \label{eqn:bd-def}
  \deg(P_N) = d \; \text{ and } \deg(P_k) \leq \frac{k}{N} d \; \text{ for all }k.
\end{equation}

We call $\cB_d$ the \emph{universal Hitchin base} of $\C$ of degree $d$ (and rank $N$).
The terminology is meant to indicate that this ``universal'' space is not the base space of a Lagrangian fibration of a \hk manifold (which would generalize the character variety of a compact surface group) but is foliated by subspaces with that property.

A further subset of $\cB_d$ contains all of the examples we study numerically:  Let $\cB_d'$ denote the subset where $\deg(P_k) < \frac{k}{N}$ for all $k \neq N$, i.e.~where the inequality of \eqref{eqn:bd-def} is strict whenever possible.  Some of the constructions we make below are simpler to state for $\cB_d'$.

\subsection{Higgs fields from the Hitchin base}

Let $E = \cO^N$ denote the trivial bundle over $\C$ of rank $N$.
A holomorphic section of $\Omega^1(\C,\End(E))$ is a \emph{Higgs field} on $E$.

A key construction we use is a map that associates to $w \in \cB$ a Higgs field $\varphi_w$ on $E$. While $\varphi_w$ can be defined for any $N$, we give the explicit formulas only for $N=2,3$.
For $N=2$ and $w = (P_2 \, \de z^2)$ define
\begin{equation} \label{eq:oper-2}
  \varphi_w = \begin{pmatrix}
    0 & -P_2 \\
    1 & 0
  \end{pmatrix} \, \de z.
\end{equation}
For $N=3$ and $w = (P_2 \, \de z^2, P_3 \, \de z^3)$ define
\begin{equation} \label{eq:oper-3}
  \varphi_w = \begin{pmatrix}
    0 & -\half P_2 & -P_3 \\
    1 & 0 & -\half P_2 \\
    0 & 1 & 0 
  \end{pmatrix} \, \de z.
\end{equation}

\subsection{Polynomial opers} \label{sec:opers}

In general, $\SL_N\C$ opers are certain holomorphic connections on filtered vector bundles over Riemann surfaces.
When the base Riemann surface is $\C$, the bundle and filtration can be holomorphically trivialized, allowing some simplification of the definition in this case.
Since this is the only case we will use, we present only the simplified definition.
Discussion of the general case can be found in e.g. \cite{Wentworth2016,Dumitrescu2016}.

Let $\de$ denote the trivial connection on $E$.  For $w \in \cB$ define
\begin{equation} \label{eq:oper-connection}
  \nabla^{\oper}_w = \de + \varphi_w
\end{equation}
where $\varphi_w$ is as defined above.
This is a flat holomorphic connection on $E$.
The family of connections $\{ \nabla^{\oper}_w \st w \in \cB_d \}$ is the set of \emph{polynomial $\SL_N\C$ opers on $\C$ of degree $d$}.

It will be convenient to extend \eqref{eq:oper-connection} to a 1-parameter family, parameterized by $\hbar \in \C^\times$:
\begin{equation} \label{eq:oper-family}
  \nabla^\oper_w(\hbar) = \de + \hbar^{-1} \varphi_w.
\end{equation}
Passing from \eqref{eq:oper-connection} to \eqref{eq:oper-family} does not bring in 
any essentially new connections: indeed $\nabla_w^\oper(\hbar)$ is equivalent
to $\nabla_{\hbar^{-1} w}^\oper$, where we define $t w \in \cB_d$ by
\begin{equation} \label{eq:scaling-differentials}
  t w = (t^2 P_2 \de z^2, t^3 P_3 \de z^3).
\end{equation}
Nevertheless the main results below are naturally phrased in the language
of the families \eqref{eq:oper-family}.

\subsection{Polynomial Hitchin sections}
\label{subsec:hitchin}

Over a compact Riemann surface, the $\SL_N\C$ Hitchin section is a collection of Higgs bundles where the vector bundle is a certain direct sum of powers of the canonical bundle, and where the Higgs field has the form of $\varphi_w$ relative to that splitting.
We refer to the original papers \cite{Hitchin1991,Hitchin1987} or the recent survey \cite{Fredrickson2019} for further discussion of this theory.
We will consider the natural analogue of this family for the base Riemann surface $\C$, with a polynomial growth condition at infinity, and will take advantage of the holomorphic triviality of the canonical bundle of $\C$ to simplify our presentation.

In this section we only consider $N=2,3$.
For $w \in \cB$, the pair $(E,\varphi_w)$ is a Higgs bundle with wild ramification at $\infty$ in the sense of \cite{wnh}.
For such a bundle, we consider a hermitian metric $h$ that with respect to the splitting $E = \cO^N$ has the matrix form
\begin{equation}
    \begin{pmatrix}
      \e^{-u} & 0 \\
      0 & \e^{u}
    \end{pmatrix}
    \, \text{ if }N=2, \qquad
    \begin{pmatrix}
      \e^{-u} & 0 & 0 \\
      0 & 1 & 0 \\
      0 & 0 & \e^{u}
    \end{pmatrix}
    \;\; \text{ if }N=3,\\
\end{equation}
for a scalar function $u$ on $\C$.
Fix a degree $d$ and restrict attention to $w \in \cB_d$ for the moment.
Also assume that $\phi_k = 0$ for $k \neq N$ (as this holds in all examples of the Hitchin section that we investigate numerically).
We say that $h$ is a \emph{harmonic metric} if it satisfies both the \emph{self-duality equation}
\begin{equation} \label{eq:self-duality}
  \Delta u = 4\left ( \e^{k u} - \e^{-2u} |P_N|^2 \right ), \quad k = \frac{2}{N-1}
\end{equation}
and the \emph{compatibility condition}
\begin{equation} \label{eq:compatibility}
  u \sim \tfrac{N-1}{N} \log |P_N| \quad \text{ as } \abs{z} \to \infty.
\end{equation}
In more invariant terms, the self-duality equation \eqref{eq:self-duality} 
is equivalent to requiring $F_{D_h} + [ \varphi, \varphi^{\dagger_h} ] = 0$ where $D_h$ is the Chern connection of $h$ and $F_{D_h}$ is its curvature, and the compatibility condition says that the metric is compatible with a certain filtration at $\infty$ (see \cite{Fredrickson2017}). 
It is known in this case \cite{wnh,wild-harmonic,mochizuki2019good} that there exists a unique harmonic metric on $(E,\varphi_w)$ for $w \in \cB_d$, which we denote by $h_w$.
We call the collection $\{ (E,\varphi_w,h_w) \st w \in \cB_d \}$ of Higgs bundles with harmonic metrics the \emph{polynomial Hitchin section of degree $d$}.

Associated to $(E,\varphi_w,h_w)$ there is the flat (non-holomorphic) connection
\begin{equation} \label{eq:hitchin-section-connection}
\nabla^{\higgs}_w = D_{h_w} + \varphi_w + \varphi_w^{\dagger_{h_w}}.
\end{equation}
As we did in \autoref{sec:opers}, we will find it convenient to extend \eqref{eq:hitchin-section-connection} 
to a family of flat connections.
In this case we introduce two parameters: $R \in \R_+$ and $\zeta \in \C^\times$.
Then we define
\begin{equation} \label{eq:hitchin-section-family}
\nabla^\higgs_w(R, \zeta) = D_{h_{Rw}} + \zeta^{-1} \varphi_{Rw} + \zeta \varphi_{Rw}^{\dagger_{h_{Rw}}}.
\end{equation}
One sees immediately from \eqref{eq:hitchin-section-family}
that the parameter $R$ can be absorbed in a rescaling of $w$.
The same is true of the phase of $\zeta$: in particular,
if $\abs{\zeta} = 1$, then we have 
$\zeta^{-1} \varphi_{Rw} = \varphi_{\zeta^{-1} R w}$,
$h_{Rw} = h_{\zeta^{-1} R w}$, and $\zeta \varphi_{Rw}^{\dagger_{h_{Rw}}} = 
 \varphi_{\zeta^{-1} Rw}^{\dagger_{h_{\zeta^{-1} Rw}}}$,
 giving altogether
  $\nabla^\higgs_w(R, \zeta) = \nabla^\higgs_{\zeta^{-1} R w}$.
In contrast, when $\abs{\zeta} \neq 1$ it cannot be absorbed in a
rescaling of $w$; the connections $\nabla^\higgs_w(R,\zeta)$ for
 $\abs{\zeta} \neq 1$ are genuinely different from those for $\abs{\zeta} = 1$.

\subsection{Stokes data}
\label{sec:stokes}

In the sequel we consider monodromy-like invariants of the flat connections $\nabla^{\oper}_w(\hbar)$
and $\nabla^\higgs_w(R,\zeta)$.
Since the base is the simply connected space $\C$, these connections have no monodromy in the traditional sense.
Instead, their generalized monodromy is defined using Stokes data --- 
concretely, growth rates of sections as $z \to \infty$.
Stokes data for irregular connections 
are in many important respects parallel to 
monodromy data for regular connections; 
in particular, the corresponding moduli spaces have natural
holomorphic-symplectic and even
\hk structures, as shown in \cite{MR1864833,wnh,MR3931781}.

Suppose $w \in \cB_d'$, and let $A$ be the leading coefficient of $\hbar^{-N} P_N$ or
$\zeta^{-N} P_N$.
Then the \emph{Stokes sectors} of $\nabla^\oper_w(\hbar)$ or $\nabla^\higgs_w(R,\zeta)$ are the sets
\begin{equation}
  \label{eqn:monic-sectors}
  \arg(z) \in \left [ \frac{\pi (2j-1-\arg(A))}{d+N}, \frac{\pi (2j+1-\arg(A))}{d+N} \right ], \; |z| > r
\end{equation}
for $1 \leq j \leq d+N$ and $r > 0$; these give a collection of evenly spaced sectors about $\infty$.
In each such sector, there is a horizontal section of $\nabla^{\oper}_w(\hbar)$ or $\nabla^\higgs_w(R,\zeta)$ that decays exponentially as $z \to \infty$ within the sector, and this section is unique up to multiplication by a complex scalar.
This is a \emph{subdominant section} for that sector, and the line containing all subdominant sections for a sector is the \emph{subdominant line}.

All of the subdominant solutions associated to $\nabla^{\oper}_w(\hbar)$ or $\nabla^\higgs_w(R,\zeta)$ live in the $N$-dimensional space of horizontal sections over $\C$.
The relative position of the subdominant lines give moduli for the connection, sometimes called \ti{Stokes data}.
While the traditional approach to the Stokes phenomenon also includes a specific encoding of such data in so-called Stokes factors and Stokes matrices, these will not be directly used in what follows.
Instead, we define certain determinantal invariants of the connection using the subdominant sections
directly.

First, fix subdominant sections $s_j$ for $j=1, \ldots, d+N$.
Now let $(i_1,\ldots,i_N)$ be a tuple of distinct integers between $1$ and $d+N$, and define
\begin{equation}
  p(i_1, \ldots, i_N) = \det \left ( s_{i_1}, \ldots, s_{i_N} \right ).
\end{equation}
In the case $N=3$ and for a sextuple $(a,b,c,d,e,f)$ of integers between $1$ and $d+N$, we also define the ``hexapod'' determinant
\begin{equation}
  q(a,b,c,d,e,f) = \det \left ( s_a \times s_b, s_c \times s_d, s_e \times s_f \right ).
\end{equation}
Here we use any identification of the space of flat sections with $\C^3$ in order to compute the cross product.

The quantities $p$ and $q$ defined above are not invariants of the connection, since for example the replacements $s_i \mapsto \lambda_i s_i$ scale these determinants by some product of the factors $\lambda_i$.
However, a ratio of products of such determinants is an invariant if each $s_i$ appears the same number of times in the numerator and denominator.
For example, the quantity
\begin{equation}
  \frac{p(a,b)p(c,d)}{p(a,d)p(c,b)}
\end{equation}
is an invariant in the $N=2$ case; it is the cross ratio of the lines spanned by $s_a, s_b, s_c, s_d$.
Similar ratios of products of determinants give invariants for $N=3$, as does the ratio
\begin{equation}
  \frac{q(a,b,c,d,e,f)}{p(a,b,c)p(d,e,f)}.
\end{equation}

The invariants for the family of connections
$\nabla^\oper_w(\hbar)$ or $\nabla^\higgs_w(R,\zeta)$ vary analytically with $\hbar$ or $\zeta$
respectively.
When $\abs{\zeta} = 1$, $\nabla^\higgs_w(R,\zeta)$ admits a reduction to $\SL_N\R$, 
which implies that the invariants are real on this locus; otherwise
they are generally complex.

\subsection{Direct numerical calculation of Stokes data}
\label{subsec:direct-method}

Given $w \in \cB_d'$, the determinantal invariants for the associated oper connection $\nabla^\oper_w(\hbar)$ are relatively easy to compute numerically.
Writing the horizontal section equation $\nabla^{\oper}_w(\hbar) s = 0$ as a system of ordinary differential equations, we can use a numerical ODE solver to compute the parallel transport operator $A_\theta(\hbar) \in \mathrm{SL}_N\C$ of $\nabla^{\oper}_w(\hbar)$ from a fixed base point $z_0 \in \C$ to $z_0 + r \exp(\I \theta)$.
For large $r$ and for $\theta$ corresponding to an interior point of one of the Stokes sectors, the eigenvector of $A_\theta(\hbar)$ with smallest eigenvalue is an approximation of the subdominant section (considered as an element of the fiber of $E$ over the basepoint $z_0$).
Solving $d+N$ numerical ODE problems thus gives a collection of subdominant solution vectors that can be directly substituted into the determinantal invariants.

Invariants of the connections $\nabla^\higgs_w(R,\zeta)$ for the Hitchin section can also be computed numerically by this approach, but an important additional complication in this case is that 
the formula \eqref{eq:hitchin-section-family} for $\nabla^\higgs_w(R,\zeta)$ involves the undetermined function $u$.
Thus we first need to determine $u$, which we do by solving the PDE and conditions \eqref{eq:self-duality}-\eqref{eq:compatibility} numerically on a region in the plane.
Once this has been done, the numerical calculation of Stokes data proceeds as in the oper case, though of course we must choose the ODE integration radius $r$ small enough so that the rays lie in the region on which $u$ has been computed.

In what follows we refer to this approach to computing Stokes data using numerical ODE/PDE solvers as the \emph{direct method} or the \emph{differential equation method} (abbreviated DE), especially when contrasting it with the conjectural integral equation method described in \autoref{sec:ieq}.
The preceding overview of the direct method omits many details involved in the actual numerical implementation used in our experiments, which are discussed in \autoref{subsec:implementation-opers} (opers) and \autoref{subsec:implementation-hitchin} (Hitchin section).

\subsection{The \hk metric} \label{sec:hk-metric}

Hitchin introduced a complete \hk metric on the moduli space of Higgs bundles over a compact
Riemann surface \cite{Hitchin1987}.
An analogous picture holds for Higgs bundles on $\C\PP^1$ 
with irregular singularity at $z = \infty$ \cite{MR1864833,wnh,MR3931781}.
We now briefly review the main facts which we will need below.

Let $\cB_{d,0} \subset \cB_d$ denote the space of tuples $(\phi_2, \dots, \phi_N)$
with $\deg(P_k) < \frac{k-1}{N}d - 1$ for all $k$.
Then for each $h \in \cB_{d} / \cB_{d,0}$ we consider the affine space
\begin{equation}
\cB_{d,h} = h + \cB_{d,0} \quad \subset \quad \cB_d.
\end{equation}
The Hitchin section described in \autoref{subsec:hitchin} realizes 
each $\cB_{d,h}$ as a subspace of
a moduli space $\cM_{d,h}$ of polynomial Higgs 
bundles. The space $\cM_{d,h}$ carries a complete \hk metric, and
by restriction one gets a canonical \kahler metric on each $\cB_{d,h}$.

The only case we will use explicitly in this paper is the case $N=2$.
In this case, given a polynomial $h(z)$ of degree $d$, 
$\cB_{d,h} \subset \cB_d$ is the affine space of 
polynomials $P_2(z)$ of the form
\begin{equation}
  P_2(z) = h(z) + l(z)
\end{equation}
where $\deg l < \frac{d}{2} - 1$. Thus
$\dim \cB_{d,h} = \lceil \frac{d}{2} - 1 \rceil$.

In this 
case we can write a concrete formula for the \kahler metric on $\cB_{d,h}$, 
as follows.
Fix a polynomial $P_2 \in \cB_{d,h}$.
Let $u$ be the solution of \eqref{eq:self-duality} for this $P_2$.
Next,
consider a tangent vector to $\cB_{d,h}$ at $P_2$; such a tangent vector is 
represented by a polynomial $\dot{P}_2$ with $\deg \dot{P}_2 < \frac{d}{2} - 1$.
Let $F$ denote the unique bounded complex function in the plane obeying the \ti{complex variation equation}:
\begin{equation} \label{eq:diff-F}
\left(\Delta - 8(\e^{2 u} + \e^{-2u} \abs{P_2}^2)\right) F + 8 \e^{-2u} \overline{P_2} \dot{{P_2}} = 0.
\end{equation}
Then the norm of the tangent vector $\dot{P}_2$ is
\begin{equation}
\label{eq:l2norm-combined}
  \norm{\dot P_2}^2 = \int_\C 4 \e^{-2u} \left(\abs{\dot P_2}^2 - \re (F P_2 \dot{\overline{P_2}})\right) \de x \de y.
\end{equation}
(For large $\abs{z}$ one has
$u \sim \frac12 \log \abs{P_2}$ and $F \sim \frac12 \frac{\dot P_2}{P_2}$,
so that the integrand $\cI$ in \eqref{eq:l2norm-combined} scales like $\abs{\dot{P}_2^2 / P_2}$,
i.e. like $\abs{z}^{2d'-d}$ if $\dot{P_2}$ has degree $d'$; thus the fact that
we imposed $d' < \frac{d}{2}-1$ is just what is needed to ensure
the integral \eqref{eq:l2norm-combined} converges.)

The integrand in the formula \eqref{eq:l2norm-combined} is the same as the integrand for 
the metric on the Hitchin section of the $\SL_2\C$-character variety of a 
compact surface. The computation leading to \eqref{eq:l2norm-combined} in that context
is given in \cite{Dumas2018}; the same computation applies in the present case.

\subsection{The semiflat approximation} \label{sec:semiflat}

In this section we only consider $N=2$.
Suppose that the polynomial $P_2$ has only simple zeros, and consider a rescaling
\begin{equation}
P_2 \to t^2 P_2, \qquad t \in \R_+.
\end{equation}
Then in the limit of large $t$ we expect a kind of ``concentration'' phenomenon:
away from small discs around the zeros of $P_2$,
which shrink to zero size as $t \to \infty$, 
we should have $u \approx u^\rmsf := \frac12 \log \abs{P_2}$
and $F \approx F^\rmsf := \frac12 \frac{\dot{P_2}}{P_2}$.
(One reason for this expectation is that an analogous concentration phenomenon on 
a compact surface was established in \cite{Mazzeo2014,Dumas2018}; 
the only difference in our case
is that instead of a compact surface we are working on the plane with 
a growth condition as $\abs{z} \to \infty$.)
This leads to a scheme for approximating the $L^2$ metric on $\cB_{d,h}$ without
solving any PDEs: we just replace $u \to u^\rmsf$ and $F \to F^\rmsf$ in 
\eqref{eq:l2norm-combined}, which leads to
\begin{equation} \label{eq:semiflat-metric}
  \norm{\dot P_2}^2 \approx \int_\C 2 \frac{ \abs{\dot P_2}^2 } {\abs{P_2}} \de x \de y.
\end{equation}
This is the \ti{semiflat approximation} to the $L^2$ metric.

As we discuss in \autoref{sec:leading},
the conjectures of \cite{Gaiotto:2008cd,Gaiotto:2009hg}
predict that the semiflat approximation is asymptotically close to the actual metric
\eqref{eq:l2norm-combined}:
the difference between the two decays exponentially in $t$.
For the $\SL_2\C$-character variety of a compact surface $C$, the analogous
statement is known to be true \cite{Mazzeo2014,Dumas2018,Fredrickson2020}.

\subsection{The conformal limit} \label{sec:conformal-limit}

In this section we have been discussing two different families of connections
associated to a point $w \in \cB$:
$\nabla^\oper_w(\hbar)$ and $\nabla^\higgs_w(R,\zeta)$.
It is expected that the family $\nabla^\higgs_w(R,\zeta)$
reduces to the simpler family $\nabla^\oper_w(\hbar)$ in a double scaling limit
known as the ``conformal limit'':
\begin{equation} \label{eq:conformal-limit}
	\lim_{R \to 0} [\nabla^\higgs_w(R, \zeta = R \hbar)] = [\nabla^\oper_w(\hbar)].
\end{equation}
Here $[\nabla]$ means the equivalence class of the connection;
so \eqref{eq:conformal-limit} does not say the actual connections have a limit, 
but that their
equivalence classes do (and thus their Stokes data do, since the
Stokes data depend on the connection only up to equivalence).
This relation was proposed in \cite{Gaiotto2014}; in our context of
polynomial Higgs bundles it has not been proven, but in the case of
Higgs bundles on a compact surface it was proven in \cite{Dumitrescu2016,collier2018conformal}.

\section{Integral equations for Stokes data}
\label{sec:ieq}

In this section we recall the twistor Riemann-Hilbert conjecture, 
a conjectural method for computing the Stokes data of
the families of connections $\nabla_w^\oper(\hbar)$ and $\nabla_w^\higgs(R,\zeta)$.
This method computes \emph{spectral coordinates} for the connections, 
which are invariants related to the determinants considered earlier.

\subsection{Spectral coordinates for \texorpdfstring{$N=2$}{N=2}}

We first consider the case $N=2$, so that $w = \phi_2 = P_2(z) \de z^2$.
We will treat the cases of opers and the Hitchin section in parallel.
For opers $\nabla_w^\oper(\hbar)$ we let $\vartheta = \arg \hbar$;
for the connections $\nabla_w^\higgs(R,\zeta)$ we let $\vartheta = \arg \zeta$.

We define the \emph{$\vartheta$-foliation} to consist of the integral curves of the line field $\ker \im e^{-\I \vartheta} \sqrt{-\phi_2}$.
This ``foliation'' has singularities at the zeros of $\phi_2$; specifically, each simple zero of $\phi_2$ is a $3$-pronged singularity.
Thus a leaf can be bi-infinite, may end at a singularity, or may be a segment between singularities.
A segment of the latter type, occurring in the $\vartheta$-foliation, 
is known as a \emph{saddle connection with angle $\vartheta$}.
When we refer simply to a ``saddle connection'' we mean a saddle connection 
with \ti{any} angle.

We assume from now on (as happens in the examples we study) that $\phi_2$ has only simple zeros.
For such $\phi_2$, an angle $\vartheta$ is \emph{BPS-free} if there are no saddle connections with angle
$\vartheta$; otherwise $\vartheta$ is \emph{BPS-ful}.

The union of the leaves of the $\vartheta$-foliation emanating from the zeros of $\phi_2$ is the \emph{$\vartheta$-critical graph}.
For BPS-free $\vartheta$, the critical graph is simply a union of $3 \deg P_2$ half-infinite leaves. 
Each of these leaves goes to infinity in an asymptotic direction lying in the middle of one of the $2 + \deg P_2$ Stokes sectors for the connection $\nabla^\oper_w(\hbar)$ or $\nabla^\higgs_w(R,\zeta)$.
(Several leaves may go to a single Stokes sector.)

The critical graph divides the plane into a collection of foliated strips and half-planes.
The configuration of these strips and half-planes is naturally encoded in a triangulated polygon $T$, depending on $\vartheta$, that is defined as follows:
The vertex set is the collection of Stokes sectors of $\nabla^\oper_w$ at infinity, and there is an edge from sector $i$ to sector $j$ if there is a bi-infinite $\vartheta$-leaf (and hence strip or half-plane) of $\phi_2$ asymptotic to both $i$ and $j$.
In this triangulated polygon, each triangle naturally corresponds to a simple zero of $\phi_2$, such that the vertices of the triangle are the sectors to which the leaves emanating from that zero are asymptotic.
Figures \ref{fig:A1A2tri} and \ref{fig:A1A3tri} show examples of critical graphs and corresponding triangulations.

The \emph{spectral curve} $Y$ is the holomorphic curve in $\C^2$ defined by
\begin{equation}
  Y = \{ (y,z) \st y^2 + P_2(z) = 0 \}.
\end{equation}
The $1$-form $\lambda = y \, \de z$ on $Y$ satisfies $\pi^*(\phi_2) = -\lambda^2$ where $\pi(y,z) = z$.
Let $\Gamma = H_1(Y,\Z)$, which is the \emph{charge lattice}.

Provided that $\vartheta$ is BPS-free, each saddle connection $e$ gives rise to an element of $\Gamma$ as follows.
Let $\vec{e}_1,\vec{e}_2$ denote the two segments in $Y$ that project isomorphically to $e$, each oriented so that $\im(\lambda_\vartheta)$ is negative when applied to the tangent vector.
Then $\vec{e}_1 + \vec{e}_2$ is a cycle on $Y$, and its homology class $[\vec{e}_1+ \vec{e}_2] \in \Gamma$ is called the \emph{$\vartheta$-lift of $e$}.
It is not hard to see from the explicit form of the spectral curve that the resulting map from the set of saddle connections to $\Gamma$ is injective.

As a particular case of this construction, if $d$ is an internal edge of the triangulated polygon for a BPS-free angle $\vartheta$, then there is a saddle connection $e_d$ naturally dual to $d$, in the sense that the two triangles adjacent to $d$ correspond to the two zeros joined by $e_d$.
Let $\gamma \in \Gamma$ be the $\vartheta$-lift of $e_d$. We define the associated \emph{spectral coordinate} as
\begin{equation}
  X_\gamma = - \frac{p(q,r)p(s,t)}{p(q,t)p(s,r)},
\end{equation}
where $q,r,s,t$ are the vertices of the quadrilateral of which $d$ is the diagonal, with $d$ joining $q$ to $s$.
That is, $X_\gamma$ is a certain cross ratio of the four subdominant solutions associated to triangles that have $d$ as an edge. (This association of a cross ratio to a triangulation
of the polygon was introduced in \cite{MR2233852}, where it was used to define a
cluster structure on an appropriate moduli space of $\PSL_2\C$-local systems.)

There is an extension of this construction, described in \cite{Gaiotto:2009hg}, whereby a spectral coordinate $X_\gamma$ can be associated to \emph{every} element of $\Gamma$, and so that the resulting coordinates satisfy
\begin{equation}
\label{eqn:functional}
  X_\alpha X_\beta = X_{\alpha + \beta}
\end{equation}
for any $\alpha,\beta \in \Gamma$.
As special cases we have that $X_0=1$ and $X_{-\gamma} = X_\gamma^{-1}$.
In the examples we study, a basis of $\Gamma$ is obtained from internal edges of $T$, hence the equation above actually determines a formula for any spectral coordinate in terms of the ones arising from internal edges.

\subsection{Integral equation for opers} \label{sec:ieq-opers}

Consider a fixed $w \in \cB'$ and varying $\hbar \in \C^\times$. This gives a 1-parameter
family of spectral coordinates $X_\gamma(\hbar)$ associated to the connections $\nabla_w^\oper(\hbar)$.
The twistor Riemann-Hilbert conjecture says that $X_\gamma(\hbar)$ can be
computed by another method, which we now describe. We will give a description of a family of
functions $\cX_\gamma(\hbar)$, which \ti{a priori} has nothing to do with flat connections; 
the conjecture is that
\begin{equation} \label{eq:gmn-conjecture-opers}
X_\gamma(\hbar) = \cX_\gamma(\hbar).
\end{equation}

The functions $\cX_\gamma$ are characterized in terms of a system of integral equations.
To state these we will need to define periods and BPS counts.
A class $\gamma \in \Gamma$ has a \emph{period} $Z_\gamma \in \C$ defined by
\begin{equation} Z_\gamma = \int_\gamma \lambda. \end{equation}
The theory of \cite{Gaiotto2012} associates to $\gamma \in \Gamma$ an integer $\Omega(w,\gamma)$, the \emph{BPS count}. In these $N=2$ examples $\Omega(w,\gamma)$ is simply given by
\begin{equation}
  \Omega(w, \gamma) = \begin{cases}
    1 & \text{if } \pm \gamma \text{ is associated to a saddle connection,}\\
    0 & \text{otherwise.}
  \end{cases}
\end{equation}
Note that $\Omega(w,\gamma) = \Omega(w,-\gamma)$.
Let $\Gamma' \subset \Gamma$ denote the set of homology classes $\gamma$ for which $\Omega(w, \gamma) \neq 0$.
In the examples we consider, $\Gamma'$ is a finite set.

Part of the twistor Riemann-Hilbert conjecture is the statement that a 
set of functions $\{ \cX_\gamma \st \gamma \in \Gamma \}$ can be uniquely determined by the system of integral equations
\begin{equation}
  \label{eq:integral-oper}
  \cX_\gamma(\hbar) = \exp \left ( \hbar^{-1} Z_\gamma + \frac{1}{4 \pi \I} \sum_{\mu \in \Gamma} \Omega(w,\mu) \langle \gamma, \mu \rangle \int_{\R_- Z_\mu} \frac{\de \xi}{\xi} \frac{\xi + \hbar}{\xi - \hbar} \log\left (1 + \cX_{\mu}(\xi) \right ) \right )
\end{equation}
where $\langle \gamma, \mu \rangle$ denotes the intersection pairing on $\Gamma = H_1(Y,\Z)$.
Note that the formally infinite sum in \eqref{eq:integral-oper} has only finitely many nonzero terms, because it includes a coefficient $\Omega(w,\mu)$ and so can be reduced to a sum over $\mu \in \Gamma'$.
Thus by considering only $\gamma \in \Gamma'$ we obtain from \eqref{eq:integral-oper} a \emph{finite} collection of coupled integral equations.

Assuming this conjecture holds, it suggests a method to 
compute the collection $\cX_\bullet = \{ \cX_\gamma \st \gamma \in \Gamma' \}$:
Let $\cF$ denote the right hand side of \eqref{eq:integral-oper}, considered as a self-map of the set of tuples of functions of $\hbar$ indexed by $\gamma \in \Gamma'$.
In terms of this function, the conjecture is that $\cX_\bullet$ is a fixed point of $\cF$, i.e.~that $\cX_\bullet = \cF( \cX_\bullet)$.
We can further optimistically conjecture that this fixed point is unique, and then
attempt to find it by iteration, starting with the initial guess $\cX_\gamma^{(0)}(\hbar) = \exp(\hbar^{-1} Z_\gamma)$ and inductively defining $\cX_\bullet^{(k)} = \cF(\cX_\bullet^{(k-1)})$ for any $k > 0$ --- or some
similar iteration with the same fixed points (see \autoref{subsec:implementation-ieq-opers} for the precise iteration we use in practice).
Finally, once $\cX_\bullet$ has been determined, we 
can easily compute $\cX_\gamma$ for any $\gamma \in \Gamma$
if desired, using \eqref{eqn:functional}.

We refer to this iterative process as the \emph{integral equation method}, which we sometimes abbreviate $\IEQ$ in tables and plots.

\subsection{Integral equation for the Hitchin section}
\label{subsubsec:inteq-hitchin}

In the previous subsection we discussed a conjectural integral equation method for computing 
the Stokes data of the family $\nabla^\oper_w(\hbar)$, for a fixed $w$.
There is a very similar conjecture for the Stokes data of the family $\nabla^\higgs_w(R,\zeta)$,
for a fixed $w$ and $R$.
Instead of \eqref{eq:integral-oper} we consider
\begin{equation}
  \label{eq:integral-hitchin-section}
  \cX_\gamma(R,\zeta) = \exp \left ( R \zeta^{-1} Z_\gamma + R \zeta \overline{Z}_\gamma +  \frac{1}{4 \pi \I} \sum_{\mu \in \Gamma} \Omega(w,\mu) \langle \gamma, \mu \rangle \int_{\R_- Z_\mu} \frac{\de \xi}{\xi} \frac{\xi + \zeta}{\xi - \zeta} \log\left (1 + \cX_{\mu}(R,\xi) \right ) \right ).
\end{equation}
The data $\Gamma$, $Z$, $\Omega$ entering \eqref{eq:integral-hitchin-section}
are exactly as they were for \eqref{eq:integral-oper}, 
so all of the discussion from \autoref{sec:ieq-opers} carries over intact
to this case. Indeed,
the only differences between \eqref{eq:integral-oper} and \eqref{eq:integral-hitchin-section} are:
\begin{itemize}
  \item in \eqref{eq:integral-hitchin-section} the variable is called $\zeta \in \C^\times$ instead of $\hbar \in \C^\times$, and is rescaled by a factor $R$ in some places,
  \item \eqref{eq:integral-hitchin-section} includes the extra term $R \zeta \overline{Z}_\gamma$ which is not present in \eqref{eq:integral-oper}.
  \end{itemize}
Also in parallel to \autoref{sec:ieq-opers}, we can attempt to produce 
functions $\cX_\gamma(R, \zeta)$ obeying
\eqref{eq:integral-hitchin-section} by iteration, starting from
the initial functions 
$\cX_\gamma^{(0)}(R, \zeta) = \exp \left(R\zeta^{-1}Z_\gamma + R \zeta \overline{Z}_\gamma \right)$.
Then the twistor Riemann-Hilbert conjecture is that this iteration 
converges and the spectral coordinates
of $\nabla^\higgs_w(R,\zeta)$ are given by
\begin{equation} \label{eq:gmn-conjecture-hitchin}
 X_\gamma(R, \zeta) = \cX_\gamma(R, \zeta).
\end{equation}

Incidentally, as observed in \cite{Gaiotto2014},
one can obtain \eqref{eq:integral-oper} 
from \eqref{eq:integral-hitchin-section} by performing the 
conformal limit $\zeta = R \hbar$, $R \to 0$, as in \autoref{sec:conformal-limit}.
In this sense the three conjectures we have discussed ---
the integral equation for opers, the integral equation for the Hitchin
section, and the conformal limit --- are compatible.

\subsection{Spectral coordinates and integral equations for \texorpdfstring{$N=3$}{N=3}}

Compared to the $N=2$ case outlined above, there are just a few differences in the predictions for $N=3$.
Since the predictions for $N=3$ are developed carefully in \cite{Neitzke17}, we omit some details here and refer to that paper for additional discussion.

In this case the spectral curve is the $3$-fold branched cover of $\C$ defined by
\begin{equation}
Y = \{  (y,z) \st  y^3 + y P_2(z) + P_3(z) = 0 \}.
\end{equation}
As before, the charge lattice $\Gamma = H_1(Y,\Z)$.
Any $\gamma \in \Gamma$ has a \emph{period} $Z_\gamma$ which is the integral of $\lambda = y \, \de z$ over that cycle.
Our experimental studies focus primarily on the \emph{cyclic} case $w = (0,\phi_3) \in \cB_d$ (i.e.~vanishing quadratic differential) in which case $\pi^*(\phi_3) = -\lambda^3$.

The rank-$3$ analogue of the critical graph is the \emph{WKB} $\vartheta$-\emph{spectral network}, a graph embedded in the plane with edges labeled by certain topological data.
We describe it in detail only in the cyclic case, and when $P_3$ has only simple zeros.
First, a $\vartheta$-trajectory of $\phi_3$ (or $w$) is an oriented curve $z(t)$ in $\C$ equipped with a pair of continuous sections $y(z), \hat{y}(z)$ of $\pi$ over the curve satisfying
\begin{equation}
  \left ( y(t) - \hat{y}(t) \right ) \frac{\de z(t)}{\de t} = \e^{\I \vartheta}.
\end{equation}
Of course, if a local labeling of the sheets of the spectral curve as $y_1,y_2,y_3$ has been given, then in this region and for some $i,j$ the functions $y$ and $\hat{y}$ are restrictions of $y_i$ and $y_j$, respectively, and we can label the trajectory according to its \emph{type} $(i,j)$.
As a global labeling of this type is generally not possible, it is necessary to introduce branch cuts that divide the plane into simply connected regions and indicate labels in each region, as well as the permutation of labels when crossing the branch cut.
Also note that reversing orientation of a trajectory of type $(i,j)$ gives a trajectory of type $(j,i)$.

A simple zero of $P_3$ has eight $\vartheta$-trajectories emerging from it, and the $\vartheta$-spectral network is defined to include the maximal extensions of these trajectories.
We also add additional curves to the spectral network:
If trajectories with local labels $(i,j)$ and $(j,k)$ meet (at $p$), then they do so at angle $\frac{\pi}{3}$ or $\frac{2\pi}{3}$.
In the latter case, there is a trajectory of type $(i,j)$ bisecting the angle between their tangent vectors at $p$, and we add the maximal extension of this trajectory to the network as well.
This may result in new intersections, and then we repeat the rule above to possibly add additional trajectories from the intersection points.
The $\vartheta$-spectral network is the union of all trajectories that arise from iterating this procedure.
In the examples we consider, it is a \emph{finite} union of trajectories.

We say that $\vartheta$ is BPS-free if no trajectory of the $\vartheta$-spectral network meets a zero of $P_3$, except possibly at its origin point; otherwise $\vartheta$ is BPS-ful.

As in rank $2$ there is a procedure which associates to $\gamma \in \Gamma$ a spectral coordinate $X_\gamma$ which is given by some combination of determinants of subdominant solutions.
The general procedure to construct this mapping from homology classes to coordinate functions is significantly more complicated than for $N=2$, and we will not describe it.
Instead, we will indicate the result of that procedure in the examples we study (referring to \cite{Neitzke17} for both the general procedure and the detailed calculations in these examples).

The functions $\cX_\gamma$ and the conjectural integral equation are defined exactly the same as before (i.e.~\eqref{eq:integral-oper} or \eqref{eq:integral-hitchin-section}), though to make sense of this equation in rank $3$ we must describe the meaning of $\Omega(w,\mu)$.
Recall that in rank $2$ we defined $\Omega(w,\mu)$ to be $1$ or $0$ depending on whether $\pm \mu$ is represented by a saddle connection.
In the cyclic rank $3$ case, there are cycles in $\Gamma$ associated to trajectories joining zeros (essentially, rank-$3$ saddle connections) \emph{and} to tripods consisting of three trajectories joining three zeros to a common point (with labels $(i,j)$, $(j,k)$, $(k,i)$).
In the examples we study\footnote{In more general rank $3$ cases it would be necessary to consider other types of ``degenerations'' of the spectral network; an algorithm for defining and computing $\Omega(w,\mu)$ in general is given in \cite{Gaiotto2012}.}, the BPS count $\Omega(w,\mu)$ is the total number of such saddle connections and tripods whose associated cycle is $\pm \mu$.
As in rank $2$ we will only study examples where $\Omega(w,\mu)=0$ for all but finitely many $\mu \in \Gamma$, allowing the same integral iteration strategy sketched above to extend naturally to rank $3$.

Finally, we allow parameterized deformations of a cyclic example, i.e.~$w = (\phi_2(\Lambda), \phi_3)$ with $\Lambda$ a complex parameter and $\phi_2(\Lambda=0)=0$.
For sufficiently small $\abs{\Lambda}$ it is not necessary to fully generalize the spectral network and spectral coordinate constructions to this case, as all of the relevant choices are locally constant: a homology basis for $\Lambda=0$ extends to the spectral curves over a neighborhood of $\Lambda=0$ using the Gauss-Manin connection, and the associated determinantal invariants also remain unchanged.

\subsection{Spectral coordinates and the \hk metric} \label{sec:spectral-hitchin}

In this section we specialize to $N=2$
and revisit the \kahler metric $g_{d,h}$ 
on the space $\cB_{d,h}$ introduced in \autoref{sec:hk-metric}.

In parallel to the usual case of Higgs bundles over a compact surface,
the metric $g_{d,h}$ is the
restriction of the \hk metric on $\cM_{d,h}$.
On a \hk manifold one has three distinguished complex structures $I_1$, $I_2$, $I_3$,
and corresponding \kahler forms $\omega_1$, $\omega_2$, $\omega_3$.
In our case, $\cB_{d,h}$ is an $I_3$-complex subspace
of $\cM_{d,h}$, and thus the restriction of $\omega_3$ to $\cB_{d,h}$
is the \kahler form for the metric $g_{d,h}$.

On any \hk manifold, the \kahler forms for any two of the 
distinguished complex 
structures can be organized into a holomorphic-symplectic form
for the third; in particular, the $I_2$-holomorphic-symplectic
form is 
\begin{equation} \label{eq:holsymp}
\Omega_2 = \omega_3 + \I \omega_1.
\end{equation}
The key reason why the spectral coordinates can be used to
understand the \hk metric on $\cM_{d,h}$ is that 
the functions
\begin{equation}
y_\gamma := X_\gamma(R=1,\zeta=1)
\end{equation}
are Darboux coordinates for the holomorphic-symplectic form $\Omega_2$ \cite{MR2233852,Gaiotto:2009hg,Gaiotto2012}. More precisely,
choosing a basis $\{\gamma_i\}$ for $\Gamma$
and setting $\eps_{ij} := \IP{\gamma_i, \gamma_j}$, $y_i := y_{\gamma_i}$,
we have
\begin{equation} \label{eq:hol-darboux}
  \Omega_2 = \half \sum_{i,j = 1}^n \eps_{ij} \, \de \log y_i \wedge \de \log y_j.
\end{equation}
When restricted to $\cB_{d,h}$ the functions $y_\gamma$
are real, so combining \eqref{eq:hol-darboux} with \eqref{eq:holsymp} gives on $\cB_{d,h}$
\begin{equation} \label{eq:kahler-spectral}
  \omega_3 = \half \sum_{i,j = 1}^n \eps_{ij} \, \de \log y_i \wedge \de \log y_j.
\end{equation}
Finally, we have
$g(\cdot, \cdot) = \omega_3(\cdot, I_3 \cdot)$, and thus
\begin{equation} \label{eq:metric-spectral}
  g = \sum_{i,j = 1}^n \eps_{ij} \, \de \log y_i \otimes \de^c \log y_j,
\end{equation}
where $\de^c f \in \Omega^1(\cB_{d,h})$ is defined by $\de^c f \cdot v = \de f \cdot I_3 v$.

Using \eqref{eq:metric-spectral}, any method of computing the functions $y_\gamma$ on $\cB_{d,h}$
gives a method of computing the metric $g_{d,h}$. In particular, we can use the
integral equation method from \autoref{subsubsec:inteq-hitchin}
to compute $y_\gamma$, and thus obtain an integral 
equation computation of $g_{d,h}$.

\subsection{Leading-order approximations} \label{sec:leading}

In this section we have been discussing a method of computing spectral coordinates 
$X_\gamma$ for the connections $\nabla^\oper_w(\hbar)$ or $\nabla^\higgs_w(R,\zeta)$, 
which involves solving either \eqref{eq:integral-oper} or \eqref{eq:integral-hitchin-section}
respectively. Although the full solutions $\cX_\gamma(\hbar)$ or $\cX_\gamma(R,\zeta)$
are complicated and do not seem to admit explicit exact formulas, we can nevertheless derive
explicit \ti{asymptotic} formulas.
The details are slightly different in the two cases:

\begin{itemize}
\item In the case of opers, if one makes the assumption
that $\cX_\gamma(\hbar) \to 0$ as $\hbar \to 0$ along the axis
$\hbar \in \R_- Z_\gamma$, then it follows that the sum of integrals in 
\eqref{eq:integral-oper} goes to zero as $\hbar \to 0$,\footnote{To verify this it is convenient to rewrite the sum of integrals by combining
the terms for $\mu$ and $-\mu$, as we do
in \autoref{subsec:implementation-ieq-opers} below.}
and thus
\begin{equation} \label{eq:wkb-asymp}
  \cX_\gamma(\hbar) \sim \exp(\hbar^{-1} Z_\gamma) \qquad \text{ as } \hbar \to 0.
\end{equation}

\item In the case of the Hitchin section, one can 
make a sharper asymptotic statement, as follows. 
If we take $\zeta \in \R_- Z_\gamma$,
the sum of integrals in 
\eqref{eq:integral-hitchin-section} is real and negative; this
implies that, for $\zeta \in \R_- Z_\gamma$, $\abs{\cX_\gamma(\zeta)}$
is bounded above by $\exp(-2 R \abs{Z_\gamma})$. This in turn implies that 
the integral term in \eqref{eq:integral-hitchin-section} 
is $O(\e^{-2 R M})$ where $M = \min \{ \abs{Z_\mu}: \mu \in \Gamma' \}$. 
We conclude that
\begin{equation} \label{eq:sf-asymp-hitchin}
  \cX_\gamma(R, \zeta) = \exp \left( R \zeta^{-1} Z_\gamma + R \zeta \overline{Z}_\gamma + O(\e^{-2RM}) \right) \qquad \text{ as } R \to \infty.
\end{equation}
In particular, at $\zeta = 1$ we have
\begin{equation} \label{eq:sf-asymp-hitchin-zeta1}
  \cX_\gamma(R, \zeta=1) = \exp \left( 2 R \re Z_\gamma + O(\e^{-2RM}) \right) \qquad \text{ as } R \to \infty.
\end{equation}

\end{itemize}

The asymptotic formula \eqref{eq:sf-asymp-hitchin-zeta1}
also leads to an asymptotic
formula for the \kahler metric on $\cB_{d,h}$, as follows.
We consider a ray on $\cB_{d,h}$ given by $w = t w_0$, in the sense 
of \eqref{eq:scaling-differentials}, and the metric along this ray, 
$g_t = g(w = tw_0)$.
$g_t$ is determined by the Darboux coordinates $x_\gamma(t w_0) = X_\gamma(t w_0, R=1, \zeta = 1) = X_\gamma(w_0, R = t, \zeta = 1)$.
Then
\begin{equation}
y_\gamma(t w_0) = \exp \left( 2 t \re Z_\gamma(w_0) + O(\e^{-2 t M}) \right) \qquad \text{as } t \to \infty.
\end{equation}
Using the relation $Z_\gamma(t w_0) = t Z_\gamma(w_0)$, we can also write this
\begin{equation}
y_\gamma(t w_0) = \exp \left( 2 \re Z_\gamma(t w_0) + O(\e^{-2 t M}) \right) \qquad \text{as } t \to \infty.
\end{equation}
Finally, substituting this in \eqref{eq:metric-spectral} gives the approximate formula
\begin{equation}
  g_t = g^\rmsf + O(\e^{-2t M}) \qquad \text{as } t \to \infty,
\end{equation}
where we define
\begin{equation} \label{eq:semiflat-spectral}
  g^\rmsf = \sum_{i,j = 1}^n \eps_{ij} \, \de \re Z_i \otimes \de^c \re Z_j.
\end{equation}

We remark that
the metric $g^\rmsf$ defined here is actually equal to the semiflat approximation
to $g$ which we described in \autoref{sec:hk-metric}.
Indeed \eqref{eq:semiflat-spectral} expresses $g^\rmsf$ in terms of a bilinear
form in periods $\dot{Z}_\gamma = \oint_\gamma \dot{\lambda}$ and their complex conjugates, 
which by Riemann bilinear identity
can be related to the integral $\int_Y \dot\lambda \wedge \dot{\bar\lambda}$;
likewise \eqref{eq:semiflat-metric} can be related to the same integral.
Thus the twistor Riemann-Hilbert conjecture leads to the prediction that
the semiflat approximation is exponentially good: it holds
up to corrections of order $\e^{-2 t M}$.
With this connection in mind we refer to the asymptotic formula \eqref{eq:sf-asymp-hitchin}
as the semiflat approximation to $\cX_\gamma(\zeta)$.

\subsection{Exact coordinates for pure flavor charges} \label{sec:pure-flavor}

In general it is difficult to write explicit exact 
formulas for the functions $\cX_\gamma$ obeying \eqref{eq:integral-oper} or \eqref{eq:integral-hitchin-section}. There is one exception, however:
this is the case when 
$\gamma$ lies in the radical of the intersection pairing $\IP{\cdot,\cdot}$.
(Such $\gamma$ are called ``pure flavor charges'' in the physics literature;
geometrically, in our examples, they arise from cycles on $\Sigma$ which are
peripheral, i.e. they lie in a small neighborhood of $\infty$.)

When $\gamma$ is pure flavor, 
the integral terms in \eqref{eq:integral-oper}, \eqref{eq:integral-hitchin-section} vanish, leaving simply
\begin{equation} \label{eq:X-pure-flavor-oper}
	\cX_\gamma(\hbar) = \exp(\hbar^{-1} Z_\gamma),
\end{equation}
and
\begin{equation} \label{eq:X-pure-flavor-hitchin}
	\cX_\gamma(R,\zeta) = \exp(R \zeta^{-1} Z_\gamma + R \zeta \overline{Z_\gamma}).
\end{equation}
In other words, when $\gamma$ is pure flavor, the asymptotic formulas
\eqref{eq:wkb-asymp}, \eqref{eq:sf-asymp-hitchin} simplify to exact formulas.
Combining these with the twistor Riemann-Hilbert conjecture gives exact 
formulas for $X_{\gamma}(\hbar)$ and $X_{\gamma}(R, \zeta)$.
These formulas are also conjectural, but should be much simpler to establish
than the conjecture for general $\gamma$, and in at least
one case they are already known: 
when $N=2$, the formula $X_\gamma(R=1,\zeta=1) = \exp(2 \re Z_\gamma)$
is proven in \cite{Gupta2017}.

\section{Experimental studies of Stokes data}
\label{sec:experiments}

\subsection{Examples} \label{sec:examples}

\subsubsection{The $N=2$ examples}
\label{subsubsec:N2ex}

In general we will refer to the examples for given $N$ and $d$ using the \emph{theory name} $(A_{N-1},A_{d-1})$, following the notation of \cite{cecotti-neitzke-vafa} for the associated generalized Argyres-Douglas quantum field theory.
For $N=2$ we consider the cases $d=3$ and $d=4$, i.e.~the $(A_1,A_2)$ and $(A_1,A_3)$ theories.

In each case we choose a base point in $\cB_d$ and a BPS-free angle $\vartheta_0$, and introduce a basis $B = \{ \gamma_1, \ldots, \gamma_{d-1} \}$ of $\Gamma$.
The homology calculations for these base points naturally extend to all polynomials in a small neighborhood of the base point (which we also parameterize explicitly, in the cases we study) and for all $\vartheta$ near $\vartheta_0$.
Fixing a homology basis allows us to limit our calculations to the spectral coordinates $X_k := X_{\gamma_k}$, which determine all others using \eqref{eqn:functional}.

In describing our homology bases use the shorthand notation $[x,y]$ to refer to the element of $\Gamma$ that is the $\vartheta_0$-lift of the saddle connection from $x$ to $y$, and $-[x,y]$ for its opposite.
The data describing these examples are summarized in Tables \ref{tab:N2ex}--\ref{tab:N2coef}, and the homology calculations are detailed below.

\begin{table}
  \makebox[\textwidth][c]{
  \begin{tabular}{@{}llllll@{}}
    \toprule
    Theory & Family & Basepoint & $\vartheta_0$ & $\Gamma$-basis & Periods at basepoint\\
    \midrule
    $(A_1,A_2)$ & $P_2 = z^3 - \Lambda z - c$ & $\Lambda=0$ & $0$ & $\gamma_1 = -[1,\omega]$ & $Z_{\gamma_1} = e^{5 \pi \I /6} M, \;\; M=\sqrt{3 \pi} \frac{\Gamma(4/3)}{\Gamma(11/6)}$ \\
           & & $c=1$ & & & $\quad\;\;\: \approx -2.52393 + 1.45719\I$ \\
           & & & & $\gamma_2 = [\omega,\omega^2]$ & $Z_{\gamma_2} = -\I M$\\
           & & & & & $\quad\;\;\: \approx -2.91438\I$ \\
    \midrule
    $(A_1,A_3)$ & $P_2 = z^4 - 1$ & --- &  $0.4$ & $\gamma_1 = [-1,1]$ & $Z_{\gamma_1} = 2 \sqrt{\pi} \frac{\Gamma(5/4)}{\Gamma(7/4)}$\\
           & & & & & $\quad\;\;\: \approx 3.49608$ \\
           & & & & $\gamma_2 = -[-1,\I]$ & $Z_{\gamma_2} = \frac12(1+\I) Z_{\gamma_1}$ \\
           & & & & & $\quad\;\;\: \approx 1.74804(1+\I)$\\
           & & & & $\gamma_3 = -[1,-\I]$ & $Z_{\gamma_3} = Z_{\gamma_2}$\\
    \bottomrule
  \end{tabular}
  }
  \vspace{0.5em}
  \caption{Summary of data for the $N=2$ examples we study.}
  \label{tab:N2ex}
\end{table}

\begin{table}
  \makebox[\textwidth][c]{
  \begin{tabular}{@{}lll@{}}
    \toprule
    Theory & $\langle \cdot,\cdot \rangle$ & $\Gamma'$\\
    \midrule
    $(A_1,A_2)$ & $\begin{pmatrix}0&1\\-1&0\end{pmatrix}$ & $(1,0)$, $(0,1)$, $(1,1)$\\
    \midrule
    $(A_1,A_3)$ & $\begin{pmatrix}0&1&1\\-1&0&0\\-1&0&0\end{pmatrix}$ & \makecell[l]{$(1,0,0)$, $(0,1,0)$, $(0,0,1)$,\\$(1,0,-1)$, $(1,-1,0)$, $(1,-1,-1)$}\\
    \bottomrule
  \end{tabular}
  }
  \vspace{1em}
  \caption{Intersection form and classes with nonzero BPS counts in the $N=2$ examples, with respect to the bases in Table \ref{tab:N2ex}.  Since $\Gamma' = -\Gamma'$, we list only one from each pair $\pm \gamma$.}
  \label{tab:N2coef}
\end{table}                               

\paragraph{Example $(A_1,A_2)$.}

We take as base point the polynomial $P_2(z) = z^3-1$ and the BPS-free angle $\vartheta_0=0$.
Here the spectral curve is the double cover of $\C$ branched over the third roots of unity $\{1,\omega,\omega^2\}$ (where $\omega = \exp(2 \pi \I/3)$), which as a Riemann surface is the hexagonal punctured torus.
For each pair of roots of unity there is a unique saddle connection joining them, and for $j=1,2,3$ we define $\gamma_j = (-1)^{j} [\omega^{j-1},\omega^j]$.
These cycles satisfy $\gamma_1+\gamma_2+\gamma_3=0$, with any two of them forming a basis for $\Gamma \simeq \Z^2$.
We fix the basis $B = \{\gamma_1, \gamma_2\}$ for our calculations, which has intersection form $\langle \gamma_1, \gamma_2\rangle=1$.
The set $\Gamma'$ (i.e.~the $\mu$ with $\Omega(w,\mu) = 1$) consists of $\pm \gamma_1, \pm \gamma_2, \pm \gamma_3 = \mp(\gamma_1+\gamma_2)$.
The periods of the basis elements are listed in \autoref{tab:N2ex}.

\begin{figure}
  \begin{center}
    \includegraphics[width=0.8\textwidth]{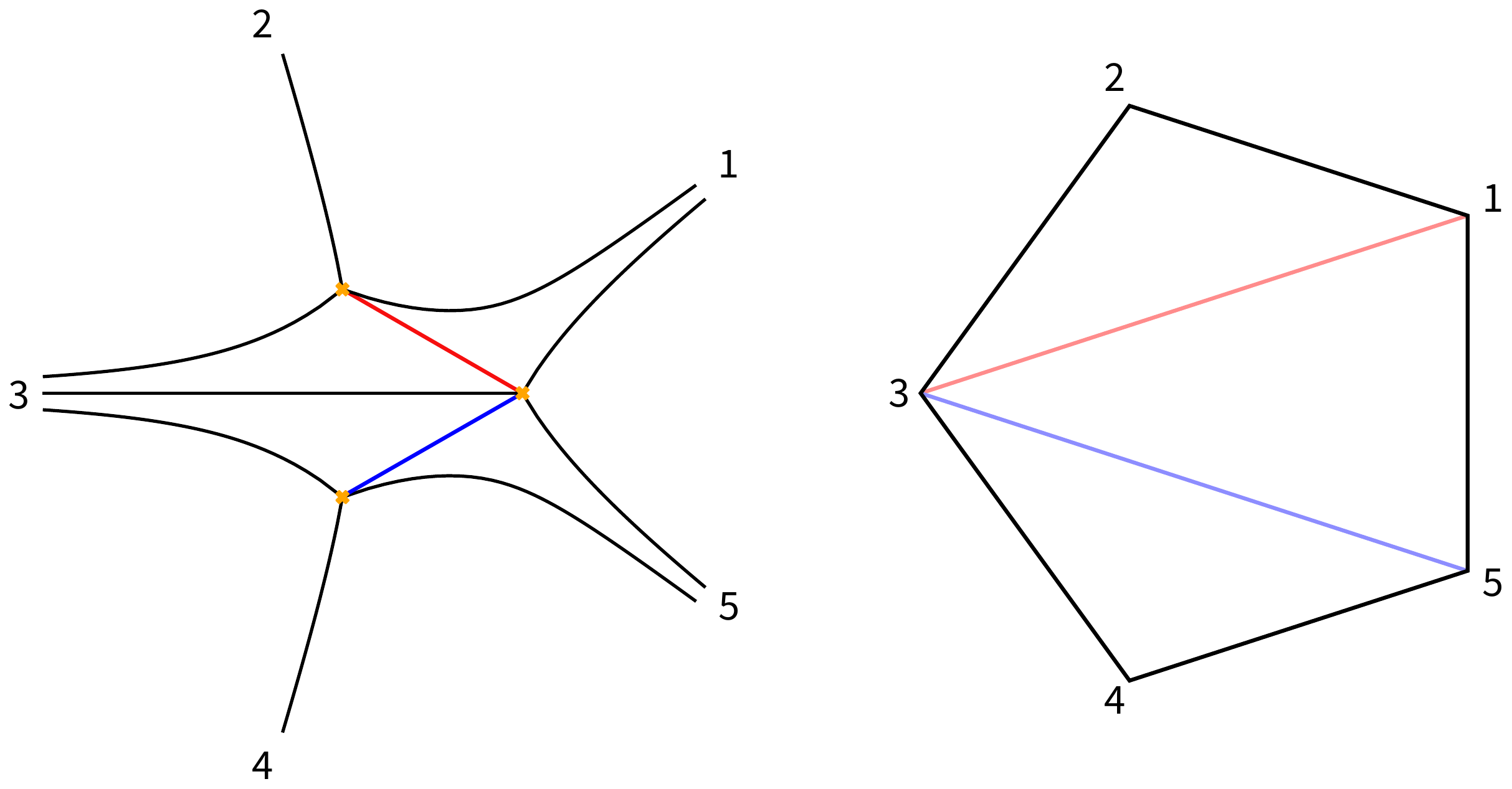}
  \end{center}
  \caption{Spectral network (left) and triangulated polygon $T$ (right) for $(A_1,A_2)$ ($P_2(z) = z^3-1$) at $\vartheta=0$.  The same triangulation arises for all $|\vartheta| < \frac{\pi}{6}$.  The labeled saddle connections correspond to $-\gamma_1$ (red) and $-\gamma_3$ (blue), and the dual edges in $T$ are correspondingly colored.}
  \label{fig:A1A2tri}
\end{figure}

The spectral network and triangulated polygon in this case are shown in \autoref{fig:A1A2tri}.
The cycles $-\gamma_1$ and $-\gamma_3$ are associated to internal edges of $T$, and so the procedure described in the previous section determines their associated spectral coordinates, $X_{-\gamma_1}$ and $X_{-\gamma_3}$.
The relation \eqref{eqn:functional} then determines $X_{\gamma_2} = X_{-\gamma_1} X_{-\gamma_3}$.
Explicitly, for the basis elements this gives
\begin{equation}
    X_1 = \frac{p(2,3) p(1,5)}{p(1,2) p(3,5)}, \;\;\;
   X_2 = \frac{p(2,3) p(4,5) p(1,3)}{p(1,2) p(3,4) p(3,5)}.
\end{equation}

In this example we also study a parameterized family within $\cB_3$ containing the base point by taking $P_2(z) = z^3 - \Lambda z - c$ where $\Lambda$ and $c$ are complex parameters.

\paragraph{Example $(A_1,A_3)$.}

We use $P_2(z) = z^4-1$ as a base point.
Here $\vartheta=0$ is BPS-ful, but angles in the range $(0,\pi/4)$ are all BPS-free and give the same triangulation; for concreteness we choose $\vartheta_0=0.4$.

\begin{figure}
  \begin{center}
    \includegraphics[width=0.8\textwidth]{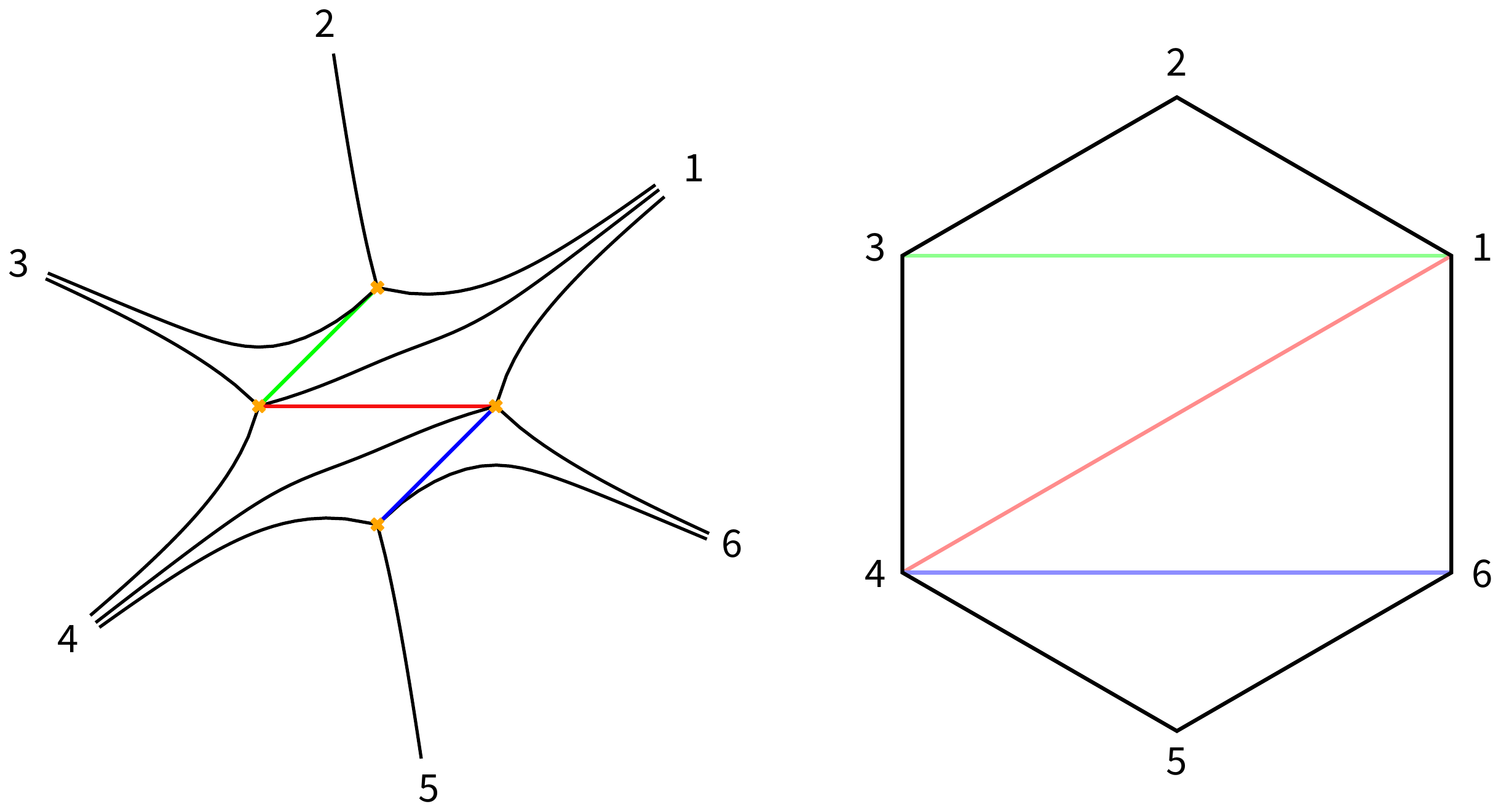}
  \end{center}
  \caption{Spectral network (left) and triangulated polygon (right) for $(A_1,A_3)$ ($P_2(z) = z^4-1$) at $\vartheta=0.4$.  The same triangulation arises for all $0 < \vartheta < \frac{\pi}{4}$. At left, the colored segments are the saddle connections corresponding to $\gamma_1$ (red), $-\gamma_2$ (green), and $-\gamma_3$ (blue), and the dual edges in $T$ are correspondingly colored.}
  \label{fig:A1A3tri}
\end{figure}

The spectral curve is the twice-punctured square torus, with $\Gamma \simeq \Z^3$.
The spectral network and triangulated polygon are shown in \autoref{fig:A1A3tri}.
As in the previous example, for any pair of zeros there is a unique BPS-ful angle giving a saddle connection joining them.
Let $\gamma_1=[-1,1]$, $\gamma_2 = -[-1,\I]$ and $\gamma_3 = -[-\I,1]$, which give a basis of $\Gamma$.
Their periods are shown in Table \ref{tab:N2ex}.

For each of these basis cycles, one of $\pm\gamma_k$ is associated to an internal edge of $T$, and the associated spectral coordinates are the cross ratios:
\begin{equation}
    X_1 = \frac{p(1,3) p(4,6)}{p(1,6) p(3,4)}, \;\;
    X_2 = \frac{p(2,3) p(1,4)}{p(1,2) p(3,4)}, \;\;
    X_3 = \frac{p(1,4) p(5,6)}{p(4,5) p(1,6)}.
\end{equation}

In this example we restrict our study to the base point, and do not consider any parametric family in $\cB_4$.

In this example the cycle $\gamma_f := \gamma_2 - \gamma_3$ lies in the kernel
of the intersection pairing, i.e. it is a pure flavor charge.
Moreover, this cycle has $Z_{\gamma_f} = 0$. Thus from
\eqref{eq:X-pure-flavor-oper}, \eqref{eq:X-pure-flavor-hitchin} 
we see that $\cX_{\gamma_f}(\hbar) = 1$ and $\cX_{\gamma_f}(R, \zeta) = 1$
identically. The corresponding spectral coordinate is
\begin{equation}
	X_{\gamma_f} = X_2 X_3^{-1} = \frac{p(2,3) p(4,5) p(1,6)}{p(1,2)p(3,4)p(5,6)}.
\end{equation}
Thus the twistor Riemann-Hilbert conjecture $\cX_{\gamma} = X_{\gamma}$ 
would imply that this combination is identically
equal to $1$, i.e. that $X_2 = X_3$.

For the particular basepoint we consider, $P_2(z) = z^4 - 1$,
we actually have an extra symmetry under the holomorphic automorphism
$z \mapsto -z$, which implies that indeed $X_2 = X_3$ in this case.
Because of the fact that $X_2 = X_3$, we omit $X_3$ when showing experimental results in this example.
(More generally, we could have considered say $P_2(z) = z^4 + a z + b$;
in this case we do not have the extra symmetry anymore, but we do
still have $Z_{\gamma_f} = 0$, and thus the conjecture would imply that 
$X_2 = X_3$ even in this case.)

\subsubsection{The $N=3$ examples}
\label{subsubsec:N3ex}

Recall that we refer to examples by theory name $(A_{N-1},A_{d-1})$.
For $N=3$ we consider the $(A_2,A_1)$ and $(A_2,A_2)$ theories, and in this section we recall the choices of base points, homology bases, and associated spectral coordinates as computed in \cite{Neitzke17}.
The results of this discussion are summarized in Tables \ref{tab:N3ex}--\ref{tab:N3coef}.

\begin{table}
  \makebox[\textwidth][c]{
  \begin{tabular}{@{}llllll@{}}
    \toprule
    Theory & Family & Basepoint & $\vartheta_0$ & $\Gamma$-basis & Periods\\
    \midrule
    $(A_2,A_1)$ & $P_3 = \frac{1}{2}(1-z^2)$ & $c=0$ & $0$ & (Fig.~\ref{fig:A2A1basis}) & $Z_{\gamma_1} = \e^{5 \pi \I/6} \frac{12 \times 2^{2/3} \times \pi^{3/2}}{5 \Gamma(-1/6) \Gamma(2/3)}$\\
           & $P_2 = c$ & & & & $\quad\;\;\: \approx - 2.00324 + 1.15657\I$\\
           & & & & & $Z_{\gamma_2} = \e^{2 \pi \I/3} Z_{\gamma_1}$\\
           & & & & & $\quad\;\;\: \approx - 2.31315\I$\\
    \midrule
    $(A_2,A_2)$ & $P_3 = \frac{1}{2}(z^3 - 3 z^2 - 2)$ & --- & $0$ & (Fig.~\ref{fig:A2A2basis}) & $Z_{\gamma_1} \approx 2.30298$\\
           & $P_2 = 0$ & & & & $Z_{\gamma_2} \approx 5.47033 + 4.48792\I$\\
           & & & & & $Z_{\gamma_3} \approx - 4.31884 + 2.49348\I$\\
           & & & & & $Z_{\gamma_4} \approx - 4.98697\I$\\
    \bottomrule
  \end{tabular}
  }
  \vspace{0.5em}
  \caption{Summary of data defining the $N=3$ examples.}
  \label{tab:N3ex}
\end{table}

\begin{table}
  \makebox[\textwidth][c]{
  \begin{tabular}{@{}lll@{}}
    \toprule
    Theory & $\langle \cdot,\cdot \rangle$ & $\Gamma'$\\
    \midrule
    $(A_2,A_1)$ & $\begin{pmatrix}0&1\\-1&0\end{pmatrix}$ & $(1,0)$, $(0,1)$, $(1,1)$\\
    \midrule
    $(A_2,A_2)$ & $\begin{pmatrix}0&1&0&0\\-1&0&0&0\\0&0&0&0\\0&0&0&0\end{pmatrix}$ & \makecell[l]{$(1,0,0,0)$, $(0,1,0,0)$, $(1,1,0,0)$,\\$(0,1,1,0)$, $(0,1,1,1)$, $(1,0,-1,0)$,\\$(1,0,-1,-1)$, $(1,-1,-1,0)$, $(1,-1,-1,-1)$,\\$(1,-1,-2,-1)$, $(1,-2,-2,-1)$, $(2,-1,-2,-1)$}\\
    \bottomrule
  \end{tabular}
  }
  \vspace{1em}
  \caption{Intersection form and classes with nonzero BPS counts in the $N=3$ examples, with respect to the bases in Table \ref{tab:N3ex}.  As before we list only one from each pair $\pm \gamma$ in $\Gamma'$.}
  \label{tab:N3coef}
\end{table}

\paragraph{Example $(A_2,A_1)$.}

We take as base point the polynomial $P_3(z) = \frac{1}{2}(1-z^2)$ and the BPS-free angle $\vartheta=0$.
Here $Y$ is a $3$-fold cover of $\C$ branched over $\pm 1$, which as a Riemann surface is the hexagonal punctured torus, i.e.~the Riemann surface obtained by gluing opposite pairs of sides of a regular hexagon, and then removing a point that corresponds to three of the original vertices.
Its homology has rank two, i.e.~$\Gamma \simeq \Z^2$, 
\begin{figure}
  \begin{center}
    \includegraphics[width=0.6\textwidth]{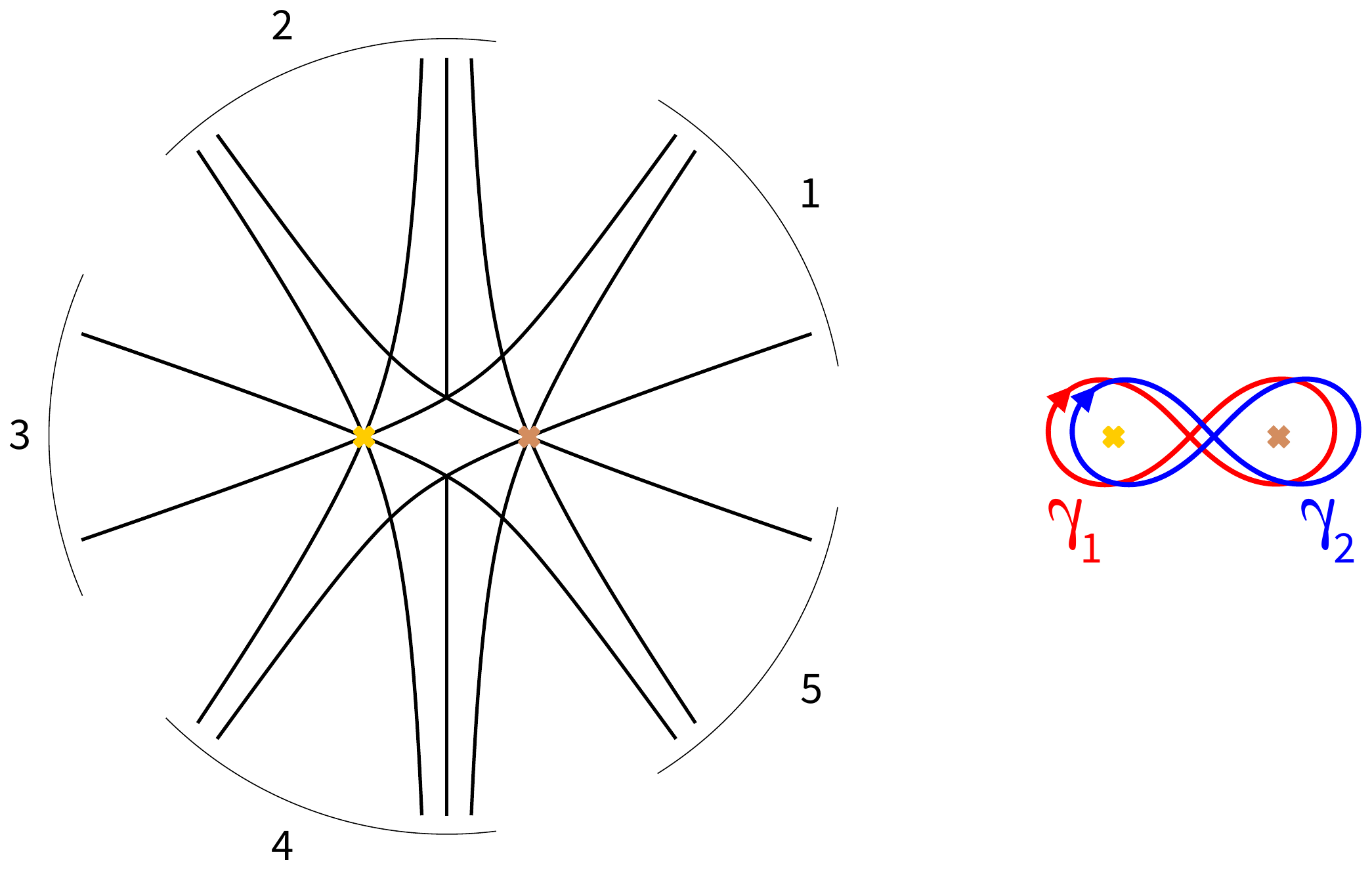}
  \end{center}
  \caption{Spectral network (left) and projected homology basis (right) for $(A_2,A_1)$, with $P_3(z) = \frac{1}{2}(1-z^2)$).}
    \label{fig:A2A1basis}
\end{figure}

To construct a homology basis for the spectral curve, we consider an oriented figure-eight curve around $\pm 1$ as shown in \autoref{fig:A2A1basis}.
This curve has three lifts to simple loops on the spectral curve, distinguished by their periods which have arguments $\frac{\pi}{6}, \frac{5 \pi}{6}, \frac{3\pi}{2}$; these correspond to three segments joining opposite pairs of sides in the hexagon model of the spectral curve.
The set $\Gamma'$ consists of these three cycles and their inverses.
We fix the basis $\gamma_1,\gamma_2$ corresponding to the lifts with period arguments $\frac{5 \pi}{6}$ and $\frac{3\pi}{2}$, respectively.

Using a correspondence between homology classes and ``abelianization trees'' described in \cite{Neitzke17}, this basis gives rise to a pair of spectral coordinates, which for $\vartheta=0$ (or more generally any $\abs{\vartheta} < \frac{\pi}{6}$) are:
\begin{equation}
  X_1 = \frac{p(2,3,4) p(1,4,5)}{p(1,2,4) p(3,4,5)}, \;\;\; X_2 = \frac{p(2,4,5) p(1,2,3) p(1,4,5)}{p(1,2,5) p(1,2,4) p(3,4,5)}.
\end{equation}

\begin{figure}
  \begin{center}
    \includegraphics[width=0.7\textwidth]{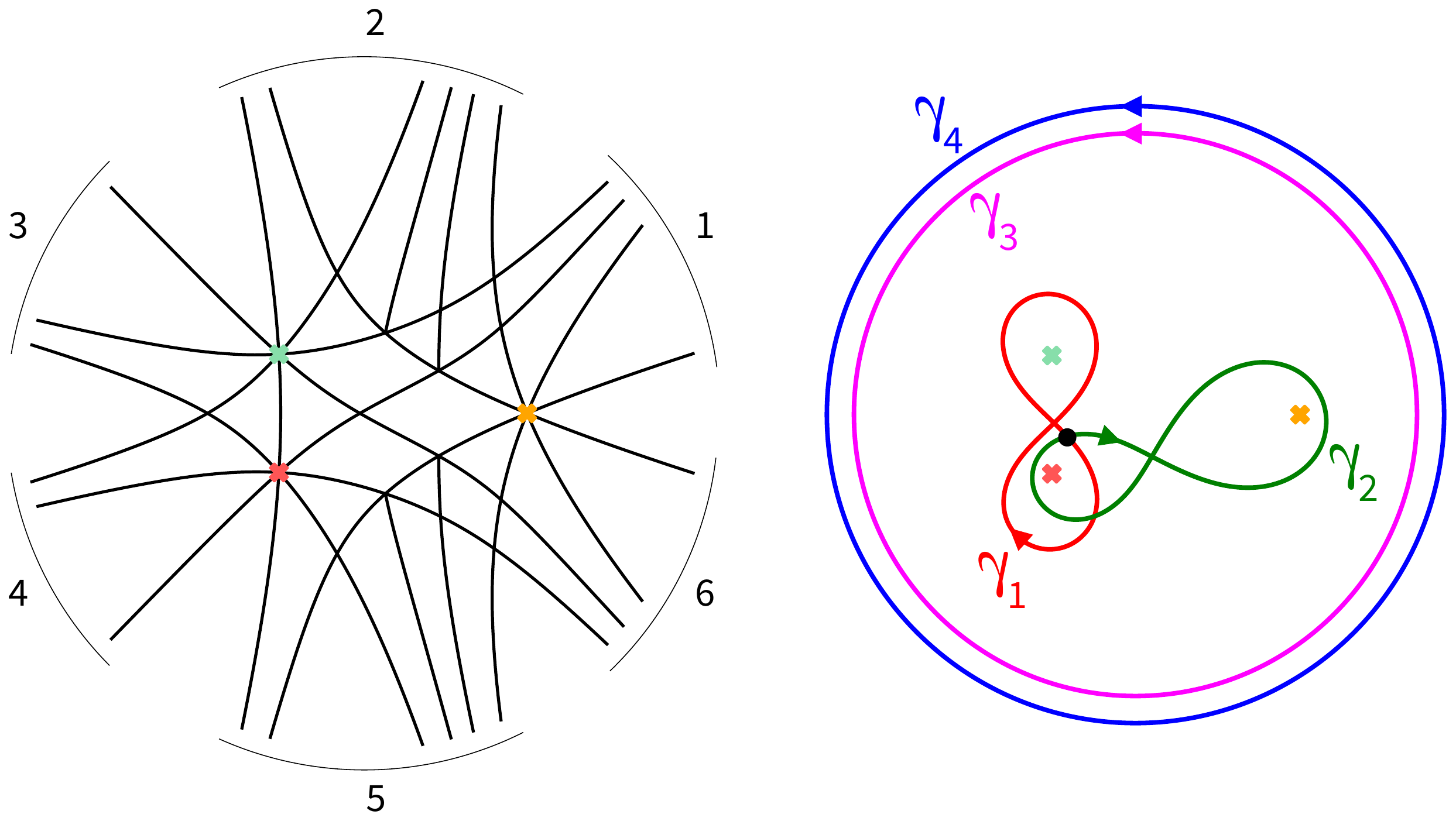}
  \end{center}
  \caption{Spectral network (left) and projected homology basis (right) for $(A_2,A_2)$, with $P_3(z) = \frac{1}{2}(z^3 - 3 z^2 - 2)$.  The cycles $\gamma_1$ and $\gamma_2$ intersect over the indicated point.}
    \label{fig:A2A2basis}
\end{figure}

\paragraph{Example $(A_2,A_2)$.}

We take as base point $P_3(z) = \frac{1}{2}(z^3 - 3 z^2 - 2)$ and BPS-free angle $\vartheta=0$.
We denote the zeros by $z_0, z_1, z_2$ with $z_1,z_2$ complex conjugates and $\im(z_1)>0$.
Here the spectral curve is a three-punctured torus, so $\Gamma \simeq \Z^4$.
We use a basis $\gamma_1, \ldots, \gamma_4$ that, as in the previous case, can be described in terms of lifts of loops around the zeros of $P_3$.
Specifically, for $\gamma_1$ and $\gamma_2$ we choose lifts of figure eight loops around $z_0,z_2$ and $z_1,z_2$, respectively, while $\gamma_3$ and $\gamma_4$ are each lifts of a large counter-clockwise circle enclosing all of the roots.
Thus, on the spectral curve, $\{ \gamma_1, \gamma_2 \}$ gives a basis of the homology of the torus obtained by filling in the punctures, while $\gamma_3$ and $\gamma_4$ are represented by small loops around punctures
(thus $\gamma_3$ and $\gamma_4$ are examples of pure flavor charges.)
These are shown in \autoref{fig:A2A2basis}.
As before the ambiguity in choice of a lift of each curve is resolved by specifying the periods, and numerical approximations of these appear in \autoref{tab:N3ex}.
In this case the set $\Gamma'$ has $24$ elements, which are indicated in \autoref{tab:N3coef}.

As explained in \cite{Neitzke17}, the associated spectral coordinates are
\begin{equation*}
  \begin{split}
  X_1 &=
  \frac{q(1,2,3,4,5,6)}{p(1,2,3)p(4,5,6)}, \;\;
  X_2 = \frac{p(1,2,5)p(3,5,6)p(4,5,6)}{p(1,5,6)p(2,5,6)p(3,4,5)}, \\  
  X_3 &= \frac{p(1,2,6)p(3,4,5)}{p(1,2,3)p(4,5,6)}, \;\;
  X_4= \frac{p(1,5,6)p(2,3,4)}{p(1,2,6)p(3,4,5)},
\end{split}
\end{equation*}
where the function $q$ in $X_1$ is the hexapod invariant discussed in \autoref{sec:stokes}.

\subsection{Results of numerical studies of opers}
\label{subsec:results-opers}

We now turn to reporting results of calculating spectral coordinates for the opers, comparing the differential equation (DE) and conjectural integral equation (IEQ) methods.
In tabulating and plotting the results for a given $w$, we fix the argument $\vartheta = \arg \hbar$ and then take $\abs{\hbar}^{-1}$ as the independent variable rather than $\abs{\hbar}$ itself.
This is convenient since $\abs{\hbar}^{-1} \to \infty$ corresponds to divergence in the moduli space, and is analogous to $R \to \infty$ in the Hitchin section results presented in the next section, thus giving the same expected behavior in the plots and tables of these two sections.

We begin with explicit numerical results in one example.
Consider the $(A_1,A_2)$ theory with $\Lambda=0.8\I$, $c=1$, 
which corresponds to taking 
\begin{equation}
w = \phi_2 = (z^3 - (0.8\I)z - 1) \,\de z^2.
\end{equation}
There are two spectral coordinates ($X_1$,$X_2$), and we denote the results of calculating the spectral coordinates by the two methods by $X_i^{\DE}(\hbar)$ and $X_i^{\IEQ}(\hbar)$.
\autoref{table:a1a2deformed-x1} and \autoref{table:a1a2deformed-x2} shows numerical results for this example for several values of $\hbar$ with $\vartheta = 0$, 
as well as the relative difference between the DE and IEQ results,\footnote{Here we define the relative difference between real or complex quantities $a$ and $b$ to be
$\mathrm{reldiff}(a,b) = \frac{2 \abs{a-b}}{\abs{a}+\abs{b}}$,
that is, $\mathrm{reldiff}(a,b)$ describes the difference as a fraction of the average of $\abs{a}$, $\abs{b}$.
}
and an estimate of the relative error in the DE results arising from numerical solution of the parallel transport ODE.
Each calculation method involves a number of internal parameters, and the calculation details and parameters used here are given in \autoref{sec:implementation}.

\begin{table}
\makebox[\textwidth][c]{
\begin{tabular}{@{}lllll@{}}
\toprule
& & & & Rel.~ODE\\
$|\hbar|^{-1}$ & $X_1^{\DE}$ & $X_1^{\IEQ}$ & $\mathrm{reldiff}(X_1^{\DE},X_1^{\IEQ})$ & err.~est. \\
\midrule
$\exp(-6)$ & $\num{6.197565467441e-01}$ & $\num{6.197565467441e-01}$ & $\num{2.5e-15}$ & $\num{1.2e-13}$\\
 & $\num{-3.078848587097e-03}$i & $\num{-3.078848587098e-03}$i &  & \\
$\exp(-3)$ & $\num{6.098896929797e-01}$ & $\num{6.098896929797e-01}$ & $\num{6.2e-14}$ & $\num{1.2e-13}$\\
 & $\num{-1.502048140005e-02}$i & $\num{-1.502048140001e-02}$i &  & \\
$\exp(0)$ & $\num{1.241494034799e-01}$ & $\num{1.241494034799e-01}$ & $\num{6.0e-13}$ & $\num{1.8e-13}$\\
 & $\num{4.675580545520e-02}$i & $\num{4.675580545526e-02}$i &  & \\
$\exp(1.5)$ & $\num{-4.541501969871e-05}$ & $\num{-4.541501970009e-05}$ & $\num{8.9e-12}$ & $\num{2.4e-11}$\\
 & $\num{1.652861158678e-04}$i & $\num{1.652861158672e-04}$i &  & \\
$\exp(2.25)$ & $\num{-7.867354971295e-09}$ & $\num{-7.867355151333e-09}$ & $\num{2.0e-08}$ & $\num{6.4e-08}$\\
 & $\num{-7.713362317977e-09}$i & $\num{-7.713362195515e-09}$i &  & \\
$\exp(3)$ & $\num{-6.050397949632e-18}$ & $\num{-6.288168615191e-18}$ & $\num{3.2e-02}$ & $\num{2.0e+00}$\\
 & $\num{1.341932452591e-17}$i & $\num{1.301568577212e-17}$i &  & \\
\bottomrule
\end{tabular}

}
\vspace{1em}
\caption{Calculated spectral coordinate $X_1$ for $(A_1,A_2)$ at $\Lambda=0.8\I$, $c=1$, $\vartheta=0$.}
\label{table:a1a2deformed-x1}
\end{table}

\begin{table}
\makebox[\textwidth][c]{
\begin{tabular}{@{}lllll@{}}
\toprule
& & & & Rel.~ODE\\
$|\hbar|^{-1}$ & $X_2^{\DE}$ & $X_2^{\IEQ}$ & $\mathrm{reldiff}(X_2^{\DE},X_2^{\IEQ})$ & err.~est. \\
\midrule
$\exp(-6)$ & $\num{9.910196824725e-01}$ & $\num{9.910196824725e-01}$ & $\num{1.7e-14}$ & $\num{1.8e-13}$\\
 & $\num{-2.459625969598e-03}$i & $\num{-2.459625969583e-03}$i &  & \\
$\exp(-3)$ & $\num{9.017178663077e-01}$ & $\num{9.017178663078e-01}$ & $\num{7.7e-14}$ & $\num{1.8e-13}$\\
 & $\num{-8.117380424180e-02}$i & $\num{-8.117380424181e-02}$i &  & \\
$\exp(0)$ & $\num{-2.797839029129e-01}$ & $\num{-2.797839029132e-01}$ & $\num{1.1e-12}$ & $\num{2.1e-13}$\\
 & $\num{-1.019844898950e-01}$i & $\num{-1.019844898951e-01}$i &  & \\
$\exp(1.5)$ & $\num{4.437252973127e-03}$ & $\num{4.437252973159e-03}$ & $\num{6.6e-12}$ & $\num{2.1e-11}$\\
 & $\num{-2.267783957724e-03}$i & $\num{-2.267783957721e-03}$i &  & \\
$\exp(2.25)$ & $\num{-1.086698308145e-05}$ & $\num{-1.086698331057e-05}$ & $\num{2.0e-08}$ & $\num{4.7e-08}$\\
 & $\num{-7.680370178260e-06}$i & $\num{-7.680370312705e-06}$i &  & \\
$\exp(3)$ & $\num{-2.017168705958e-11}$ & $\num{-2.112193509873e-11}$ & $\num{4.1e-02}$ & $\num{8.1e-01}$\\
 & $\num{-4.428933036588e-11}$i & $\num{-4.255436289723e-11}$i &  & \\
\bottomrule
\end{tabular}

}
\vspace{1em}
\caption{Calculated spectral coordinates $X_2$ for $(A_1,A_2)$ at $\Lambda=0.8\I$, $c=1$, $\vartheta=0$.}
\label{table:a1a2deformed-x2}
\end{table}

A pattern seen in these tables is present in all of the computations we report:
For sufficiently small $\abs{\hbar}^{-1}$ the two methods are in close agreement, but for larger $\abs{\hbar}^{-1}$ the relative difference grows quickly.
This is to be expected, 
since the relative numerical error in the results of the $\DE$ computation is
expected to grow with $\abs{\hbar}^{-1}$.

Plots of the spectral coordinates and relative errors for all of the examples discussed in \autoref{subsubsec:N2ex} and \autoref{subsubsec:N3ex} are shown on the next several pages (\operfigs).
Each of these ``four-pane'' plots has the following structure:
The top row of plots shows results for $\abs{\hbar}^{-1} < 0.1$ (``small parameter'') and the bottom show results for $\abs{\hbar}^{-1} \geq 0.1$ (``large parameter'').
In each case, the left plot shows the spectral coordinates themselves (as computed by both methods), and the right shows the relative difference between the two methods, as well as the relative difference in the $\DE$ result corresponding to an estimate of the error in that calculation.
The error model used in this estimate is described in \autoref{subsec:error-opers}.
The upper limit of $\abs{\hbar}^{-1}$ in each set of experiments is chosen as a point where the error estimate for the direct method becomes comparable to the spectral coordinate itself; 
beyond that point, the $\DE$ results are dominated by numerical error and comparison with $\IEQ$ is meaningless.
The relative differences are always shown on a logarithmic $y$-axis scale, and all relative error plots use the same $y$-axis limits ($5 \times 10^{-14}$ to $5$).
The scales for the other axes are adapted to the different regions:
For small $\abs{\hbar}^{-1}$, the $\abs{\hbar}^{-1}$ axis uses a logarithmic scale, as is suited to the exponential spacing of the sample points.
For large $\abs{\hbar}^{-1}$, the $\abs{\hbar}^{-1}$ axis uses a linear scale and the $|X_i|$ axis is logarithmic; this has the effect of making the leading-order WKB asymptotic \eqref{eq:wkb-asymp} a linear function, which is shown as a dashed line.
For small $\abs{\hbar}^{-1}$, the WKB asymptotic is not expected to be accurate and is not shown, except for the pure flavor coordinates $X_3$ and $X_4$ of the $(A_2,A_2)$ example where it is expected to give an exact formula.

\newpage

\vspace*{\floatpagetop}

\newpage
\vspace*{\floatpagetop}

\newpage
\vspace*{\floatpagetop}

\newpage
\vspace*{\floatpagetop}

\newpage
\vspace*{\floatpagetop}

\newpage
\vspace*{\floatpagetop}

\newpage
\vspace*{\floatpagetop}
\newpage

\subsection{Results of numerical studies for the Hitchin section}
\label{subsec:results-hitchin}

Now we turn to the case of the flat connections $\nabla^\higgs_w(R,\zeta)$ 
discussed in \autoref{subsec:hitchin}, and reporting the results of computing the associated spectral coordinates by the DE and IEQ methods.
As with our numerical results for opers, we begin by presenting a table of computed values in the theory $(A_1,A_2)$.
In this case we tabulate results for the basepoint $P_2(z) = z^3-1$ (i.e.~parameters $\Lambda=0$, $c=1$) and $\zeta = \exp(\I \vartheta)$ where $\vartheta = 0.1$.
Spectral coordinates calculated for several values of $R$ are shown in \autoref{table:a1a2-hitchin}.
Recall that the parameter $R$ in this case is analogous to $\abs{\hbar}^{-1}$ for opers.
Analogously to the results in the oper case, we see that the relative difference of the $\DE$ and $\IEQ$ methods grows as $R$ increases.

The rightmost column of \autoref{table:a1a2-hitchin} shows an error estimate for some of the $\DE$ calculations, which is also included in all of the plots to follow.
This estimate is not based on a theoretical error analysis, but rather on testing the dependence of the $\DE$ results on the grid spacing $\Delta x$ in the discretization of the PDE and applying Richardson extrapolation to predict a limit value as $\Delta x \to 0$.
When the dependence on the spacing is approximately quadratic in $\Delta x$ (the expected form), the difference between the extrapolated value and the one calculated with the finest grid is taken as an estimate of PDE discretization error.
In other cases the dependence on $\Delta x$ does not exhibit the expected form, and no error estimate is obtained;
this would be expected to happen when, for example, discretization error is not the dominant source of error in the $\DE$ calculation.
This error estimation technique is described in more detail in \autoref{subsec:error-hitchin}.

\begin{table}
\makebox[\textwidth][c]{
\begin{tabular}{@{}lllll@{}}
\toprule
& & & & Rel.~PDE\\
$R$ & $X_1^{\DE}$ & $X_1^{\IEQ}$ & $\mathrm{reldiff}(X_1^{\DE},X_1^{\IEQ})$ & error est. \\
\midrule
$\exp(-6)$ & $\num{6.161191878e-01}$ & $\num{6.161191872e-01}$ & $\num{1.0e-09}$ & $\num{2.3e-10}$\\
$\exp(-3)$ & $\num{5.511806552e-01}$ & $\num{5.511806344e-01}$ & $\num{3.8e-08}$ & $\num{8.5e-09}$\\
$\exp(0)$ & $\num{8.800765041e-03}$ & $\num{8.800772243e-03}$ & $\num{8.2e-07}$ & --\\
$\exp(1.5)$ & $\num{6.111579918e-10}$ & $\num{6.171310177e-10}$ & $\num{9.7e-03}$ & $\num{9.8e-03}$\\
$\exp(1.875)$ & $\num{1.418541282e-13}$ & $\num{3.981788142e-14}$ & $\num{1.1e+00}$ & $\num{2.0e+00}$\\
\midrule
& & & & Rel.~PDE \\
$R$ & $X_2^{\DE}$ & $X_2^{\IEQ}$ & $\mathrm{reldiff}(X_2^{\DE},X_2^{\IEQ})$ & error est.\\
\midrule
$\exp(-6)$ & $\num{9.994046408e-01}$ & $\num{9.994046406e-01}$ & $\num{2.0e-10}$ & $\num{4.8e-11}$\\
$\exp(-3)$ & $\num{9.787121061e-01}$ & $\num{9.787120991e-01}$ & $\num{7.1e-09}$ & $\num{1.7e-09}$\\
$\exp(0)$ & $\num{5.589112961e-01}$ & $\num{5.589110919e-01}$ & $\num{3.7e-07}$ & $\num{5.7e-07}$\\
$\exp(1.5)$ & $\num{7.369086545e-02}$ & $\num{7.368779561e-02}$ & $\num{4.2e-05}$ & $\num{4.3e-05}$\\
$\exp(1.875)$ & $\num{2.252230986e-02}$ & $\num{2.249418541e-02}$ & $\num{1.2e-03}$ & $\num{1.2e-03}$\\
\bottomrule
\end{tabular}

}
\vspace{1em}
\caption{Calculated spectral coordinates for $(A_1,A_2)$ Hitchin section, $\Lambda=0$, $c=1$, $\vartheta=0.1$.}
\label{table:a1a2-hitchin}
\end{table}

Plots of the spectral coordinates and relative errors for this example and the others introduced above are shown on the next several figures (\hitchinfigs).
Each of these ``four-pane'' plots has the same structure described in \autoref{subsec:results-opers}, with the additional complication that error estimates are only shown for values of $R$ where the Richardson extrapolation is successful.
Generally, the extrapolation succeeds for most large $R$ and yields an error estimate that increases with $R$.
The upper limit of $R$ in each experiment is chosen to be a point where the resulting error estimate first becomes comparable in size to the spectral coordinates themselves, i.e.~the largest $R$ for which this error estimate suggests the $\DE$ results are meaningful.
Analogously to the WKB asymptotic in the opers results, the semiflat approximation to $X_i$ is shown as a dashed line in the large-$R$ plots (where it is expected to be a good approximation) and in all $R$ for the pure flavor coordinates $X_3$ and $X_4$ of the $(A_2,A_2)$ example (where it is expected to be exact).

\clearpage 

\vspace*{\floatpagetop}

\newpage
\vspace*{\floatpagetop}

\newpage
\vspace*{\floatpagetop}

\newpage
\vspace*{\floatpagetop}

\newpage
\vspace*{\floatpagetop}

\newpage
\vspace*{\floatpagetop}

\newpage

\section{Experimental studies of the \hk metric} \label{sec:hk-experiments}

We consider the metric discussed in \autoref{sec:hk-metric} 
and \autoref{sec:spectral-hitchin}, in 
the $(A_1,A_2)$ example. Take
\begin{equation} \label{eq:def-hl}
	h(z) = z^3 - \Lambda z \ (\Lambda \in \C), \quad l(z) = -c \ (c \in \C).
\end{equation}
This gives a 1-parameter family of \kahler manifolds $\cB_\Lambda := \cB_{d,h}$,
indexed by $\Lambda \in \C$; each $\cB_\Lambda$ is 1-dimensional 
(coordinatized by $c \in \C$) and carries a \kahler metric
\begin{equation}
	g = g(c) \abs{\de c}^2.
\end{equation}

\subsection{Direct PDE computation} \label{sec:hk-direct}

Our first approach to computing the metric $g$ is to use 
the definition as an $L^2$ norm directly.
Thus, given the polynomial 
\begin{equation} \label{eq:hkmetric-polys}
P_2(z) = z^3 - \Lambda z - c,	
\end{equation}
and the tangent vector corresponding to $\partial / \partial c$,
\begin{equation}
	\dot{P}_2(z) = -1,
\end{equation}
we first solve the nonlinear PDE \eqref{eq:self-duality} for $u$,
then solve the linear PDE \eqref{eq:diff-F} for $F$,
then compute the integral \eqref{eq:l2norm-combined} to
get the desired metric coefficient $g(c)$.

\subsection{Integral equation computation} \label{sec:hk-ieq}

Our second method of computing the metric $g$ is the integral
equation method discussed in \autoref{sec:spectral-hitchin}.
In the $(A_1,A_2)$ example there are just two independent spectral 
coordinates $X_1$, $X_2$, and the formula
\eqref{eq:metric-spectral} specializes to
\begin{equation}
  g(c) = \de \log y_1 \otimes \de^c \log y_2 - \de \log y_2 \otimes \de^c \log y_1.
\end{equation}
Even more concretely, if we write $c = a + \I b$, then
\begin{align}
  g(c) &= g(\partial_a, \partial_a) = (\partial_a \log y_1) (\partial_b \log y_2) - (\partial_a \log y_2) (\partial_b \log y_1). \label{eq:gc-concrete}
\end{align}
We use the integral equation method from \autoref{subsubsec:inteq-hitchin} to compute $y_1$ and $y_2$ at various values of $c$, thus compute the derivatives appearing \eqref{eq:gc-concrete} by finite differences, and finally use \eqref{eq:gc-concrete} to compute $g(c)$.

\subsection{Experimental comparison}

We have described two methods of computing the \kahler metric
on $\cB_{\Lambda}$. We applied both of these methods to
compute $g(c)$ in the case $\Lambda = 0$ and $c \in \R_+$.
When $\Lambda = 0$ there is a rotational symmetry, so that 
$g(c)$ depends only on $\abs{c}$; thus the values
of $g(c)$ for $c \in \R_+$ determine the full $g$.
The result is shown in \autoref{fig:paperfighkmetric};
the observed difference $\abs{g^{\DE} - g^{\IEQ}} < 1.2 \times 10^{-4}$
over the range of $c$ we studied.

\autoref{fig:paperfighkmetric} also shows the semiflat approximation
$g^\rmsf(c)$ discussed in \autoref{sec:semiflat} and \autoref{sec:leading}. In this 
example we can compute $g^\rmsf$ 
in closed form, using \eqref{eq:semiflat-spectral} and the fact
that $Z_\gamma(c) = Z_\gamma(c=1) c^{\frac56}$, with $Z_\gamma(c=1)$ tabulated
in \autoref{tab:N2ex}; the result is
\begin{equation}
g^\rmsf(c) = \frac{25 M^2}{6 \sqrt{3}} \abs{c}^{-\frac13} \approx 20.4325 \abs{c}^{-\frac13}, \qquad M = \sqrt{3 \pi} \frac{\Gamma\left(\frac43\right)}{\Gamma\left(\frac{11}{6}\right)}.
\end{equation}
The figure shows that the semiflat approximation is increasingly
accurate for large $\abs{c}$ and not at all accurate
for small $\abs{c}$, as expected: 
It could hardly be accurate near $c=0$ since $g^\rmsf(c)$ 
has a singularity at that point while $g(c)$ is smooth.

\begin{figure}
  \inswide{{paperfighkmetric}.pdf}
  \caption{Left: The metric coefficient $g(c)$ for $\Lambda = 0$
  and $c \in \R_+$.
  The blue marks show values of $g(c)$ 
  computed using two methods: the direct PDE approach and the integral equations.
  The dashed line shows the semiflat approximation. Right: The absolute difference
  $g^{\mathrm{DE}} - g^{\mathrm{IEQ}}$.}  
  \label{fig:paperfighkmetric}
\end{figure}

\section{Gallery}
\label{sec:gallery}

\subsection{The \hk metric integrand}

\begin{figure}
  \begin{center}
    \includegraphics[width=\textwidth]{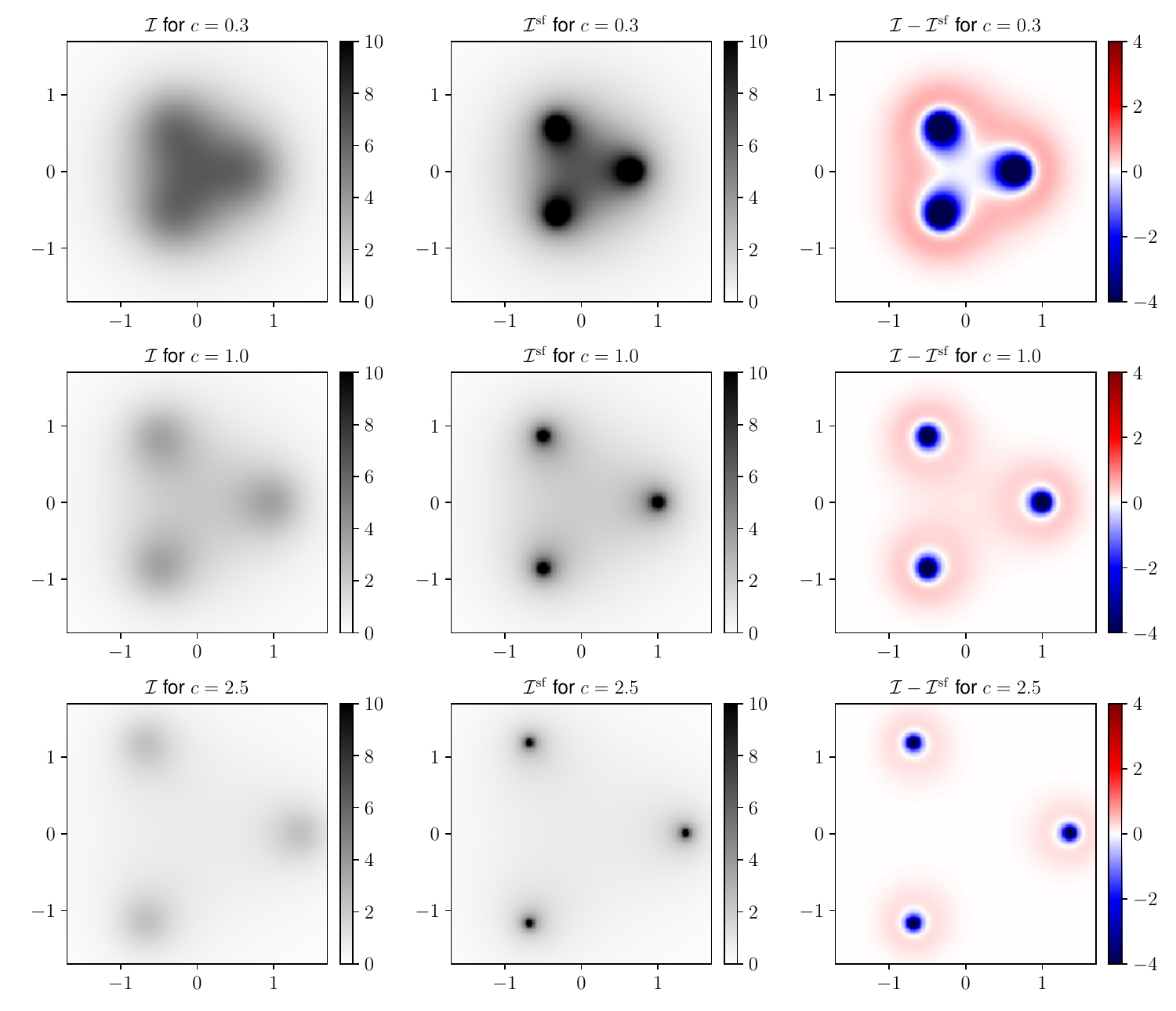}
  \end{center}
  \caption{The integrand in the computation of the \hk metric
  as an $L^2$ norm, when
  $P_2(z) = z^3 - c$, $\dot{P}_2(z) = -1$.
  Left: the integrand $\cI$ of the $L^2$ metric \eqref{eq:l2norm-combined}, 
  plotted in the $z$-plane.
  Middle: the integrand $\cI^{\rmsf}$ of the semiflat approximation to the $L^2$ metric.
  Right: the difference $\cI - \cI^{\rmsf}$.}
  \label{fig:paperfigheatmaps}
\end{figure}

In \autoref{fig:paperfigheatmaps} we illustrate some features of the
numerical computation of the metric \eqref{eq:l2norm-combined}
in the simple case $N=2$, $d=3$, $P_2(z) = z^3 - c$, $\dot{P}_2(z) = -1$.
The theoretical expectation based on \cite{Gaiotto:2009hg,Mazzeo2014,Dumas2018} is that
\begin{itemize}
\item the pointwise difference $\cI - \cI^{\rmsf}$ decays exponentially as 
a function of the distance 
from the zeros of $P_2(z)$ (measured in the metric $\abs{P_2 \, \de z^2}$),
\item the integral of $\cI - \cI^{\rmsf}$ over a disc in the metric $\abs{P_2 \, \de z^2}$,
centered on a zero of $P_2$ and not containing any other zero, 
decays exponentially as a function of the radius
of the disc.
\end{itemize}
In other words, $\cI-\cI^{\rmsf}$ can be large near the zeros, but there is 
a local cancellation around each zero which makes its \ti{integral} nevertheless
small; see \cite{Dumas2018} for the precise statement.
We see this phenomenon in the figure: 
near each zero we have
$\cI - \cI^{\rmsf}$ large and negative, but there is a halo a bit further out, where
$\cI - \cI^{\rmsf}$ is positive. 

We also observe that as $\abs{c}$ increases the error
$\cI - \cI^{\rmsf}$ becomes more concentrated around
the zeros, as expected since the distances in the metric $\abs{P_2 \de z^2}$
grow as $\abs{c}$ increases.
Moreover, in the limit of large $\abs{c}$ the individual zeros
effectively decouple from one another; indeed the solution in a neighborhood of
each zero approaches a standard ``fiducial'' solution \cite{Cecotti:1991me,Gaiotto:2009hg} 
when written in the coordinate
$w = \int \sqrt{P_2} \de z$.

\subsection{The functions \texorpdfstring{$\cX_\gamma$}{Xgamma}}

We consider the $(A_1,A_2)$ example.
Let $x_\gamma^\inst$ denote the integral term 
in \eqref{eq:integral-oper} (for opers) 
or \eqref{eq:integral-hitchin-section} (for the Hitchin section).

In this case the integral equations \eqref{eq:integral-oper}, \eqref{eq:integral-hitchin-section}
coincide with ones which have been studied extensively in the
literature on the thermodynamic Bethe ansatz, beginning
with \cite{Zamolodchikov:1990cf}, and also
in the context of the ODE/IM correspondence beginning
with \cite{Dorey:1999uk}.
All of the main features of the $x_\gamma^\inst$ which we discuss in this section are also
noted in \cite{Zamolodchikov:1990cf}; we present them here for completeness
and for readers not familiar with that reference.

\begin{figure}
  \begin{center}
    \includegraphics[width=0.8\textwidth]{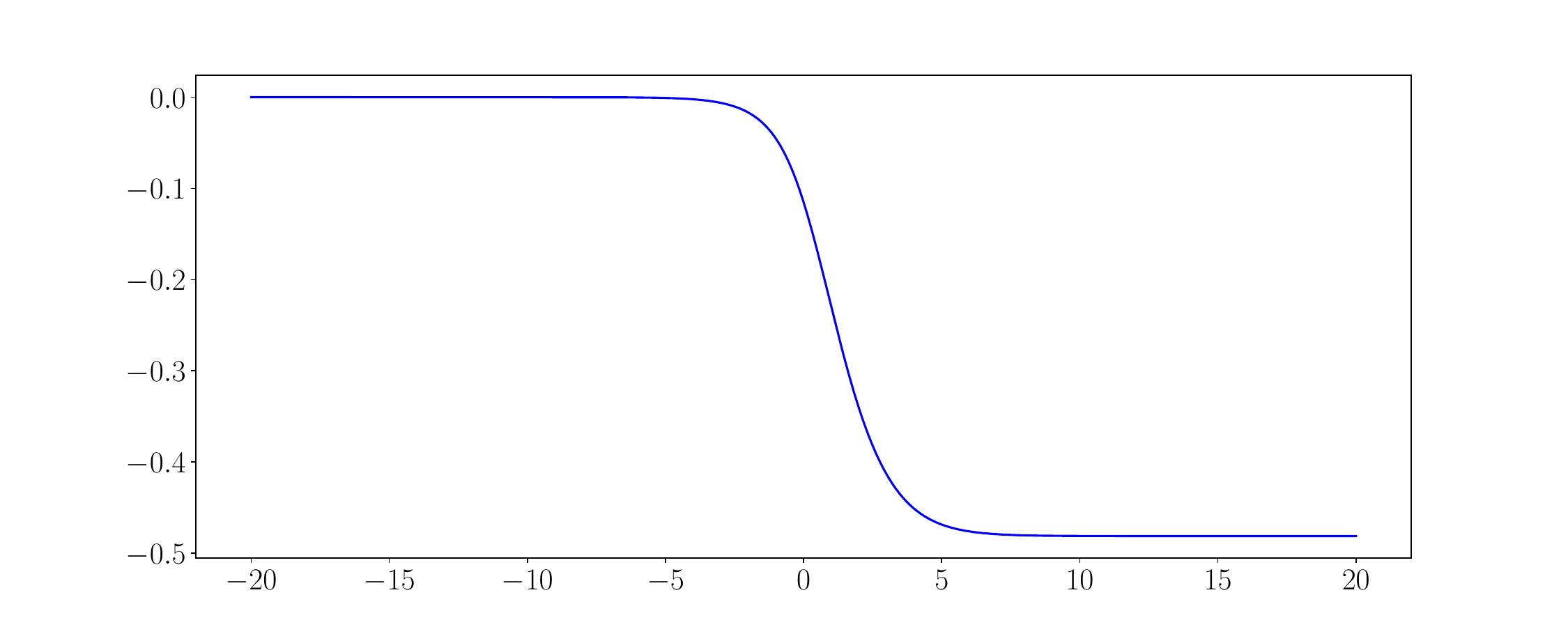}
  \end{center}
  \caption{$x^{\mathrm{inst}}_{\gamma_1}(\hbar)$ in the $(A_1,A_2)$
  example, evaluated at the values $\hbar = - \exp(t + \I \arg Z_\gamma)$, for $t \in [-20,20]$ 
  (horizontal axis).}
  \label{fig:paperfigxaroper}
\end{figure}

In \autoref{fig:paperfigxaroper} we show the function $x^\inst_{\gamma_1}(\hbar)$, evaluated along the ray $\hbar \in \R_- Z_{\gamma_1}$.
A few features are worthy of comment:
\begin{itemize}
\item As $\hbar \to 0$, $x^\inst_{\gamma_1}(\hbar)$ approaches $0$.
This confirms our expectation from \autoref{sec:leading},
and the consequence that in this limit the full
$x_{\gamma_1}(\hbar) = \hbar^{-1} Z_{\gamma_1} + x^\inst_{\gamma_1}$
is asymptotic to $\hbar^{-1} Z_{\gamma_1}$.

\item As $\hbar \to \infty$, $x^\inst_{\gamma_1}(\hbar)$ approaches
a nonzero finite limit, and hence so does 
the full $x_{\gamma_1}(\hbar) = \hbar^{-1} Z_{\gamma_1} + x^\inst_{\gamma_1}$.
This limit corresponds to the polynomial $P_2(z) = z^3$, for which the
oper has a $\Z / 5\Z$ symmetry 
which determines its Stokes data as
\begin{equation}
	\lim_{\hbar \to \infty} x_{\gamma_1}(\hbar) = x_* := \log \left( \frac{\sqrt{5} - 1}{2} \right) \approx -0.4812,
\end{equation}
matching the asymptotic value in the figure.
\end{itemize}

\begin{figure}
  \begin{center}
    \includegraphics[width=0.8\textwidth]{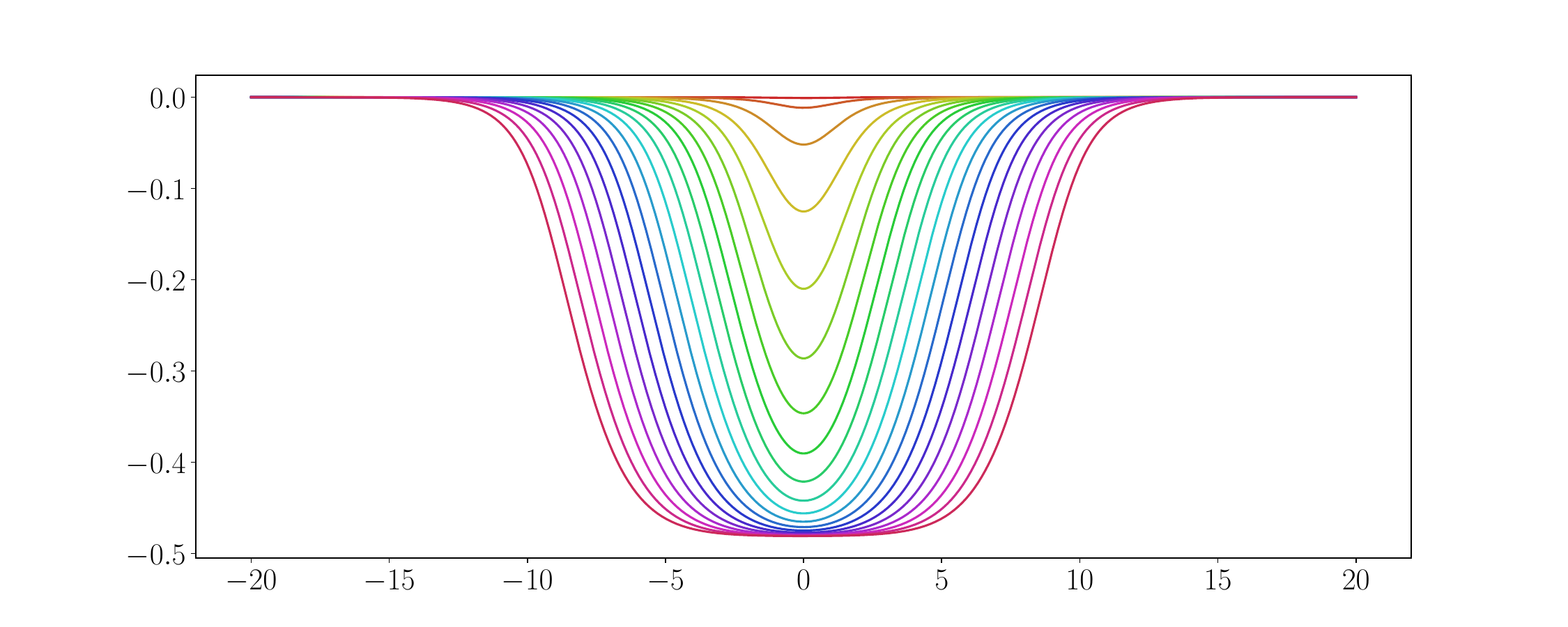}
  \end{center}
  \caption{$x^{\mathrm{inst}}_{\gamma_1}(R, \zeta)$ in the $(A_1,A_2)$
  example, evaluated at the values $\zeta = - \exp(t + \I \arg Z_\gamma)$, for $t \in [-20,20]$ 
  (horizontal axis). The $20$ curves correspond to the values 
  $R = \exp(-k/2)$ for $k = 0, 1, \dots, 19$.}
  \label{fig:paperfigxarshitchin}
\end{figure}
In \autoref{fig:paperfigxarshitchin} we show the function $x^\inst_{\gamma_1}(R, \zeta)$, evaluated along the ray $\zeta \in \R_- Z_{\gamma_1}$,
for various values of $R$.
Some features apparent from \autoref{fig:paperfigxarshitchin} are:
\begin{itemize}
\item For all $R$, $x^\inst_{\gamma_1}(R,\zeta) \to 0$ as $\zeta \to 0$
or $\zeta \to \infty$.
\item For small $R$, 
$x^\inst_{\gamma_1}(R, \zeta)$ has an approximate plateau at the value $x_*$, 
in a neighborhood of $\abs{\zeta} = 1$; the length of this plateau grows
as $R \to 0$, so that for any fixed $\zeta$, $\lim_{R \to 0} x^\inst_{\gamma_1}(R, \zeta) = x_*$
(but not uniformly in $\zeta$).
\item For small $R$, the crossover region between $x \approx 0$
and $x \approx x_*$ has a universal shape, which moreover looks
like the graph of $x^\inst_{\gamma_1}(\hbar)$ in \autoref{fig:paperfigxaroper},
where we make the substitution $\zeta = \hbar R$.
\end{itemize}

This last feature is a manifestation of the ``conformal limit''
which we discussed in \autoref{sec:conformal-limit};
indeed the Stokes data of $\nabla^\higgs_w(R,\zeta)$
should converge to those of $\nabla^\oper_w(\hbar)$ in that limit,
which would imply
\begin{equation}
	\lim_{R \to 0} x_{\gamma_1}^\inst(R,\zeta = \hbar R) = x_{\gamma_1}^\inst(\hbar),
\end{equation}
and this is what we observe in the figures.

\section{Implementation details}
\label{sec:implementation}

In this section we discuss the implementation of the experiments presented in sections \ref{sec:experiments}--\ref{sec:gallery} in more detail.
The source code is available at \cite{code}.

\subsection{Direct method for opers}
\label{subsec:implementation-opers}

\begin{table}
\begin{tabular}{@{}lll@{}}
    \toprule
    Name & Value & Description\\
    \midrule
    \texttt{ode\_method} & \texttt{dopri5} & ODE Solver from \texttt{scipy.ode}\\
    \texttt{ode\_thresh} & $10^{-14}$ & Relative error goal for ODE solver \\
    \texttt{ode\_rstep} & $10^{-4}$ & Initial ODE step size\\
    \bottomrule
\end{tabular}

\vspace{1.5em}

\caption{Parameters used for the direct method for opers (ODE solver) in the calculations presented in \autoref{subsec:results-opers}.}
\label{tab:oper-params}
\end{table}

The parameter values used in computing Stokes data for opers by the direct method (as reported in \autoref{subsec:results-opers}) are shown in Table \ref{tab:oper-params}.
The meanings of these parameters are described below.
In our Python implementation, direct method calculations for opers are performed by the \texttt{framedata} module.

Recall that the direct method computes parallel transport matrices for $\nabla_w^{\oper}$ from a basepoint $z_0$ to $d+N$ points on a circle of radius $r$, and then uses the eigenvectors of these matrices of smallest eigenvalue to approximate the subdominant solutions.

More precisely, given the tuple of differentials $w$ (as coefficient vectors of the associated polynomials) and the parameter $\hbar$, we first compute the scaled tuple $\hbar^{-1} w = (\hbar^{-2} \phi_2)$ or $\hbar^{-1} w = (\hbar^{-2} \phi_2, \hbar^{-3} \phi_3)$ 
and then apply a holomorphic change of coordinates to make the connection behave as much like the one associated to $z^d \de z^N$ as possible.
Specifically, we find $a \in \R_+$ and $b \in \C$ so that pulling back by the coordinate change $f(z) = az+b$ has the effect of making the leading coefficient of $\hbar^{-N} P_N$ have unit modulus and so that the coefficient of $z^{d-1}$ vanishes (properties we call ``quasi-monic'' and ``centered'', respectively).
After this change, it is natural to take $z_0 = 0$ as the basepoint for parallel transport, and we use the bisectors of the Stokes sectors as the directions for the $d+N$ rays.
We select the radius $r$ by finding a disk $|z| < r_0$ containing the roots of all of the nonzero polynomials $P_k$ and then setting $r = \max(8, 8r_0)$.
In practice only $P_N$ is considered here, as the examples we consider have $P_k$ constant for $k<N$.
This choice for $r$ is based on the heuristic that deviation of the connection from its asymptotic behavior is concentrated near the zeros of $P_N$, and so we select a radius significantly larger than those of the zeros.

With the entries of the connection form $\varphi_{\hbar^{-1} w}$ given explicitly by
\eqref{eq:oper-2} or \eqref{eq:oper-3}, the computation of the parallel transport along a segment (which we parameterize by $[0,1]$) now reduces to solving an explicit ODE; for this we use the \texttt{scipy.ode} module with the \texttt{dopri5} integration method (an implementation of the Dormand-Prince method of order 4(5)).
This ODE solver is applied with a fixed relative precision goal \texttt{ode\_thresh}, an initial step size \texttt{ode\_step}, and a maximum step size $2*$\texttt{ode\_step}.
Such a solution is computed for a segment $[0,r e^{\I \eta}]$ bisecting each Stokes sector, resulting in $d+N$ \emph{frame matrices} $F_i$.
The eigenvectors of $F_i$ are then computed and the normalized eigenvector with minimum eigenvalue is selected, giving the \emph{subdominant vectors} $v_i$.
These vectors represent the values in the fiber over $0$ of horizontal sections that approximate the subdominant solutions for $\nabla_w^{\oper}$.
Finally, the spectral coordinates $X_i^{\DE}$ are computed by taking ratios of products of determinants formed from the subdominant vectors.

While this method of calculation is simple to implement, it suffers from significant loss of relative precision when the determinants involved in $X_i$ are close to zero, as these determinants are sums of floating-point numbers of approximately unit norm.
Unfortunately this is the generic case for large $\abs{\hbar}^{-1}$: The asymptotic behavior of the coordinates $X_i$ is exponential in $\abs{\hbar}^{-1}$, and the individual determinants are bounded, so the generic situation of $X_i \to 0$ or $X_i \to \infty$ requires at least one determinant to approach zero.
Thus it is expected that this method of calculation will be accurate only for sufficiently small $\abs{\hbar}^{-1}$.

\subsection{ODE error estimate}
\label{subsec:error-opers}

We now explain the error estimate that is included in \operfigs (results of calculations for opers).
Recall that this estimate concerns the effect on the spectral coordinates of the limited accuracy of the numerical solution of the parallel transport ODE.
We expect this to be the dominant source of error for large $\abs{\hbar}^{-1}$.

We first consider the calculations that apply to a single Stokes sector, which involve a frame matrix (numerical approximation of parallel transport) $F$ and its eigenvectors $v_1, \ldots, v_N$.
Let $\hat{F}$ denote the corresponding exact parallel transport matrix for the same points.
More generally in this section we use a hat decoration to indicate exact quantities, in contrast to computed approximations.
The ODE solver is given a requested relative tolerance \texttt{ode\_thresh} and absolute tolerance $0$.
Assuming that the solver produces an approximate solution satisfying this request, the result is that the error in the frame matrix $\delta F := F - \hat{F}$ satisfies
\begin{equation}
\label{eq:frame-variation}
|(\delta F)_{ij}| \leq \texttt{ode\_thresh} \cdot |F_{ij}|.
\end{equation}
To analyze the propagation of this error to subsequent calculations, we will freely use linearizations of the functions applied to the frame matrices, and will then derive upper bounds on the resulting expressions.
The results are therefore estimates for an upper bound on the error, but they do \emph{not} constitute rigorous upper bounds on the error due to the use of linearization.
For brevity we will use the term \emph{estimated error bound} to refer to such an upper bound on the linearized error, and write $\delta A \lesssim B$ to mean that $B$ is such an estimated error bound for $\delta A$.

Let $v_1, \ldots, v_N$ denote the normalized eigenvectors of $F$, ordered so that the eigenvectors $\lambda_1, \ldots, \lambda_N$ increase in magnitude.
The subsequent calculations involve only the lowest eigenvector $v_1$ (the subdominant vector), the error in which can be estimated in terms of $\delta F$ using the first-order variation formula for eigenvectors \cite[Section 2.10]{wilkinson65}:
\begin{equation}
\label{eq:vertex-variation}
\delta v := v - \hat{v} \approx \sum_{j=2}^{N} \frac{v_j^* (\delta F) v_1}{\lambda_1 - \lambda_j}
\end{equation}
Here $v_1^*, \ldots, v_N^*$ denotes the dual basis.
Using \eqref{eq:frame-variation}, each component of the vector on the right hand side of \eqref{eq:vertex-variation} has absolute value bounded by the corresponding component of
\begin{equation}
\label{eq:vertex-variation-bound}
\texttt{ode\_thresh} \cdot \sum_{j=2}^{N} \frac{|v_j^*| |F| |v_1|}{|\lambda_1 - \lambda_j|}
\end{equation}
In this expression, the absolute value of a vector or matrix represents the result of applying the absolute value componentwise, i.e.~the components of $|v|$ are the absolute values of the components of $v$, and the entries of $|F|$ are $|F_{ij}|$.
The expression above is thus used as the estimated error bound for the components of each subdominant vector.

Turning now to the calculation of the determinantal invariants and spectral coordinates, it is convenient to change notation slightly and denote by $v_1, \ldots, v_{d+N}$ the collection of subdominant vectors of the frame matrices for all of the Stokes sectors.
Having estimated the componentwise error in each of the vectors $v_i$, we now promote this to an estimated error bound for the relevant invariants $p(i,j,k)$ and $q(i,j,k,l,m,n)$.
To do this we compute the partial derivative of the invariant at $(v_1, \ldots, v_{d+N})$ with respect to each vector component and contract this with \eqref{eq:vertex-variation-bound}.
In our implementation, the partial derivatives of the determinantal invariants are numerically approximated by finite differences with a fixed step size of $10^{-12}$.

Finally, the spectral coordinates have the form $X_i = \frac{A_1 \cdots A_k}{A_{k+1} \cdots A_r}$ where each quantity $A_i$ is one of the determinental invariants discussed above.
Using logarithmic differentiation we arrive at an estimated error bound for $X_i$ in terms of those of $A_i$,
\begin{equation}
|\delta X_i| \lesssim |X_i| \sum_{j} \frac{|\delta A_j|}{|A_j|}
\end{equation}
Our final linearized estimate for $\delta X_i$ is obtained by substituting the error estimate for each invariant $A_i$ obtained above.

\subsection{Direct method for the Hitchin section}
\label{subsec:implementation-hitchin}

The parameter values used in computing Stokes data for the Hitchin section by the direct method (as reported in \autoref{subsec:results-hitchin}) are shown in Table \ref{tab:hitchin-params}, and the meanings of these parameters are described below.
In our Python implementation, direct method calculations for the Hitchin section are performed by the \texttt{framedata} module.

\begin{table}
\begin{tabular}{@{}lll@{}}
    \toprule
    Name & Value & Description\\
    \midrule
    \texttt{method} & \texttt{fourier} & PDE solver strategy (\texttt{euler} or \texttt{fourier})\\
    \texttt{pde\_nmesh} & $8191$ & PDE mesh size for presented results\\
     & $2047,4095,8191$ & PDE mesh sizes for Richardson extrapolation error estimate \\
    \texttt{pde\_thresh} & $1 \times 10^{-9}$ & Absolute error goal for PDE solver \\
    \multicolumn{3}{c}{\emph{\& ODE parameters from Table \ref{tab:oper-params}}}\\
    \bottomrule
\end{tabular}

\vspace{1.5em}

\caption{Parameters used for the direct method for Hitchin section (PDE solver) in the calculations presented in \autoref{subsec:results-hitchin}.}
\label{tab:hitchin-params}
\end{table}

Recall (from \autoref{subsec:direct-method}) that the direct method for the Hitchin section builds on the same ODE integration technique applied to opers, and hence it involves all of the same parameters and solution steps used there, as well as an important additional step:
For the Hitchin section, the connection matrix involves the density function $u$ of the harmonic metric, which is computed by numerically solving the self-duality equation (\autoref{eq:self-duality}).

To do this, we first discretize the problem by introducing a uniform rectangular grid in $\{ \re(z) \leq r, \im(z) \leq r \}$ of size \texttt{pde\_nmesh}$\times$\texttt{pde\_nmesh}.
The same radius $r$ is used for the subsequent ODE solution step, and as in the case of opers we choose $r$ to exceed the magnitude of the roots of $P_N$ by a significant margin.
In this case the precise algorithm to select the radius is slightly more complicated, incorporating a heuristic to balance two potential sources of error in the final results; the algorithm itself is documented in the source (\texttt{approx\_best\_rmax()} in \cite[\texttt{polynomial\_utils.py}]{code}).

Next, we compute an approximation to the harmonic metric $u$ as a function on the grid.
Rather than working with $u$ directly, we introduce a smooth function $u_0$ (the \emph{model}) that is computed directly from $P_N$ in closed form, and which has the same asymptotic behavior as $u$.
We then consider the difference $v = u-u_0$ and the PDE equivalent to \eqref{eq:self-duality} that it satisfies.
In our implementation the model is given by
\begin{equation}
u_0(z) = \frac{1}{N+2} \log \left ( |P_N(z)|^2 + \sigma(z) \exp(-|P_N(z)|^4) \right )
\end{equation}
where $\sigma$ is a smooth function that is positive on a disk $|z| < 0.9 r$ and vanishes elsewhere.\footnote{
We take $\sigma(z) = s(1-\abs{z}/(0.9r))$ where $s(t)$ is the smoothed step function $s(t) = 0$ for $t < 0$, $s(t) = 1$
for $t > 1$, and $s(t) = \frac12(1 - \cos (\pi t))$ for $t \in [0,1]$.}
(For $N=2$ this model should be compared to the function $u^{\rmsf}$ approximating $u$ that was discussed in \autoref{sec:semiflat}; indeed we have $u_0(z) = u^{\rmsf}(z)$ for all large $|z|$, and in general these functions are close except near the zeros of $P_2$ where $u_0$ is smooth while $u^{\rmsf}$ has logarithmic singularities.)

We then solve for $v$ on the grid with the Dirichlet condition $v = 0$ at the boundary, 
which is a reasonable approximation as $v$ is expected to decay exponentially in some power of $|z|$.
As \autoref{eq:self-duality} is nonlinear, we use Newton's method, i.e.~iteratively solving the linearization of the equation to improve an initial guess.
The iteration terminates when the size of the residual $\Delta u - 4 (e^{k u} - e^{-2u} |P_N|^2 )$ is less than a parameter \texttt{pde\_thresh}, in terms of a specific norm (which depends on the method, as described below).

For the core step of solving the linearization of \eqref{eq:self-duality} at a given point, our implementation offers two methods, named \texttt{euler} and \texttt{fourier}, which can be selected at runtime.
Both methods are rather elementary and were selected for simplicity of implementation.
Though most results in \autoref{sec:examples} use the \texttt{fourier} method, it is helpful to first explain the simpler \texttt{euler} method since \texttt{fourier} can be understood as a more complicated analogue of it with different trade-offs.

The \texttt{euler} solver uses the a finite-difference Laplacian based on the standard five-point stencil (see e.g.~\cite[Section~4.2]{thomas95}).
The linearized equation thus becomes a linear system of rank $(\texttt{pde\_nmesh})^2$ which is solved using \texttt{scipy.linalg.lstsq}.
In this method we measure the size of the residual in the Newton iteration using the $C^0$ norm.
This method stores several dense $(\texttt{pde\_nmesh})^2 \times (\texttt{pde\_nmesh})^2$ matrices and hence suffers from high memory consumption (approximately $12$ GiB for $\texttt{pde\_nmesh}=2000$).
Since we have found that increasing \texttt{pde\_nmesh} significantly improves the accuracy of these calculations, this limitation of the \texttt{euler} solver prompted development of a less memory-intensive alternative.

The \texttt{fourier} solver uses a simple spectral method, based on the $2D$ discrete Fourier transform, to avoid storing any dense $(\texttt{pde\_nmesh})^2 \times (\texttt{pde\_nmesh})^2$ matrices or solving any linear systems in the iterative step.
However, to make this possible, the \texttt{fourier} solver does not solve the linearization of \eqref{eq:self-duality} itself, which has the form $(\Delta - \kappa)u=f$ for a scalar function $\kappa$.
Rather, following Concus-Golub \cite{concus-golub} we replace the linearization with an approximating Helmholtz equation $(\Delta - C)v=f$ (for a constant $C$) which therefore has a closed-form solution in frequency space.
Here the constant $C$ is chosen to approximate the (non-constant) function $\kappa$ in the true linearization.
As in \cite{concus-golub} we use the ``minimax'' value $C = \frac{1}{2} \left (\sup \kappa + \inf \kappa \right )$ where both sup and inf are taken over the grid points.

To implement the desired Dirichlet boundary conditions on $v$ in the \texttt{fourier} solver we use the $2D$ discrete sine transform (DST), which is equivalent to extending $v$ as a doubly-periodic function which is odd with respect to reflections in the grid boundaries.
Also, in this method we measure the size of the residual in the Newton iteration using the $L^2$ norm, since this can be computed directly in frequency space, thus avoiding an additional Fourier transform step in the iteration.
The fact that the Helmholtz equation is a poor approximation of the true linearization has the effect of requiring many more iterations of Newton's method to reach a desired accuracy (in comparison to the \texttt{euler} solver).
However, in practice the high iteration count is more than compensated by the high speed of the Fast Fourier Transform when $\texttt{pde\_nmesh} > 1000$ and when the size has the optimal form for DST, i.e.~$\texttt{pde\_nmesh}=2^j-1$
for some positive integer $j$.

After solving the discretized self-duality equation, the result is a vector of values for $u$ at the grid points.
Since the next step of solving the ODE for parallel transport of the flat connection (\autoref{eq:hitchin-section-connection}) requires evaluation of the self-dual metric density $u$ at arbitrary points, the interpolation scheme from \texttt{scipy.RectBivariateSpline} is used with order $3$ in $x$ and $y$ to produce an approximation to the function $u$ on the bounding rectangle of the grid.
The same interpolation is applied to the finite difference approximations of the partial derivatives $u_x$ and $u_y$, which also appear in the connection form.
Finally, with a means of evaluating the connection form in hand, the process of solving the parallel transport ODE and computing Stokes data proceeds by the same process described in \autoref{subsec:implementation-opers}.

\subsection{PDE error estimate}
\label{subsec:error-hitchin}

We now explain the error estimates that are included in \hitchinfigs (results of calculations for the Hitchin section), which concern the discretization error introduced by solving an analogue of the self duality equation on a grid, using finite differences or discrete Fourier transforms instead of differential operators.
While in general this error is expected to be bounded by a multiple of $(\Delta x)^2$, where $\Delta x = 2 r / \texttt{pde\_nmesh}$ is the spacing of the grid points (along both the $x$- and $y$-axes), it would require a more subtle theoretical analysis of the method to derive the constant of proportionality.
Rather than conducting such a theoretical analysis, as mentioned in \autoref{subsec:results-hitchin} we derive an empirical error estimate using Richardson extrapolation.
The general theory of Richardson extrapolation is discussed in more detail in e.g.~\cite[Section 8.3]{oberkampf-roy}).

We recall the basic principle of this method:
First, we assume that the final quantity $X_i$ is subject to error that is approximately proportional to $(\Delta x)^p$ for some real $p$ (rather than being merely bounded by a quantity of this form).
Then, using results of calculations for three values of $\Delta x$ (equivalently, of $\texttt{pde\_nmesh}$), it is possible to both recover the value of $p$ that best fits the observed results, and to extrapolate to obtain a refined estimate $X_i^{\text{Rich}}$ for the limit as $\Delta x \to 0$.
In presenting our results we do not use this extrapolated value for the spectral coordinates directly, and instead we show the value computed with the largest $\texttt{pde\_nmesh}$.
However, the difference $|X_i^{DE} - X_i^{\text{Rich}}|$ is taken as the approximation of the discretization error in $X_i^{DE}$.
In addition, this method of empirical error estimation also allows a test for a fit between the error model and the observed results; we expect the best fit exponent $p$ to be approximately $2$, and significantly different values indicate that the hypothesis of discretization error being dominant and proportional to $(\Delta x)^2$ is not consistent with the results.
For this reason, PDE error estimates are only shown in \hitchinfigs in cases where the best fit exponent lies in the interval $[1.6,2.4]$.

\subsection{Integral equation method for opers}
\label{subsec:implementation-ieq-opers}

\begin{table}
\begin{tabular}{@{}lll@{}}
    \toprule
    Name & Value & Description\\
    \midrule
    \texttt{L} & $200$ & Interval size: $t = \log \abs{\hbar}$ runs over $[-L,L]$ \\
    \texttt{steps} & $2^{17}$ & Number of sampling points \\
    \texttt{tolerance} & $2 \times 10^{-15}$ & Target $L^\infty$ norm of difference between iterations \\
    \texttt{method} & \texttt{fourier} & Method for numerical integration ({\tt fourier} or {\tt simps}) \\
    \texttt{damping} & 0.3 & Damping factor in the iteration step \\
    \bottomrule
\end{tabular}

\vspace{1.5em}

\caption{Parameters used for the integral equation method for opers and the Hitchin section in the calculations presented in \autoref{subsec:results-opers} and \autoref{subsec:results-hitchin}.}
\label{tab:ieq-oper-params}
\end{table}

Here we describe some implementation details of our 
integral equation computation of spectral coordinates for
the opers $\nabla_w^\oper$.
In our Python implementation, these calculations are performed by the \texttt{integralequations} module.

As explained in \autoref{sec:ieq-opers} the problem is to find a solution $\cX_\bullet$ 
of the system of integral equations \eqref{eq:integral-oper}.
Since \eqref{eq:integral-oper} represents the desired functions 
$\cX_\gamma$ as exponentials, it is convenient
to write $\cX_\gamma = \exp x_\gamma$ and study $x_\gamma$
directly.

We begin with the initial guess
\begin{equation} \label{eq:initial-guess}
  x_\gamma^{(0)}(\hbar) = \hbar^{-1} Z_\gamma
\end{equation}
and then, defining
\begin{equation} \label{eq:F-concrete-1}
  (\cF(x))_\gamma(\hbar) = \hbar^{-1} Z_\gamma + \frac{1}{4 \pi \I} \sum_{\mu \in \Gamma'} \Omega(\mu) \langle \gamma, \mu \rangle \int_{\R_- Z_\mu} \frac{\de \xi}{\xi} \frac{\xi + \hbar}{\xi - \hbar} \log\left (1 + \exp (x_{\mu}(\xi)) \right),
\end{equation}
our iteration step is
\begin{equation} \label{eq:iteration-step}
  x_\gamma^{(n+1)} = (1-p) \cF(x^{(n)})_\gamma + p x_\gamma^{(n)}
\end{equation}
where $p = {\tt damping}$ is a damping parameter. Note that for any $p \in [0,1)$ the
fixed points of the iteration are the same as the fixed points of $\cF$; nevertheless
the rate of convergence can depend on $p$.
Importantly, in the iteration process we do not need to use the values
of $x_\gamma(\hbar)$ for arbitrary $\gamma$ and $\hbar$: all we need is
the values of $x_\gamma(\hbar)$ when $\gamma$ lies in the finite set $\Gamma'$
and $\hbar$ lies on the ray $\R_- Z_\gamma$.

We approximate the iteration numerically as follows.
For each charge $\gamma \in \Gamma'$ we parameterize the ray
$\hbar \in \R_- Z_\gamma$ by a parameter $t \in \R$, 
related to $\hbar$ by $\hbar = - \e^{\I \arg(Z_\gamma) + t}$.
We work with discrete approximations to functions $x^{(n)}_\gamma(t)$,
sampled at $M = {\tt steps}$ evenly spaced points $t \in [-L,L]$, where $L = {\tt L}$ 
is a cutoff parameter. We first construct $x_\gamma^{(0)}(t)$ by sampling the
function $\hbar^{-1} Z_\gamma = - \abs{Z_\gamma} \e^{-t}$.
Then, to construct $x_\gamma^{(n+1)}$, 
we numerically evaluate the right side of \eqref{eq:iteration-step} at each sampling point $t$.

To do the numerical integrations by Simpson's rule requires
$\sim M^2$ work at each iteration, because we need to evaluate $\sim M$ different integrals
each of which involves summing over $\sim M$ sampling points.
We reduce this work
to $\sim M \log M$ as follows.
First, the map $\cF$ preserves the property $x_\gamma(t) = x_{-\gamma}(t)$
(to see this we use the fact that $\Omega(\gamma) = \Omega(-\gamma)$).
Using this symmetry we can reduce our work in several ways.
First, we only need to compute $x_\gamma$ for half of the $\gamma \in \Gamma'$,
say all the ones where $Z_\gamma$ lies in some chosen half-plane.
Second, using the symmetry to average the term for $\mu$ with the term
for $-\mu$, we can rewrite \eqref{eq:F-concrete-1} as
\begin{equation} \label{eq:F-concrete-2}
  \cF(x)_\gamma(t) = - \abs{Z_\gamma} \e^{-t} + \frac{1}{4 \pi \I} \sum_{\mu \in \Gamma'} \Omega(\mu) \IP{\mu,\gamma} \int_{-\infty}^\infty \de t' \frac{\log(1 + \exp x_\mu(t'))}{\sinh((t-t') + \I \delta_{\gamma \mu})}, \quad \delta_{\gamma \mu} = \arg Z_\gamma - \arg Z_\mu,
\end{equation}
and the integral on the RHS is a convolution $F \ast G$ where
$F(t) = \log(1 + \exp x_\gamma(t))$ and $G(t) = \frac{1}{\sinh(t + \I \delta_{\gamma \mu})}$.
Such a convolution can be evaluated approximately, at all sampling points $t$ at once,
by the device
of transforming $F$ and $G$ to Fourier space (with a periodic boundary condition
at the ends of the interval $[-L,L]$), multiplying them, and then
transforming back. Using the Fast Fourier Transform this takes work
$\sim M \log M$.
(Our implementation offers this method as well as the slower, more direct method 
by Simpson's rule.)

We repeat the iterative process until the maximum value of
$\abs{x^{(n+1)}_\gamma (t) - x^{(n)}_\gamma (t)}$ over all sampling points
and all $\gamma \in \Gamma'$ is smaller than the constant {\tt tolerance}
for $5$ consecutive iterations, and then stop.
The final $x_\gamma(t)$ so obtained give a discrete approximation
to a solution of the integral equation \eqref{eq:integral-oper}, with each $x_\gamma$
computed along the ray $\hbar \in \R_- Z_\gamma$.
Finally, in order to evaluate the desired $\cX_\gamma(\hbar)$
for more general $\hbar$ (as we must do in order to test the conjecture), 
we just numerically evaluate the
RHS of \eqref{eq:integral-oper} once more
(this time using Simpson's rule instead of Fourier transforms, 
since we only need to evaluate \ti{one} integral
instead of $M$ of them.)

\subsection{Integral equation method for the Hitchin section}
\label{subsec:implementation-ieq-hitchin}

The methods we use for the integral equation computation for the Hitchin section are
essentially identical to those we use for opers (and are implemented in the same module, \texttt{integralequations}).
The same parameter values (from \autoref{tab:ieq-oper-params}) are used as well.
There are just a few slight changes in the formulas, since we now need to solve \eqref{eq:integral-hitchin-section} rather than \eqref{eq:integral-oper}.
The initial guess \eqref{eq:initial-guess} is replaced by
\begin{equation} \label{eq:initial-guess-h}
  x_\gamma^{(0)}(\zeta) = R \zeta^{-1} Z_\gamma + R \zeta \overline{Z}_\gamma,
\end{equation}
\eqref{eq:F-concrete-1} is replaced by
\begin{equation}
  (\cF(x))_\gamma(\zeta) = R \zeta^{-1} Z_\gamma + R \zeta \overline{Z}_\gamma + \frac{1}{4 \pi \I} \sum_{\mu \in \Gamma'} \Omega(\mu) \langle \gamma, \mu \rangle \int_{\R_- Z_\mu} \frac{\de \xi}{\xi} \frac{\xi + \zeta}{\xi - \zeta} \log\left (1 + \exp (x_{\mu}(\xi)) \right),
\end{equation}
and \eqref{eq:F-concrete-2} is replaced by
\begin{equation}
  \cF(x)_\gamma(t) = - 2 R \abs{Z_\gamma} \cosh t + \frac{1}{4 \pi \I} \sum_{\mu \in \Gamma'} \Omega(\mu) \IP{\mu,\gamma} \int_{-\infty}^\infty \de t' \frac{\log(1 + \exp x_\mu(t'))}{\sinh((t-t') + \I \delta_{\gamma \mu})}.
\end{equation}

\subsection{Direct method for the \hk metric}
\label{subsec:implementation-hk-direct}

\begin{table}
\begin{tabular}{@{}lll@{}}
    \toprule
    Name & Value & Description\\
    \midrule
    \texttt{method} & \texttt{euler} & PDE solver strategy (\texttt{euler} or \texttt{fourier})\\
    \texttt{rmax} & $10$ & Size of region for PDE solver \\
    \texttt{pde\_nmesh} & $1400$ & PDE mesh size \\
    \texttt{pde\_thresh} & $5 \times 10^{-11}$ & Absolute error goal for PDE solver\\
    \bottomrule
\end{tabular}

\vspace{1.5em}

\caption{Parameters used for the direct method for the \hk metric in the calculations presented in \autoref{sec:hk-experiments}.}
\label{tab:hk-direct-params}
\end{table}

The direct PDE computations of the \hk metric, as reported in \autoref{sec:hk-direct}, are handled by the module \texttt{hkmetric} (in particular, the function \texttt{fdcomputeG}).
Here, we first need to solve the PDE \eqref{eq:self-duality} for 
(a discrete approximation of) $u$,
then solve the linear PDE \eqref{eq:diff-F} for 
(a discrete approximation of) $F$.
The computation of $u$ 
is done using the same routines described in \autoref{subsec:implementation-hitchin}; this time we use
the {\tt euler} solver instead of {\tt fourier}.
The computation of $F$ uses the same code,
but runs for only one iteration since the equation is linear.
The parameters we use are listed in \autoref{tab:hk-direct-params}.

Once $F$ and $u$ have been computed, the last remaining step 
is to evaluate the integral \eqref{eq:l2norm-combined}.
For this we divide the plane into the regions $\abs{z} < r$ and $\abs{z} > r$,
where $r = {\tt rmax}$. In the region $\abs{z} < r$ we just approximate
the integral by a Riemann sum, using our ${\tt pde\_nmesh} \times {\tt pde\_nmesh}$ 
grid of sampling points; call this sum $I_{in}$. In the region $\abs{z} > r$
we do not have a numerical solution of the PDEs available, so we have
to make do with the asymptotic formulas 
$u \sim \frac12 \log \abs{P_2}$ and $F \sim \frac12 \frac{\dot P_2}{P_2}$;
when $P_2$ has leading term $z^n$ and $\dot{P}_2 = 1$, this gives
for the integrand $\cI \sim 2 \abs{z}^{-n+1}$. 
Using this we estimate the integral over the region $\abs{z} > r$ to be 
$I_{out} = \frac{4 \pi}{(n-2) r^{n-2}}$. (In \autoref{sec:hk-direct} we
have $n=3$, so $I_{out} = \frac{4\pi}{r}$.)
Our final result for the integral is then
\begin{equation}
	I = I_{in} + I_{out}.
\end{equation}

\subsection{Integral equation method for the \hk metric}
\label{subsec:implementation-ieq-hk}

Here we describe some details of the integral equation
computation of the \hk metric reported in 
\autoref{sec:hk-ieq}.
These calculations are performed by the module \texttt{hkmetric} (in particular, the function \texttt{ieqcomputeG}).

Using \eqref{eq:gc-concrete} we see that what we need is to to compute
the first derivatives of the quantities $y_i = \log \cX_i(R=1, \zeta=1)$ with
respect to the real and imaginary parts of the parameter $c$, at various
values of $c$.

The first step is to use the integral equation method 
(as described in \autoref{subsec:implementation-ieq-hitchin}) to compute $y_i$ itself
at various values of $c$. 
For this purpose
we need to know the periods $Z_\gamma$ as functions of $c$; 
fortunately, in this particular example
the periods are of the simple form $Z_\gamma = \mathrm{const} \times c^{\frac56}$.
Next we approximate the desired derivatives by the method of finite differences, choosing
a small $\eps = {\tt eps}$ and computing $\eps^{-1} (y_i(c+\eps) - y_i(c))$
and $\eps^{-1} (y_i(c+\I\eps) - y_i(c))$.
Finally we plug into \eqref{eq:gc-concrete} to get our estimate of $g(c)$.

The parameters we use are listed in \autoref{tab:ieq-hk-params}.

\begin{table}
\begin{tabular}{@{}lll@{}}
    \toprule
    Name & Value & Description\\
    \midrule
    \texttt{L} & $200$ & Interval size: $t = \log \abs{\zeta}$ runs over $[-L,L]$ \\
    \texttt{steps} & $2^{17}$ & Number of sampling points \\
    \texttt{tolerance} & $10^{-15}$ & Target $L^\infty$ norm of difference between iterations \\
    \texttt{method} & \texttt{fourier} & Method for numerical integration ({\tt fourier} or {\tt simps}) \\
    \texttt{damping} & 0.3 & Damping factor in the iteration step \\
    \texttt{eps} & $10^{-6}$ & Step size in finite difference estimate of derivatives \\ 
    \bottomrule
\end{tabular}

\vspace{1.5em}

\caption{Parameters used for the integral equation computation of the \hk metric in the calculations presented in \autoref{sec:hk-experiments}.}
\label{tab:ieq-hk-params}
\end{table}

\section{Sample numerical calculation}
\label{sec:sample}

Here is a sample session using our code to compute the cross-ratio invariants quoted in the introduction, and which can be run on a laptop:
{ \small
\begin{verbatim}
$ ipython3 harmonic.py -i
\end{verbatim}
[ ... ]
\begin{verbatim}
In [1]: xar = integralequations.computeXar(theoryname="A1A2")
\end{verbatim}
[ ... ]
\begin{verbatim}
[MainProcess] [integralequations] [INFO] Finished in 13.5 s

In [2]: xar.getCluster()
Out[2]: [-0.006415703123337184, -1.0]

In [3]: fd = framedata.computeFrames(theoryname="A1A2", pde_nmesh = 1023)
\end{verbatim}
[ ... ]
\begin{verbatim}
[MainProcess] [framedata] [INFO] Finished in 192.6 s

In [4]: fd.getCluster()
Out[4]: [-0.006415963850395493, -0.9999999999999987]
\end{verbatim}
}
The result from the PDE calculation (the last line of the output above) differs a bit from the result quoted in the
introduction; the result in the introduction is obtained if one sets {\tt pde\_nmesh = 8191},
which however requires more RAM than is available on a typical laptop in 2020.

A note about signs: the code uses a different convention in defining
the $\cX_\gamma$ than this paper does. For charges $\gamma \in \Gamma'$
the relation is simply
\begin{equation}
\cX_\gamma^{\mathrm{code}} = - \cX_\gamma^{\mathrm{paper}}.
\end{equation}
For more general charges $\gamma$ the relation
is
$\cX_\gamma^{\mathrm{code}} = \sigma(\gamma) \cX_\gamma^{\mathrm{paper}}$
where $\sigma(\gamma) = \pm 1$ is a certain quadratic refinement of the 
pairing $(-1)^{\IP{,}}$.

\section{Discussion}
\label{sec:discussion}

\subsection{Interpretation of results}

The quantities $X_i^{\DE}$ and $X_i^{\IEQ}$ reported in \autoref{sec:experiments} are numerical approximations of two quantities that the main conjecture asserts to be equal.
Therefore, assuming the conjecture, the difference between $X_i^{\DE}$ and $X_i^{\IEQ}$ would reflect only the numerical error in the two computations.
We have analyzed one source of error in $X_i^{\DE}$ for each class of experiments (ODE error for opers, PDE discretization error for Hitchin section) and obtained an estimate that in practice grows rapidly with the relevant parameter ($R$ or $\abs{\hbar}^{-1}$).
Broadly, our experiments show that the relative difference of $X_i^{\DE}$ and $X_i^{\IEQ}$ is either small (for small parameter values), or else is comparable to the single-source error estimate (for large parameter values).
This is consistent with the hypothesis that the conjecture holds, and that the source of error considered in our estimates is the dominant one for large parameter values, while other sources of error become dominant for small parameters.
Notably absent from our results are any strong candidates for counterexamples to the conjecture, such as examples where a large relative difference between $X_i^{\DE}$ and $X_i^{\IEQ}$ is both stable under variation of the parameter and greatly exceeds the single-source error estimate.

In our experiments computing spectral coordinates for opers, the relative difference between $X_i^{\DE}$ and $X_i^{\IEQ}$ is less than the ODE error estimate whenever $\abs{\hbar}^{-1} > 2$, and is less than $10^{-9}$ (1 PPB) for $\abs{\hbar}^{-1} \leq 2$.
For the Hitchin section, an additional complication is that the error estimate is not always available; recall that the empirical PDE error estimate is applied to each spectral coordinate separately, and it requires the results for varying mesh sizes to match a theoretical model.
Nevertheless the results exhibit similar behavior to the case of opers, though with relative error that is roughly $1000$ times larger:
For $R$ less than an example-dependent threshold (approximately $0.5$ for $(A_2,A_2)$ or $1$ for the other examples) the relative difference remains below $10^{-6}$ (1 PPM), and the error estimate often does not apply in this region.
This is consistent with the hypothesis that other sources of error dominate here, preventing the Richardson extrapolation from seeing the expected dependence on the mesh size.
For $R$ above the threshold, the error estimate succeeds in most cases, and gives a result that closely tracks (and often exceeds) the observed difference.
The $(A_1,A_3)$ experiment is the least consistent with this overall pattern, as the error estimate fails for a significant number of large-$R$ data points.
In general the $(A_1,A_3)$ Hitchin section computations were also the most computationally expensive for the direct method, requiring significantly more iterations than other examples at comparable parameter values (e.g.~those which ultimately give similar relative difference with $X_i^{\IEQ}$), and we expect this is related to this example having the highest-degree polynomial differential among all those considered.

Our experiments with the \hk metric for the $(A_1,A_2)$ case show that using the integral equation method to calculate this metric through \autoref{eq:metric-spectral} has reasonably close agreement with the results of a direct calculation of the metric by solving the relevant PDE.
The agreement between the two \hk metric calculation methods is not as close (in terms of relative difference) as that of the spectral coordinates themselves, however this is to be expected.
Both methods of calculating the \hk metric involve an additional layer of approximation that builds on the spectral coordinate computations:
For the direct method, the numerical approximation of the self-dual metric $u$ is used as input to another PDE solver which computes the complex variation $F$, allowing error in $u$ itself to propagate to a larger error in $F$.
In the integral equation method, we take finite differences of the spectral coordinates, amplifying the error in the individual coordinate calculations and introducing additional truncation error.

\subsection{Limitations}

Our experimental approach initially prioritized simplicity of implementation, and was incrementally improved to explore additional aspects of the twistor Riemann-Hilbert conjecture.
If a substantial revision or full redesign were considered, it would be natural to attempt to address the following limitations of our approach.

\smallskip

\begin{itemize}

\item \textbf{Numerical instability of the direct method.}
This significant limitation of our approach was already mentioned in \autoref{subsec:implementation-opers}:
As $\abs{\hbar}^{-1} \to \infty$ or $R \to \infty$, a typical spectral coordinate approaches $0$ or $\infty$, and hence one of the determinants $p(i,j)$ or $p(i,j,k)$ involved in its calculation must approach $0$.
For large parameters we are therefore computing the determinant of a nearly-singular matrix whose entries are of size $\sim 1$ and which are known to fixed relative precision, resulting in a significant loss of precision.
This phenomenon is the reason our current implementation is unable to test the twistor Riemann-Hilbert conjecture for larger parameter values.

\item \textbf{High memory usage of the direct method (Hitchin section).}
While the \texttt{fourier} method requires much less memory than the \texttt{euler} method, memory was still the main factor constraining the values of $\texttt{pde\_nmesh}$ we were able to study.
This is significant since discretization error is expected to be the dominant source of numerical error in the direct method calculations for large $R$ or $\abs{\hbar}^{-1}$, and in our current approach, discretization error is roughly proportional to $(\texttt{pde\_nmesh})^{-2}$.
However, rather than attempting to reduce the memory requirements of the current method, we believe a better approach would be to use a PDE solver that can achieve higher accuracy at a given $\texttt{pde\_nmesh}$.

\item \textbf{Lack of error bounds.}
As noted earlier, our study of the numerical error in our calculations is limited:
We do not consider the error in the integral equation method at all, and our analysis of error in the direct method is an estimate, not a rigorous bound.
While working out precise error bounds for these methods would be a substantial undertaking, it would be necessary in order to make rigorous positive statements about the bearing of the results on the main conjecture, e.g.~to say that the predictions of the conjecture in a given case hold to within numerical error.

\end{itemize}

\subsection{Future directions}

Our work suggests several directions that could be pursued in future experiments, though each brings its own challenges.

\smallskip

\begin{itemize}

\item \textbf{Changing the base Riemann surface.}
All of our experiments involve Higgs bundles or opers on the complex plane.
The use of a single global coordinate chart and the expression of the holomorphic differentials as polynomials are deeply embedded in our implementation, and so a change of this type would require substantial revision.
The cylinder $\C^*$ is the case accessible by the least modification of our code; there it would be possible to work on the universal cover, imposing suitable mixed boundary conditions on a rectangular region projecting to an annulus $\eps \leq |z| \leq 1/\eps$.

\item \textbf{Changing the rank.}
We have considered Higgs bundles of rank $N \leq 3$, because these are exactly the cases where the self-duality equation reduces to a single scalar equation.
It would be relatively straightforward to generalize our study of opers to higher rank, but it would require major changes to study the Hitchin section in higher rank due to the need to handle systems of PDE.

\item \textbf{Moving off the Hitchin section.}
It would also require considerable changes to our implementation to study Higgs bundles that are not in the Hitchin section (even for rank $N=2$).
However, generalizing the code in this direction would be especially appealing; the full \hk metric on $\cM_{d,h}$
is a considerably richer object than the \kahler metric on $\cB_{d,h}$, and while the twistor Riemann-Hilbert conjecture
does extend to the whole $\cM_{d,h}$, we are not aware of any numerical studies in this direction.
The asymptotics of the full metric on $\cM_{d,h}$ have been the subject of extensive study recently
(see e.g.~\cite{Mazzeo2019,Fredrickson2020} and the survey \cite{Fredrickson2019}); the results support the
twistor Riemann-Hilbert conjecture.

\item \textbf{Integral equation prediction for self-dual metric.}
While we have focused on the conjecture's predictions for spectral coordinates (and the \hk metric, as computed from derivatives of spectral coordinates), in the case of Higgs bundles it is also possible to refine the twistor Riemann-Hilbert conjecture to give a conjectural formula for the self-dual metric itself in terms of integral equations described in \cite{Gaiotto:2011tf}.
If this extension of the integral equation approach were implemented, it would be very convenient to compare the results to the direct methods already implemented in our code.
Indeed, such a comparison would avoid all of the complexity of computing parallel transports and determinantal invariants, and quite possibly be easier to analyze rigorously to obtain bounds on numerical error.

\end{itemize}

\medskip

We hope to pursue some of these directions in the future.


\vfill
\pagebreak

\bibliographystyle{hyperplain}
\bibliography{numerics}

\vspace{1.5em}

\noindent Department of Mathematics, Statistics, and Computer Science\\
University of Illinois at Chicago\\
Chicago, IL\\
\texttt{emily@dumas.io}

\vspace{1em}

\noindent Department of Mathematics\\
Yale University\\
New Haven, CT\\
\texttt{andrew.neitzke@yale.edu}

\end{document}